\documentclass{amsart}
\usepackage{amssymb}
\newtheorem{Thm}{Theorem}[section]
\newtheorem{Prop}[Thm]{Proposition}
\newtheorem{Lem}[Thm]{Lemma}
\newtheorem{Cor}[Thm]{Corollary}
\newtheorem{Def}[Thm]{Definition}
\newtheorem{Rem}[Thm]{Remark}

\numberwithin{equation}{section}

\begin{document}

\title{On BC Type Basic Hypergeometric Orthogonal Polynomials}

\author{J.V. Stokman}
\address{Department of Mathematics, University of Amsterdam,
Plantage Muidergracht 24, 1018 TV  Amsterdam, The Netherlands. Email:
jasper@wins.uva.nl}


\date{July, 1997}


\thanks{The author received financial support by NWO/Nissan}

\begin{abstract}
The five parameter family of multivariable 
Askey-Wilson polynomials is studied with four parameters
generically complex. The multivariable Askey-Wilson 
polynomials form an orthogonal system with respect to 
an explicit (in general complex) measure.
A partially discrete orthogonality measure is obtained 
by shifting the contour to the torus while picking up residues. A parameter
domain is given 
for which the partially discrete orthogonality measure is positive.  
The orthogonality relations and norm evaluations for multivariable
$q$-Racah polynomials and multivariable big and little $q$-Jacobi polynomials
are proved by taking suitable limits in the orthogonality
relations for the multivariable Askey-Wilson polynomials.
In particular new proofs of several well known $q$-analogues of the Selberg
integral are obtained.

\end{abstract}

\maketitle

\section{Introduction}

In \cite{M1} Macdonald introduced a remarkable family of 
multivariable ($q$-)orthogonal polynomials associated with root systems 
(the so called Macdonald polynomials). The Macdonald polynomials interpolate
various families of special functions associated with groups, such as
spherical functions on real semisimple Lie groups ($q=1$),
spherical functions on semisimple p-adic groups ($q=0$) and 
characters of compact semisimple Lie groups. 
They also arise as  spherical functions on certain quantum symmetric spaces  
(see for instance \cite{N}, \cite{NS}). 

Recent work by Cherednik (see for instance
\cite{C}, \cite{C0}, \cite{C1}) has shown that these
polynomials are closely related with certain representations of 
affine Hecke algebras. 
Many properties of the polynomials can be derived by this
approach while they were first untractable or could 
only be proved for very special parameter values
from the (quantum) group interpretation.
In particular the explicit evaluations for the quadratic norms of the Macdonald
polynomials, which were conjectured by Macdonald \cite{M1}, have been 
proved by Cherednik \cite{C} using this approach. 

In this paper 
we will consider families of multivariable $q$-orthogonal polynomials 
associated with the non-reduced root system $BC$. 
Koornwinder \cite{K1} 
extended the three parameter families of BC type Macdonald
polynomials to a five parameter multivariable generalization of
the well-known four parameter family of Askey-Wilson polynomials \cite{AW}.
Van Diejen \cite{vD1} evaluated the quadratic norms
of the multivariable Askey-Wilson polynomials using so called Pieri formulas.

The one-variable Askey-Wilson polynomials contain
various interesting families of basic hypergeometric orthogonal polynomials
as special cases or as limit cases (these families are collected in the 
Askey-Wilson scheme, see \cite{AW}, \cite{KS}). 
In particular the Askey-Wilson scheme contains the families
of $q$-Racah polynomials \cite{AW0}, big $q$-Jacobi polynomials \cite{AA2}
and little $q$-Jacobi polynomials \cite{AA1}. These three families were
recently also studied in the multivariable setting.
 
The multivariable $q$-Racah polynomials \cite{vDS} can be considered
as multivariable Askey-Wilson polynomials for which the parameters satisfy a  
particular truncation condition.
In \cite{vDS} it was shown that they are orthogonal 
with respect to a finite, discrete orthogonality measure. 
The multivariable big and little $q$-Jacobi polynomials \cite{S1}
can be considered as 
limit cases of rescaled  multivariable Askey-Wilson polynomials 
in which some of the parameters tend to infinity (see \cite{SK1}). 
In \cite{S1} it was shown that they are orthogonal with
respect to infinite 
discrete measures (which can most conveniently be expressed in
terms of multidimensional Jackson integrals). 
The orthogonality measure for the three
limit cases was obtained by repeating the techniques 
used by Macdonald for the Macdonald polynomials 
(and similarly by Koornwinder for the multivariable Askey-Wilson polynomials).
The quadratic norm evaluations for the multivariable $q$-Racah polynomials 
were obtained in 
\cite{vDS} by use of Pieri formulas (see \cite{vD2} for the constant term
identity, i.e. the quadratic norm evaluation of the unit polynomial).

In this paper the orthogonality relations and norm evaluations
for the multivariable 
$q$-Racah polynomials, big and little $q$-Jacobi polynomials
are proved by taking suitable limits in the orthogonality relations
and norm evaluations for the Askey-Wilson polynomials. We proceed as follows. 
In section 2 the orthogonality relations 
for the Askey-Wilson polynomials obtained by Koornwinder \cite{K1} 
and the quadratic norm evaluations 
obtained by van Diejen \cite{vD1} are extended 
to the case where four of the five parameters are generically complex. 
The extended orthogonality measure is an absolutely continuous
complex measure with weight function identical to the weight function
considered by Koornwinder \cite{K1}, 
but with a suitably deformed integration contour. 
If all the four parameters have moduli $<1$ then 
the $n$-torus $T^n$ may be chosen as integration contour
and Koornwinder's \cite{K1} orthogonality measure is recovered.

In section 3 a residue calculus is developed for the complex 
orthogonality measure of section 2. The orthogonality
relations and norm evaluations for the multivariable Askey-Wilson polynomials
can be reformulated 
with respect to partially discrete orthogonality measures by use of this calculus.

In section 4 the limit from multivariable Askey-Wilson polynomials 
to multivariable $q$-Racah polynomials is studied 
on the level of the orthogonality measures.
A suitable partially discrete orthogonality measure for the 
multivariable Askey-Wilson polynomials is rescaled such that 
certain common poles of the completely discrete weights 
become zeros for the continuous parts of the orthogonality measure. 
These zeros cause the vanishing of the continuous parts of the
orthogonality measure in the limit. We end up with a finite discrete 
measure, which is easily recognized as the measure 
of the multivariable $q$-Racah polynomials. We obtain new proofs of the 
orthogonality relations and norm evaluations for the multivariable 
$q$-Racah polynomials, simply by taking the limit in the 
corresponding results for the multivariable Askey-Wilson polynomials. 
See \cite{SK2} for the use of this idea in the one-variable case.

In section 5 the residue calculus of section 3 is used to deform the 
integration contour of the orthogonality measure for the
multivariable Askey-Wilson polynomials (cf. section 2) 
to the $n$-torus $T^n$. We obtain a partially discrete 
orthogonality measure which is positive for a large parameter domain. 
In particular the parameter values which occur in   
the limits from multivariable Askey-Wilson polynomials to multivariable
big and little $q$-Jacobi polynomials lie in this parameter domain.

In sections 6 and 7 these two limits are taken in  
the (suitably rescaled) positive, partially discrete orthogonality measure 
for the Askey-Wilson polynomials. The partially continuous contributions
(which are supported on subtori of $T^n$) disappear in these limits
because the corresponding weight functions tend to zero, 
while the support of the completely discrete
part of the measure blows up to an infinite set.
We end up with the orthogonality measures of the multivariable little and
big $q$-Jacobi polynomials. The orthogonality relations 
and norm evaluations for the multivariable little  
and big $q$-Jacobi polynomials are obtained by taking limits 
of the corresponding results for the multivariable Askey-Wilson
polynomials. 
The rigorous proofs of the limits of the orthogonality measures
are postponed to section 8 and 9.

The constant term identities which are obtained as special cases of the 
quadratic norm evaluations for the multivariable little and 
big $q$-Jacobi polynomials 
reduce to well known $q$-analogues of the Selberg integral. 
The constant term identity 
for the multivariable little $q$-Jacobi polynomials (Corollary \ref{gevolg}) 
is known as the Askey-Habsieger-Kadell formula (see \cite{A},\cite{H},\cite{K})
and was proved in full generality by Aomoto \cite{Ao}. 
The constant term identity for the multivariable big $q$-Jacobi polynomials 
with one of the parameters discrete (Corollary \ref{AE}) 
was conjectured by Askey \cite{A} and proved by Evans \cite{E}. 
The constant term identity in the general form  (Corollary \ref{gevolgBB}) 
is equivalent to Tarasov's and Varchenko's constant term identity \cite[Theorem
(E.10)]{TV}. 

The method of proving the orthogonality relations 
and norm evaluations for the one-variable $q$-Racah, big and little $q$-Jacobi 
polynomials by considering them as limit cases of the 
Askey-Wilson polynomials was discussed in detail
in \cite{SK2}. Although the computations are much more
involved in the multivariable setting, the techniques we employ 
here are essentially the same as in the one-variable case (cf. \cite{SK2}).

Finally let me note that the starting points in this article are 
the orthogonality relations and the Pieri formulas 
for the five parameter
family of multivariable Askey-Wilson polynomials, together
with an explicit second order $q$-difference operator which diagonalizes the 
multivariable Askey-Wilson polynomials (see the proof of Theorem \ref{concl}). 
A complete proof of the Pieri formulas is at present only published for a four 
parameter subfamily \cite{vD1}. On the other hand, 
the affine Hecke-algebraic approach to
the study of orthogonal polynomials 
related to root systems as mentioned earlier
in the introduction can in fact be worked out for the complete five parameter
family of multivariable Askey-Wilson polynomials (this is announced in \cite{M}
and worked out in more detail in \cite{N1}). As van Diejen remarked in 
\cite[note added in proof]{vD1}, 
this approach turns out to imply all the results 
in \cite{vD1} for the complete five parameter family of 
multivariable Askey-Wilson polynomials.
In the proof of Theorem \ref{concl} the explicit 
Pieri formulas for the complete five parameter family of multivariable 
Askey-Wilson polynomials 
are therefore used without further addressing this matter.

{\it Notations and conventions:} We write
${\mathbb{N}}=\lbrace 1,2,\ldots\rbrace$,
${\mathbb{N}}_0={\mathbb{N}}\cup\lbrace 0 \rbrace$,
${\mathbb{C}}^*= {\mathbb{C}}\setminus \lbrace 0 \rbrace$
and ${\mathbb{R}}^*={\mathbb{R}}\setminus\lbrace 0 \rbrace$.
Empty sums are equal to $0$, empty products are equal to $1$. 
In order to keep the notations transparent, we will omit in formulas
the dependance on parameters if it is clear from the context.
In this paper $n$ denotes the rank of the BC type root system  
(i.e. the multivariable orthogonal polynomials which we study in this paper 
depend on $n$ variables $z=(z_1,\ldots,z_n)$).
We use the convention that a function $h(z)$ for which a definition
with respect to the $n$ variables $z=(z_1,\ldots,z_n)$ is given, 
should be read with $n=m$ if it follows from the context 
that $z\in {\mathbb{C}}^m$. If $h(z)$ appears in formulas with
$z\in {\mathbb{C}}^m$ and $m=0$, then $h(z)$ should be read as $1$.
Throughout this paper, we fix a $q\in (0,1)$.

\section{Askey-Wilson polynomials for generic parameters values}

We first introduce some notations.  
The q-shifted factorial is given by
\begin{equation}\label{qshift}
\left(a;q\right)_b := 
\frac{\left(a;q\right)_{\infty}}{\left(aq^b;q\right)_{\infty}},
\quad \left(a;q\right)_{\infty} := \prod_{j=0}^{\infty}\left(1-q^ja\right),
\end{equation}
provided that $aq^b\notin\lbrace q^{-k}\rbrace_{k\in{\mathbb{N}}_0}$. For
$b=k\in {\mathbb{N}}_0$, 
we have $\left(a;q\right)_k=\prod_{i=0}^{k-1}(1-aq^i)$,
which is well defined for all $a\in {\mathbb{C}}$.
We write furthermore
\begin{equation}\label{qshiftprod}
\left(a_1,\ldots,a_m;q\right)_b:=
\prod_{j=1}^m\left(a_j;q\right)_b
\end{equation}
for products of $q$-shifted factorials.

Let $S$ be the group of permutations of 
the set $\lbrace 1,\ldots,n\rbrace$ and $W=S\ltimes \lbrace \pm 1 \rbrace^n$
the Weyl group of type $BC_n$. 
Let $z_1,\ldots,z_n$ be independent
variables, then $W$ acts in a natural way 
on the algebra $A:={\bf C}[z_1^{\pm 1},\ldots,z_n^{\pm 1}]$. We denote
$A^W$ for the subalgebra
of $W$-invariant functions in $A$. A basis for $A^W$ is given by 
the monomials $\lbrace m_{\lambda}\, | \, \lambda\in \Lambda \rbrace$, where
$\Lambda:=\lbrace \lambda\in {\mathbb{N}}_0^n \, | \, \lambda_1\geq\ldots\geq
\lambda_n \rbrace$, 
and
\[m_{\lambda}(z):=\sum_{\mu\in W\lambda}z^{\mu}\]
with $z^{\mu}=z_1^{\mu_1}\ldots z_n^{\mu_n}$.
The $W$-orbit of $\lambda\in \Lambda\subset {\mathbb{Z}}^n$ is with 
respect to the 
natural action of $W$ on ${\mathbb{Z}}^n$.

We write $\underline{t}=(t_0,t_1,t_2,t_3)$ for the four tuple
of parameters $t_0,t_1,t_2,t_3$. We assume in this section that
$\underline{t}\in V$, where $V\subset {\mathbb{C}}^4$ is the following
parameter domain.
\begin{Def}\label{WAW}
Let $V$ be the set of parameters $\underline{t}\in ({\mathbb{C}}^*)^4$
for which
\[
\#\lbrace \hbox{arg}(t_i),\hbox{arg}(t_i^{-1}) \, | \, i=0,1,2,3\rbrace=8
\]
and for which $t_0t_1t_2t_3\notin {\mathbb{R}}_{\geq 1}$.
Here $\hbox{arg}(u)\in [0,2\pi )$ is the argument of $u\in {\mathbb{C}}^*$
and ${\mathbb{R}}_{\geq 1}:=\lbrace r\in {\mathbb{R}} \, | \, r\geq 1\rbrace$.
\end{Def}
For $\underline{t}\in V$ 
we define $\alpha^{\pm}_i=\alpha^{\pm}_i(\underline{t})$
by
\begin{equation}\label{alpha}
\alpha_i^{\pm}:=\frac{\hbox{arg}(t_i^{\pm 1})}{2\pi}, \qquad i\in\lbrace 
0,1,2,3\rbrace.
\end{equation}
Note that 
$\alpha_i^{\pm}\not=0,1/2$ and that $\alpha_i^-=1-\alpha_i^+$ for all $i$.

The (in general complex) measure which we will introduce in this
section  is supported on a suitably 
deformed $n$-torus $C^n\subset {\mathbb{C}}^n$, where 
$C\subset {\mathbb{C}}$ is the following deformation
of the unit circle $T$.
\begin{Def}\label{contour}
We call a continuous rectifiable
Jordan curve $C=\phi_{C}([0,1])\subset {\mathbb{C}}$ 
a deformed circle if $C$ has a  
parametrization $\phi_C$ of the form 
\begin{equation}\label{param}
\phi_{C}(x)=r_{C}(x)e^{2\pi i x}\,\, (x\in [0,1]), 
\quad r_C: [0,1]\rightarrow (0,\infty)
\end{equation}
and if $C$ is invariant under inversion, i.e.
$C^{-1}:=\lbrace z^{-1} \, | \, z\in C \rbrace=C$.
For $\underline{t}\in V$, we call a deformed circle $C$ a 
$\underline{t}$-contour if 
the four parameters $t_0,t_1,t_2,t_3$ are in the interior of $C$.
\end{Def}
For a deformed circle $C=\phi_C([0,1])$ the radial function $r_C$ satisfies 
$r_C(1-x)=(r_C(x))^{-1}$ since $C=C^{-1}$. Since a deformed circle $C$ is 
by definition a closed contour, we furthermore have that
$r_C(0)=r_C(1/2)=r_C(1)=1$.
Note that the unit circle $T$ is a deformed circle
with $r_T\equiv 1$.
If $\underline{t}\in V$ and $C$ is a $\underline{t}$-contour, then the radial
function $r_C$ satisfies the extra conditions
\[r_C(\alpha_i^+)>|t_i|,\qquad i\in\lbrace 0,\ldots,3\rbrace
\]
since $t_i$ is in the interior of $C$ for all $i$.
In particular the unit circle $T$ is a $\underline{t}$-contour
if $|t_i|<1$ for all $i\in\lbrace 0,\ldots,3\rbrace$.
We will use the convention that a deformed circle $C$ is counterclockwise
oriented (i.e. has the orientation induced from the parametrization $\phi_C$)
when we integrate over $C$.
 
Let $t\in (0,1)$, $\underline{t}\in V$ and let $C$ be a deformed
circle such that $t_iq^j\notin C$ for $i\in\lbrace 0,\ldots,3\rbrace$ and
$j\in {\mathbb{N}}_0$.
Let $d\nu(z;\underline{t};t)$ be the measure on $C^n$ given by
\begin{equation}\label{measure}
d\nu(z;\underline{t};t):=\Delta(z;\underline{t};t)\frac{dz}{z}
\end{equation}
for $z\in C^n$ with $\frac{dz}{z}:=\frac{dz_1}{z_1}\ldots\frac{dz_n}{z_n}$ and
with weight function $\Delta(z;\underline{t};t)$ given
by
\begin{equation}\label{opdeling+}
\Delta(z;\underline{t};t)=\left(\prod_{j=1}^nw_c(z_j;\underline{t})\right)
\delta(z;t),
\end{equation}
with $w_c(x;\underline{t})$ given by
\begin{equation}\label{wc}
w_c(x;\underline{t}):=\frac{\bigl(x^2,x^{-2};q\bigr)_{\infty}}
{\bigl(t_0x,t_0x^{-1},
t_1x,t_1x^{-1},t_2x,t_2x^{-1},t_3x,t_3x^{-1};q\bigr)_{\infty}}
\end{equation}
and $\delta(z;t)$ given by
\begin{equation}\label{deltatau}
\delta(z;t):=\prod_{1\leq i<j\leq
n}\bigl(z_iz_j,z_i^{-1}z_j,z_iz_j^{-1},z_i^{-1}z_j^{-1};q\bigr)_{\tau},
\quad t=q^{\tau}.
\end{equation}
The factor $w_c(x;\underline{t})$ is exactly the 
weight function 
which occurs in the continuous part of the orthogonality measure
for the one-variable Askey-Wilson polynomials \cite{AW}.
The interaction factor $\delta(z;t)$ is only present in the multivariable
setting (i.e. when $n>1$), so in particular the measure
is independent of the deformation parameter $t$ when $n=1$.

The measure $d\nu(.;\underline{t};t)$ on $C^n$ is well defined, since
the poles of $\Delta(.;\underline{t};t)$ do not lie on the integration
domain $C^n$. Indeed, 
the poles of the weight function $\Delta(z;\underline{t};t)$
lie on hyperplanes
\begin{equation}\label{poles1}
z_i=t_mq^j\,\, \hbox{ or } \,\, z_i=t_m^{-1}q^{-j}
\end{equation}
with $m\in \lbrace 0,\ldots,3\rbrace$,  $j\in {\mathbb{N}}_0$ and
$i\in\lbrace 1,\ldots,n\rbrace$ (the poles coming from $w_c(z_i)$) 
and on hypersurfaces
\begin{equation}\label{poles2}
z_k^{\epsilon_k}z_l^{\epsilon_l}=t^{-1}q^{-j}
\end{equation}
with $1\leq k\not=l\leq n$, $j\in {\mathbb{N}}_0$ and
$\epsilon_k,\epsilon_l\in \lbrace -1,1\rbrace$ (the poles coming from
$\delta(z;\underline{t})$). We have $z\notin C^n$ for a pole $z$
of the form \eqref{poles1} since $C^{-1}=C$ and the assumption that
$t_iq^j\notin C$ ($i\in\lbrace 0,1,2,3\rbrace, j\in {\mathbb{N}}_0$)
and it follows from the definition of a deformed circle $C$ that
$z\notin C^n$ for a pole $z$ of the form \eqref{poles2}. 

In this section we will study the orthogonal polynomials related to
the complex measure $\bigl(C^n,d\nu(.;\underline{t};t)\bigr)$ 
where $C$ is an arbitrary $\underline{t}$-contour.
We first show that the measure
$d\nu(.;\underline{t};t)$ is independent of the $\underline{t}$-contour $C$
when integrating against $W$-invariant Laurent polynomials. 
In order to obtain this result, 
we first define specific subsets of $({\mathbb{C}}^*)^n$
on which the interaction factor $\delta(.;t)$ is analytic.

Let $C$, ${\mathfrak{C}}$ be deformed circles,
with parametrization given by $\phi_C(x)=r_C(x)e^{2\pi ix}$ respectively
$\phi_{\mathfrak{C}}(x)=
r_{\mathfrak{C}}(x)e^{2\pi ix}$. Let $A^+(C,{\mathfrak{C}})$
be the open subset 
\begin{equation}\label{A+}
A^+(C,{\mathfrak{C}}):=
\lbrace x\in [0,1] \, | \, r_{C}(x)>r_{\mathfrak{C}}(x) \rbrace
\subset (0,1).
\end{equation}
Set 
\begin{equation}\label{Omega}
\Omega(C,{\mathfrak{C}}):=\Omega^+(C,{\mathfrak{C}})\cup
\Omega^+({\mathfrak{C}},C)\cup C,
\end{equation}
where $\Omega^+(C,{\mathfrak{C}})$ is given by
\begin{equation}
\Omega^+(C,{\mathfrak{C}}):=\bigcup_{x\in A^+(C,{\mathfrak{C}})}
\lbrace y(x)e^{2\pi ix} \, | \, r_C(x)\geq y(x)\geq r_{\mathfrak{C}}(x)\rbrace.
\end{equation}
The following properties of $\Omega(C,{\mathfrak{C}})\subset
{\mathbb{C}}^*$ follow easily from the definitions:\\
{\it (i)} $\Omega(C,{\mathfrak{C}})=\Omega({\mathfrak{C}},C)$.\\
We will use, 
in view of {\it (i)}, the notation $\Omega=\Omega(C,{\mathfrak{C}})$
when it is clear 
from the context which pair of contours $C$, ${\mathfrak{C}}$ is
meant.\\
{\it (ii)} $\Omega^{-1}=\Omega$.\\
{\it (iii)} The contour $C$ can de deformed homotopically to ${\mathfrak{C}}$
within $\Omega$.

We call $\Omega\subset
{\mathbb{C}}^*$ {\it the domain associated with the pair} $(C,{\mathfrak{C}})$.
We write ${\mathcal{O}}_W(\Omega^n)$ 
for the ring of $W$-invariant functions $f$
which are analytic on $\Omega^n$. We have now the following crucial lemma.

\begin{Lem}\label{deltaanalytic}
Let $t\in (0,1)$ and
let $C$, ${\mathfrak{C}}$ be deformed circles satisfying the condition  
$tr_C(x)<r_{\mathfrak{C}}(x)$ for all $x\in A^+(C,{\mathfrak{C}})$. 
Then $\delta(.;t)\in {\mathcal{O}}_W(\Omega^n)$.
\end{Lem}
\begin{proof}
Let $C$, ${\mathfrak{C}}$ 
be deformed circles satisfying $tr_C(x)<r_{\mathfrak{C}}(x)$
for all $x\in A^+(C,{\mathfrak{C}})$. Let $z\in ({\mathbb{C}}^*)^n$
such that $z_k^{\epsilon_k}z_l^{\epsilon_l}=t^{-1}q^{-j}$ for some
$j\in {\mathbb{N}}_0$, 
some $k\not=l$ and some $\epsilon_k, \epsilon_l\in\lbrace
\pm 1\rbrace$. Write $\beta_k:=\hbox{arg}(z_k)/2\pi$ and
$\beta_l:=\hbox{arg}(z_l)/2\pi$. For the proof of the lemma it suffices to 
show that either $z_k\not\in \Omega$ or $z_l\not\in\Omega$. 

As an example, 
let us check that either $z_k\not\in\Omega$ or $z_l\not\in\Omega$
when $\beta_k\in A^+(C,{\mathfrak{C}})$ 
and $z_kz_l=t^{-1}q^{-j}$ for some $j\in
{\mathbb{N}}_0$. We then have $\beta_k=1-\beta_l$ and $\beta_l\in
A^+({\mathfrak{C}},C)$, so in particular
$r_C(\beta_k)=r_C(\beta_l)^{-1}$,
$r_{\mathfrak{C}}(\beta_k)=r_{\mathfrak{C}}(\beta_l)^{-1}$.
Suppose that $z_l\in\Omega$, then
\[|z_k|=t^{-1}q^{-j}|z_l^{-1}|\geq q^{-j}t^{-1}r_{\mathfrak{C}}(\beta_k)
> q^{-j}r_C(\beta_k)\geq r_C(\beta_k),
\]
hence $z_k\not\in\Omega$. All the other cases are checked similarly.
\end{proof}

\begin{Lem}\label{independent}
Let $\underline{t}\in V$, $t\in (0,1)$ 
and $f\in A^W$. Then
\begin{equation}\label{Nf}
\iint\limits_{z\in C^n}f(z)d\nu(z;\underline{t};t) 
\end{equation}
is independent of the choice of $\underline{t}$-contour $C$.
\end{Lem}
\begin{proof}
Write $N_f(C)$ for the integral \eqref{Nf}.
We have to show that $N_f(C)=N_f({\mathfrak{C}})$ for arbitrary pairs
of $\underline{t}$-contours $(C,{\mathfrak{C}})$. 

Let ${\mathcal{L}}$ be the collection of pairs of $\underline{t}$-contours
$(C,{\mathfrak{C}})$ for which 
$A^+(C,{\mathfrak{C}})$ is a finite disjoint union
of open intervals and for which $tr_C(x)<r_{\mathfrak{C}}(x)$ for all
$x\in A^+(C,{\mathfrak{C}})$.  Fix a pair
$(C,{\mathfrak{C}})\in {\mathcal{L}}$ 
and let $\Omega$ be the associated domain.
Since the four parameters $t_0,t_1,t_2,t_3$ are in the interior of $C$ and
${\mathfrak{C}}$, we have $w_c(.;\underline{t})\in \mathcal{O}_{\lbrace \pm
1\rbrace}(\Omega)$, and by Lemma \ref{deltaanalytic} we have $\delta(.;t)\in
{\mathcal{O}}_W(\Omega)$. So Cauchy's Theorem implies  that 
$N_f(C)=N_f({\mathfrak{C}})$.

Suppose now that $(C,{\mathfrak{C}})$ is 
an arbitrary pair of ${\underline{t}}$-contours.
Then there exists a finite sequence of $\underline{t}$-contours 
$C_0,C_1,\ldots, C_s$ such that $C_0={\mathfrak{C}}$, $C_s=C$ and such that 
$(C_i,C_{i-1})\in {\mathcal{L}}$ for 
all $i\in\lbrace 1,\ldots,s\rbrace$. 
It follows that $N_f(C)=N_f({\mathfrak{C}})$.
\end{proof}
We define for parameters $\underline{t}\in V$ and $t\in (0,1)$
a symmetric bilinear form $\bigl( .,.\bigr)_{\underline{t},t}$ on $A^W$ by
\begin{equation}\label{symmform}
\bigl( f,g \bigr)_{\underline{t},t}:=\frac{1}{(2\pi i)^n}
\iint\limits_{z\in C^n}f(z)g(z)d\nu(z;\underline{t};t),\quad f,g\in A^W
\end{equation}
where $C$ is an arbitrary $\underline{t}$-contour. 
The bilinear form \eqref{symmform} is independent of the choice of
$\underline{t}$-contour $C$ by Lemma \ref{independent}.
An important tool 
for studying orthogonal polynomials with respect to the bilinear
form $\bigl( .,. \bigr)_{\underline{t},t}$ 
is an explicit second order $q$-difference operator $D=D_{\underline{t},t}$
which preserves the algebra $A^W$ and which is symmetric with respect to
the bilinear form $\bigl( .,. \bigr)_{\underline{t},t}$. The second order
$q$-difference operator $D$
was introduced by Koornwinder \cite{K1} 
and it is explicitly given by 
\begin{equation}\label{secondorderqdiff}
D :=
\sum_{j=1}^n\left(\phi_{j}^+(z)(T_{j}^+-\hbox{Id})+\phi_{j}^-(z)(T_{j}^--
\hbox{Id})\right),
\end{equation}
where $T_{j}^{\pm}$ is the $q^{\pm 1}$-shift in the $j$th coordinate,
\[\left(T_{j}^{\pm}f\right)(z) := 
f(z_1,\ldots,z_{j-1},q^{\pm 1}z_j,z_{j+1},\ldots,z_n),\]
and the functions 
$\phi_{j}^+(z;\underline{t};t)$ and $\phi_{j}^-(z;\underline{t};t)$ 
are given by
\[
\phi_{j}^+(z;\underline{t};t):=
\frac{\prod_{i=0}^3(1-t_iz_j)}{(1-z_j^2)(1-qz_j^2)}
\prod_{l\not=j}\frac{\left(1-tz_lz_j\right)\left(1-tz_l^{-1}z_j\right)}
{\left(1-z_lz_j\right)\left(1-z_l^{-1}z_j\right)},
\]
\[
\phi_{j}^-(z;\underline{t};t):=\phi_j^+(z^{-1};\underline{t};t),
\] 
where we have used the notation $z^{-1}=(z_1^{-1},\ldots,z_n^{-1})$. 
The $BC$ type dominance order on $\Lambda$ is defined by
\begin{equation}
\mu\leq\lambda \quad \Leftrightarrow 
\quad\sum_{j=1}^i\mu_j\leq\sum_{j=1}^i\lambda_j \quad (i=1,\ldots,n)
\end{equation}
for $\lambda,\mu\in \Lambda$.
Koornwinder proved the following triangularity property of $D$
(see \cite[Lemma 5.2]{K1} and the remark after \cite[Proposition 4.1]{SK1}).
\begin{Prop}\label{triaAW}
Let $\lambda\in \Lambda$. For arbitrary $\underline{t}\in {\mathbb{C}}^4$ and 
$t\in {\mathbb{C}}$ we have
\begin{equation}\label{triangulariteit}
Dm_{\lambda} = \sum_{\mu\leq\lambda}E_{\lambda,\mu}m_{\mu}
\end{equation}
with $E_{\lambda,\mu}(\underline{t};t)\in {\mathbb{C}}$ depending polynomially
on $\underline{t}$ and $t$. 
The leading term $E_{\lambda,\lambda}(\underline{t};t)$ will be denoted by
$E_{\lambda}(\underline{t};t)$ and is given by 
\begin{equation}\label{AWb}
E_{\lambda}(\underline{t};t) := 
\sum_{j=1}^n\left(q^{-1}t_0t_1t_2t_3t^{2n-j-1}(q^{\lambda_j}-1)+
t^{j-1}(q^{-\lambda_j}-1)\right).
\end{equation}
\end{Prop}
In particular $D$ maps $A^W$ into itself. The other property of $D$ which we 
already mentioned is the symmetry of $D$ with respect to the bilinear form
$\bigl( .,. \bigr)$, i.e.
\begin{equation}\label{selfadjointAW}
\bigl(D_{\underline{t},t}f,g\bigr)_{\underline{t},t}=
\bigl(f,D_{\underline{t},t}g\bigr)_{\underline{t},t},\quad
f,g\in A^W
\end{equation}
for parameters $\underline{t}\in V$ and $t\in (0,1)$. 
Koornwinder \cite[Lemma 5.3]{K1} proved \eqref{selfadjointAW}
for parameters $\underline{t}$ with $|t_i|<1$ (then the unit circle $T$
can be chosen as $\underline{t}$-contour). By Proposition \ref{triaAW},
\eqref{selfadjointAW} follows for $\underline{t}\in V$ by analytic 
continuation.

We define explicit expressions 
${\mathcal{N}}(\lambda;\underline{t};t)$ for $\lambda\in\Lambda$ by
\begin{equation}\label{normexp}
{\mathcal{N}}(\lambda;\underline{t};t):=2^nn!{\mathcal{N}}^+
(\lambda;\underline{t};t){\mathcal{N}}^-(\lambda;\underline{t};t)
\end{equation}
where ${\mathcal{N}}^+(\lambda):={\mathcal{N}}^+(\lambda;\underline{t};t)$
is given  by
\begin{equation}\label{hatdelta+}
\begin{split}
{\mathcal{N}}^+(\lambda):=
&\prod_{i=1}^n
\frac{\bigl(q^{2\lambda_i-1}t^{2(n-i)}t_0t_1t_2t_3;q\bigr)_{\infty}}
{\bigl(q^{\lambda_i-1}t^{n-i}t_0t_1t_2t_3,
q^{\lambda_i}t^{n-i}t_0t_1, q^{\lambda_i}t^{n-i}t_0t_2,
q^{\lambda_i}t^{n-i}t_0t_3;q\bigr)_{\infty}}\\
&.
\prod_{1\leq j<k\leq n}
\frac{\bigl(q^{\lambda_j+\lambda_k-1}t^{2n-j-k}t_0t_1t_2t_3,
q^{\lambda_j-\lambda_k}t^{k-j};q\bigr)_{\infty}}
{\bigl(q^{\lambda_j+\lambda_k-1}t^{2n-j-k+1}t_0t_1t_2t_3,
q^{\lambda_j-\lambda_k}t^{k-j+1};q\bigr)_{\infty}},\\
\end{split}
\end{equation}
and ${\mathcal{N}}^-(\lambda):={\mathcal{N}}^-(\lambda;\underline{t};t)$
is given by 
\begin{equation}\label{tildedelta+}
\begin{split}
{\mathcal{N}}^-(\lambda):=
&\prod_{i=1}^n
\frac{\bigl(q^{2\lambda_i}t^{2(n-i)}t_0t_1t_2t_3;q\bigr)_{\infty}}
{\bigl(q^{\lambda_i+1}t^{n-i},q^{\lambda_i}t^{n-i}t_1t_2,
q^{\lambda_i}t^{n-i}t_1t_3,q^{\lambda_i}t^{n-i}t_2t_3;q\bigr)_{\infty}}\\
&.\prod_{1\leq j<k\leq n}
\frac{\bigl(q^{\lambda_j+\lambda_k}t^{2n-j-k}t_0t_1t_2t_3,
q^{\lambda_j-\lambda_k+1}t^{k-j};q\bigr)_{\infty}}
{\bigl(q^{\lambda_j+\lambda_k}t^{2n-j-k-1}t_0t_1t_2t_3,
q^{\lambda_j-\lambda_k+1}t^{k-j-1};q\bigr)_{\infty}}.\\
\end{split}
\end{equation}
The following theorem 
extends the results of  Koornwinder in \cite{K1} 
(the orthogonality relations for the multivariable Askey-Wilson polynomials) 
and van Diejen in \cite{vD1} 
(the quadratic norm evaluations 
for the multivariable Askey-Wilson polynomials) 
to parameters $\underline{t}\in V$
(in \cite{K1} and \cite{vD1} the results were obtained for a parameter domain
such that $|t_i|<1$ for all $i$).
\begin{Thm}\label{concl}
Let $\underline{t}\in V$ and $t\in (0,1)$.
There exists a unique basis 
$\lbrace P_{\lambda}(.;\underline{t};t)\rbrace_{\lambda\in\Lambda}$ of $A^W$ 
such that
\begin{equation}\label{propertya}
P_{\lambda}(.;\underline{t};t)=m_{\lambda}+\sum_{\mu<\lambda}
c_{\lambda,\mu}(\underline{t},t)m_{\mu}, \quad 
c_{\lambda,\mu}(\underline{t},t)\in {\mathbb{C}}
\end{equation}
\begin{equation}\label{propertyb}
\bigl(P_{\lambda}(.;\underline{t};t),P_{\mu}(.;\underline{t};t)\bigr)_{
\underline{t},t}=0 \quad \hbox{ if }\, \mu\not=\lambda.
\end{equation}
Furthermore, $P_{\lambda}(.;\underline{t};t)$   
is an eigenfunction of $D_{\underline{t},t}$ with eigenvalue
$E_{\lambda}(\underline{t};t)$ and we have the explicit evaluation formula
\begin{equation}\label{normAWdomainW}
\bigl(P_{\lambda}(.;\underline{t};t),P_{\lambda}(.;\underline{t};t)\bigr)_{
\underline{t},t}
=\mathcal{N}(\lambda;\underline{t};t)
\end{equation}
for the quadratic norms of the polynomials $P_{\lambda}$.
\end{Thm}
\begin{Def}\label{AWdefinitie}
We call $P_{\lambda}(.;\underline{t};t)$ the monic multivariable
Askey-Wilson polynomial of degree $\lambda$.
\end{Def}
We end this section with a sketch of the proof for Theorem \ref{concl}
using the techniques of Koornwinder \cite{K1} 
and van Diejen \cite{vD1}. For more details, we refer
the reader to these two papers. 

We fix arbitrary $0\not=\nu\in\Lambda$.  It is sufficient to 
prove the existence and uniqueness of a
set of $W$-invariant Laurent polynomials $\lbrace
P_{\lambda}(.;\underline{t};t)\rbrace_{\lambda\leq\nu}$
satisfying \eqref{propertya} and \eqref{propertyb} for $\lambda,\mu\leq\nu$
and to prove the remaining assertions of the theorem  
for the polynomials 
$\lbrace P_{\lambda}(.;\underline{t};t)\rbrace_{\lambda\leq\nu}$.

We first define the 
polynomials $\lbrace P_{\lambda}\rbrace_{\lambda\leq\nu}$
for a dense parameter domain $U_{\nu}\subset V\times (0,1)$. 
The subset $U_{\nu}$ is by definition the set of parameters 
$(\underline{t},t)\in V\times (0,1)$
such that $E_{\mu}(\underline{t};t)\not=E_{\lambda}(\underline{t};t)$ for all
$\lambda,\mu\leq\nu$, $\lambda\not=\mu$. Note that
$U_{\nu}\subset V\times (0,1)$ is open and dense since the eigenvalues
$\lbrace E_{\lambda}(\underline{t};t)\rbrace_{\lambda\in\Lambda}$
are mutually different as polynomials in the parameters $\underline{t},t$.

The polynomials $P_{\lambda}(.;\underline{t};t)\in A^W$ for
$(\underline{t},t)\in U_{\nu}$ and $\lambda\leq\nu$ are defined by
\begin{equation}\label{prodform}
P_{\lambda}(.;\underline{t};t):=\left(\prod_{\mu<\lambda}
\frac{D_{\underline{t},t}-E_{\mu}(\underline{t};t)}
{E_{\lambda}(\underline{t};t)-E_{\mu}(\underline{t};t)}\right)m_{\lambda}
\end{equation} 
(cf. \cite{M1}, \cite{SK1}). 
By Proposition \ref{triaAW}, $P_{\lambda}(.;\underline{t};t)$ is the
unique eigenfunction of $D_{\underline{t},t}$ with eigenvalue 
$E_{\lambda}(\underline{t};t)$ which satisfies \eqref{propertya}, 
and by the symmetry of $D$  \eqref{selfadjointAW} we obtain 
the orthogonality relations \eqref{propertyb} for the polynomials $\lbrace
P_{\lambda}(.;\underline{t};t)\rbrace_{\lambda\leq\nu}$ 
and for parameters $(\underline{t},t)\in U_{\nu}$.

Next we establish the quadratic norm evaluations \eqref{normAWdomainW}
for $\lbrace P_{\lambda}(.;\underline{t};t)\rbrace_{\lambda\leq\nu}$.
In the special case $\lambda=0$, \eqref{normAWdomainW} reduces to 
\begin{equation}\label{norm0}
\begin{split}
\bigl(1,1\bigr)_{\underline{t},t}&={\mathcal{N}}(0;\underline{t};t)\\
&=2^nn!\prod_{i=1}^n\frac{\bigl(t,t^{2n-i-1}t_0t_1t_2t_3;q\bigr)_{\infty}}
{\bigl(q,t^{n-i+1};q\bigr)_{\infty}\prod_{0\leq j<k\leq
3}\bigl(t^{n-i}t_jt_k;q\bigr)_{\infty}}\\
\end{split}
\end{equation}
which was proved by Gustafson \cite{G2}
for parameters $t\in (0,1)$ and $\underline{t}\in {\mathbb{C}}^4$ with
$|t_i|<1$ (since then the torus $T$ can be chosen as $\underline{t}$-contour). 
The second equality follows by a straightforward computation
(see also \cite{M2}). By analytic continuation, 
\eqref{norm0} is valid for parameters $t\in (0,1)$
and $\underline{t}\in V$.

For general $\lambda$, van Diejen \cite{vD1} proved explicit Pieri formulas
for the renormalized Askey-Wilson polynomials 
\begin{equation}\label{renorm}
p_{\lambda}(.;\underline{t};t):=
c(\lambda;\underline{t};t)P_{\lambda}(.;\underline{t};t)
,\quad c(\lambda;\underline{t};t):=
\frac{{\mathcal{N}}^{+}(0)}{{\mathcal{N}}^+(\lambda)}
\prod_{j=1}^n(t_0t^{n-j})^{\lambda_j}.
\end{equation}
Note that the 
renormalization constant $c(\lambda;\underline{t};t)$ is a rational
expression in the parameters 
$\underline{t},t$. The Pieri formulas for $p_{\lambda}$
are explicit expressions for the coefficients 
$d_{\lambda}^{(r)}(\mu;\underline{t};t)$ in the expansions
\begin{equation}\label{expansionPieri}
E_r(z;\underline{t};t)p_{\lambda}(z;\underline{t};t)=
\sum_{\mu\leq\lambda+\omega_n}
d_{\lambda}^{(r)}(\mu;\underline{t};t)p_{\mu}(z;\underline{t};t),\quad
(r=1,\ldots,n)
\end{equation}
where $\lbrace E_r(z;\underline{t};t)\rbrace_{r=1}^n$ 
are explicit algebraic generators of the algebra of $W$-invariant 
Laurent polynomials $A^W$ and where $\omega_n:=(1,\ldots,1)\in\Lambda$
is the $n$th fundamental weight (see \cite{vD1} for the explicit formulas,
or \cite[Appendix B]{vDS} where the notations are closer to the ones used
in this paper). 
Since everything in \eqref{expansionPieri} depends rationally on the parameters
$(\underline{t},t)$, we may view \eqref{expansionPieri} as an identity
in the algebra of $W$-invariant Laurent polynomials over the quotient field
${\mathbb{C}}(\underline{t},t)$, where $\underline{t},t$ are considered as
indeterminates.

The Pieri formulas and the orthogonality relations for
the renormalized Askey-Wilson polynomials with real parameters $t_i$ and 
$|t_i|<1$ allowed van Diejen to reduce the 
the norm computation for arbitrary $\lambda$ to the case 
$\lambda=0\in \Lambda$. 
Gustafson's evaluation \eqref{norm0} then completes the evaluation 
for general $\lambda$. 

Exactly the same reduction can now be done for the norm evaluations of the
polynomials $\lbrace p_{\lambda}(.;\underline{t};t) \rbrace_{\lambda\leq\nu}$
for parameters $(\underline{t},t)\in U_{\nu}$.
Indeed, by taking a dense subset of $U_{\nu}$ if necessary,
we may assume 
that for all $(\underline{t},t)\in U_{\nu}$ the following conditions are
fulfilled,\\
(i) the Askey-Wilson polynomials $P_{\lambda}(.;\underline{t};t)$ 
are well defined and mutually
orthogonal with respect to $\bigl( .,. \bigr)_{\underline{t},t}$ 
for $\lambda\leq\nu+\omega_n$,\\
(ii) the renormalization constants
$c(\lambda;\underline{t};t)$ are well defined and non zero for
$\lambda\leq\nu+\omega_n$,\\ 
(iii) the Pieri formulas 
\eqref{expansionPieri} are valid for all $\lambda\leq\nu$.\\
The quadratic norm evaluations for the polynomials
$\lbrace p_{\lambda}(.;\underline{t};t)\rbrace_{\lambda\leq\nu}$ 
with parameter values $(\underline{t},t)\in U_{\nu}$ 
can then be reduced to the case $\lambda=0$ using exactly the
same arguments as in \cite[Theorem 4]{vD1}. The extension of
Gustafson's result 
\eqref{norm0} then completes the proof of \eqref{normAWdomainW}
for $\lbrace P_{\lambda}(.;\underline{t};t)\rbrace_{\lambda\leq\nu}$
and $(\underline{t},t)\in U_{\nu}$. 

Now note that the functions
\begin{equation}\label{contfunct}
(\underline{t},t)\mapsto {\mathcal{N}}^{\pm}(\lambda;\underline{t};t) : 
\quad V\times (0,1)\rightarrow {\mathbb{C}}
\end{equation}
are well defined 
continuous functions which do not have zeros on the domain $V\times
(0,1)$. This is immediately clear except for 
${\mathcal{N}}^+(\lambda)$ with $\lambda\in\Lambda$ and
$\lambda_n=0$. But then the expression for ${\mathcal{N}}^+(\lambda)$
can be simplified, similarly as the simplification of the expression for 
${\mathcal{N}}(0)$ in  \eqref{norm0}, from which
it follows that ${\mathcal{N}}^+(\lambda;\underline{t};t)$ is a well defined, 
continuous function of $(\underline{t},t)\in V\times (0,1)$ without zeros.
It follows in particular that the quadratic norms of the polynomials 
$\lbrace P_{\lambda}(.;\underline{t};t)
\rbrace_{\lambda\leq\nu}$ are non zero for $(\underline{t},t)\in U_{\nu}$.
Hence \eqref{propertya} and \eqref{propertyb} for $\lambda,\mu\leq\nu$ 
uniquely characterize the set of polynomials 
$\lbrace P_{\lambda}(.;\underline{t};t)\rbrace_{\lambda\leq\nu}$
when $(\underline{t},t)\in U_{\nu}$. 

It follows that the Askey-Wilson 
polynomial $P_{\lambda}$ satisfies the following Gram-Schmidt formula,
\begin{equation}\label{clm}
P_{\lambda}(z;\underline{t};t)=
m_{\lambda}(z)-\sum_{\mu<\lambda}
\frac{\bigl(m_{\lambda},P_{\mu}(.;\underline{t};t)\bigr)_{
\underline{t},t}}
{\mathcal{N}(\mu;\underline{t};t)}
P_{\mu}(z;\underline{t};t)
\end{equation}
for $(\underline{t},t)\in U_{\nu}$ and $\lambda\leq\nu$. 
By induction, we conclude from \eqref{clm} that the coefficients
$c_{\lambda,\mu}: U_{\nu}\rightarrow {\mathbb{C}}$ in \eqref{propertya}
uniquely extend to continuous functions
$c_{\lambda,\mu}: V\times (0,1)\rightarrow {\mathbb{C}}$ for all
$\mu<\lambda\leq\nu$.
(Compare with the proof of \cite[Proposition 6.3]{vDS}.)
Hence existence and uniqueness of $\lbrace
P_{\lambda}(.;\underline{t};t)\rbrace_{\lambda\leq\nu}$ as well as 
the other assertions 
follow now by continuity for all $(\underline{t},t)\in V\times
(0,1)$. This completes the proof of the theorem.
\begin{Rem}
For $n=1$
the polynomials $\lbrace P_{\lambda}(z;\underline{t}) 
\, | \, \lambda\in {\mathbb{N}}_0\rbrace$ are independent of $t$ and
are known as (monic, one-variable) Askey-Wilson polynomials.
Theorem \ref{concl} reduces 
to the orthogonality relation and quadratic norm evaluation 
stated in \cite[Theorem 2.3]{AW}. 
The polynomials $P_{\lambda}(z;\underline{t})$
$(\lambda\in {\mathbb{N}}_0)$ 
can then be given explicitly in terms of the basic
hypergeometric series
\[
{}_{s+1}\phi_s\left( \begin{array}{c}
             a_1,\ldots ,a_{s+1} \\ b_1,\ldots ,b_s
            \end{array} ; q,z \right) =
\sum_{m=0}^\infty
\frac{(a_1,\ldots ,a_{s+1};q)_m}{(b_1,\ldots ,b_s,q;q)_m}z^m
\]
as
\begin{equation}\label{explicitexpression}
P_\lambda (z;\underline{t}) =
\frac{(t_0t_1, t_0t_2, t_0t_3; q)_\lambda }{t_0^{\lambda}
(t_0t_1t_2t_3q^{\lambda -1};q)_\lambda }\;
{}_4\phi_3
\left(
\begin{matrix}
q^{-\lambda },\; q^{\lambda -1}t_0t_1t_2t_3,\;
t_0z,\; t_0z^{-1} \\ [0.5ex]
t_0t_1 ,\; t_0t_2 ,\; t_0t_3
\end{matrix} \; ; q,q \right) 
\end{equation}
(see \cite{AW} for details). The renormalized Askey-Wilson polynomial
$p_{\lambda}(z;\underline{t})$ \eqref{renorm} is then exactly the ${}_4\phi_3$
part of \eqref{explicitexpression}. 
\end{Rem}
The renormalization constant $c(\lambda;\underline{t};t)$ \eqref{renorm}
is easily seen 
to be regular and non zero at $(\underline{t},t)\in V\times (0,1)$
for all $\lambda\in\Lambda$. Hence the renormalized Askey-Wilson polynomials
$\lbrace p_{\lambda}(z;\underline{t};t)\rbrace_{\lambda\in\Lambda}$
form an orthogonal basis of $A^W$ with respect to the bilinear form $\bigl( .,.
\bigr)_{\underline{t},t}$ 
for all parameter values $(\underline{t},t)\in V\times
(0,1)$.

\section{Residue calculus for the orthogonality measure $d\nu$}
In this section 
we develop for specific pairs of deformed circles $(C,{\mathfrak{C}})$
a residue calculus for integrals of the form
\begin{equation}\label{uitgangspunt}
\frac{1}{(2\pi i)^n}
\iint\limits_{z\in C^n}f(z)d\nu(z)=
\frac{1}{(2\pi i)^n}\iint\limits_{z\in C^n}f(z)\Delta(z)\frac{dz}{z},\quad
f\in {\mathcal{O}}_W(\Omega^n)
\end{equation}
when $C^n$ is shifted to ${\mathfrak{C}}^n$. Here $\Omega\subset
{\mathbb{C}}^*$ is the domain associated with the pair $(C,{\mathfrak{C}})$.
We will develop the residue calculus 
for pairs of contours $(C,{\mathfrak{C}})$
such that the poles only depend on $q,t$
and one of the four parameters $t_0,t_1,t_2,t_3$.
By the symmetry 
of $\Delta(z;\underline{t};t)$ in the parameters $t_0,t_1,t_2,t_3$,
we may assume 
without loss of generality that the poles only depend on $q,t$ and $t_0$.
We use in this section the notations of Definition \ref{contour}. 
So for a deformed circle $C$,
we write $\phi_C(x)=r_C(x)e^{2\pi i x}$ for its parametrization.

We fix in this section $t\in (0,1)$ and $t_1,t_2,t_3\in {\mathbb{C}}^*$ 
such that $\#\lbrace\alpha_i^+,\alpha_i^- \, |
\, i=1,2,3\rbrace=6$ ($\alpha_i^{\pm}$ given by \eqref{alpha}). 
We simplify notations by omiting the dependance on the parameters $t_1,t_2,t_3$
and $t$ in this section. For instance, we will write $w_c(x;t_0)$ instead of
$w_c(x;\underline{t})$, etc.
In the next definition, we introduce the pairs of contours for which we
will develop the residue calculus.
\begin{Def}\label{Definitionpairs}
Let $t_0\in {\mathbb{C}}^*$ such that
$\underline{t}=(t_0,t_1,t_2,t_3)\in V$ ($V$ given by Definition \ref{WAW}).
A pair of contours $(C,\mathfrak{C})$ is called
a $(n,t_0)$-residue pair if $C$ and $\mathfrak{C}$ 
are deformed circles  such that\\
{\bf (i)} the subset $A^{+}(C,\mathfrak{C})$ \eqref{A+} 
is an open interval for which $\alpha_0^+\in A^+(C,{\mathfrak{C}})$
but $\alpha_i^{\pm}\notin A^+(C,{\mathfrak{C}})$\, ($i=1,2,3$);\\
{\bf (ii)} $t_iq^r\notin C\cup {\mathfrak{C}}$ 
for $r\in {\mathbb{N}}_0$ and
$i\in \lbrace 0,\ldots,3\rbrace$;\\
{\bf (iii)} 
$t_0t^pq^r\notin \mathfrak{C}$ for $p\in \lbrace -1,\ldots,n-1\rbrace$ 
and $r\in {\mathbb{Z}}$.\\
\end{Def}
Note that a $(n,t_0)$-residue pair $(C,{\mathfrak{C}})$ is a
$(r,t_0)$-residue pair for all $r\in {\mathbb{N}}_0$ with $r\leq n$.
We define now the measures $d\nu_r(z;t_0)$\, ($r=1,\ldots,n$) 
which will appear when the contour $C^n$ in \eqref{uitgangspunt}
is shifted to ${\mathfrak{C}}^n$. First we need to
introduce some more notations.

We set 
\begin{equation}\label{Pr}
P(r):=
\lbrace \lambda\in {\mathbb{N}}_0^r \, | \, \lambda_1\leq\lambda_2\leq\ldots
\leq\lambda_r \rbrace
\end{equation}
and we set $\lambda_0:=0$ for arbitrary $\lambda\in P(r)$. We write
$\rho_i=\rho_i(t_0;t):=t_0t^{i-1}$ for $i\in
{\mathbb{Z}}$ and, for $\lambda\in P(r)$,
\begin{equation}
\rho q^{\lambda}:=(\rho_1q^{\lambda_1},\rho_2q^{\lambda_2},
\ldots, \rho_rq^{\lambda_r})=(t_0q^{\lambda_1},t_0tq^{\lambda_2},\ldots,
t_0t^{r-1}q^{\lambda_r}).
\end{equation}
Define
$D(r)=D(r;C,{\mathfrak{C}};t_0)$ for $r=1,\ldots,n$ by
\begin{equation}\label{Dparameter}
D(r):=\lbrace \rho q^{\lambda} \, | \, \lambda\in P(r) \hbox{ and }
r_C(\alpha_0^+)>|\rho_iq^{\lambda_i}|>r_{\mathfrak{C}}(\alpha_0^+) \quad
(i=1,\ldots,r) \rbrace.
\end{equation}
Note that for $\omega\in D(r)$, 
we have $\omega_i\in\hbox{int}\bigl(\Omega\bigr)$
for all $i$, where $\hbox{int}\bigl(\Omega\bigr)$ is the interior of $\Omega$. 
For $\rho q^{\lambda}\in D(r)$, we set
\begin{equation}\label{discreteweights}
\Delta^{(d)}(\rho q^{\lambda};t_0):=
\left(\prod_{j=1}^rw_d\bigl(\rho_jq^{\lambda_j};
\rho_jq^{\lambda_{j-1}}\bigr)\right)\delta_d\bigl(\rho q^{\lambda}\bigr)
\end{equation}
with
\begin{equation}\label{resonevar}
w_d(\rho_jq^{\lambda_j};\rho_jq^{\lambda_{j-1}}):=
\hbox{res}_{x=\rho_jq^{\lambda_j}}
\left(\frac{w_{c}(x;\rho_jq^{\lambda_{j-1}},t_1,t_2,t_3)}{x}\right)
\end{equation}
where $w_c$ is given by \eqref{wc} 
and with interaction factor
\begin{equation}\label{discreteinteraction}
\delta_{d}(\rho q^{\lambda}):=
\prod_{1\leq k<l\leq r}\frac{\bigl(\rho_k^{-1}\rho_lq^{\lambda_l-\lambda_k},
\rho_k^{-1}\rho_l^{-1}q^{-\lambda_k-\lambda_l};q\bigr)_{\tau}}
{\bigl(\rho_k\rho_lq^{\lambda_{k-1}+\lambda_l},
\rho_k\rho_l^{-1}q^{\lambda_{k-1}-
\lambda_l};q\bigr)_{\lambda_k-\lambda_{k-1}}} \quad (t=q^{\tau}).
\end{equation}
The discrete parts of the measure which we will 
obtain by deforming  $C^n$ to ${\mathfrak{C}}^n$ in \eqref{uitgangspunt}
will involve the weights $\Delta^{(d)}$. 
Note that for $r=1$ and $\rho q^{\lambda}=t_0q^i\in D(1)$, we have
$\Delta^{(d)}(t_0q^i;t_0)=w_d(t_0q^i;t_0)$.  
Note furthermore that 
the weight function $w_c(x;\rho_jq^{\lambda_{j-1}})/x$ has a simple pole
in $\rho_jq^{\lambda_j}$ since $\lambda_j-\lambda_{j-1}\in {\mathbb{N}}_0$
and $\underline{t}\in V$. In fact,
for $\underline{\tau}=(\tau_0,\tau_1,\tau_2,\tau_3)\in V$ 
and for $i\in {\mathbb{N}}_0$, the
discrete weight
\begin{equation}\label{residuen}
w_{d}(\tau_0q^i;\tau_0)=
w_{d}\left(\tau_0q^i;\tau_0;\tau_1,\tau_2,\tau_3\right)=
\hbox{ res}_{x=\tau_0q^i}\left(\frac{w_{c}(x;\underline{\tau})}{x}\right)
\end{equation}
can be explicitly given by the formula
\begin{equation}\label{weightfunctiondisc}
\begin{split}
w_{d}(\tau_0q^i;\tau_0) 
=&\frac{\left(\tau_0^{-2};q\right)_{\infty}}
{\left(q,\tau_0\tau_1, \tau_1/\tau_0, \tau_0\tau_2,
\tau_2/\tau_0, \tau_0\tau_3, \tau_3/\tau_0;q \right)_{\infty}}\\
&.\frac{\left(\tau_0^2,\tau_0\tau_1,\tau_0\tau_2,\tau_0\tau_3;q\right)_i} 
{\left(q,\tau_0q/\tau_1,\tau_0q/\tau_2,
\tau_0q/\tau_3;q\right)_i}\frac{\left(1-\tau_0^2q^{2i}\right)}
{\left(
1-\tau_0^2\right)}\left(\frac{q}{\tau_0\tau_1\tau_2\tau_3}
\right)^i\\
\end{split}
\end{equation}
(see \cite[Theorem 2.4]{AW} or \cite[(7.5.22)]{GR} to avoid a small misprint). 
For $\rho q^{\lambda}\in D(r)$ and $z\in {\mathfrak{C}}^{n-r}$, we set
\begin{equation}\label{dnur}
d\nu_r(\rho q^{\lambda},z;t_0)=
\Delta_r(\rho q^{\lambda},z;t_0)\frac{dz}{z}
\end{equation}
with weight function $\Delta_r(\rho q^{\lambda},z;t_0)$ 
given by
\begin{equation}\label{weightalgt}
\Delta_r(\rho q^{\lambda},z;t_0):=\Delta^{(d)}\bigl(\rho q^{\lambda};t_0\bigr)
\Delta(z;t_0)\delta_c\bigl(\rho q^{\lambda};z\bigr)
\end{equation}
where $\Delta(z;t_0)$ is the weight function \eqref{opdeling+}
in the variables $z=(z_1,\ldots,z_{n-r})$ and where
$\delta_c(\rho q^{\lambda};z)$ is an interaction factor given by
\begin{equation}\label{continuousinteraction}
\delta_c(\rho q^{\lambda};z):=
\prod_{\stackrel{1\leq k\leq r}{1\leq l\leq n-r}}
\bigl(\rho_kq^{\lambda_k}z_l,\rho_kq^{\lambda_k}z_l^{-1},
\rho_k^{-1}q^{-\lambda_k}z_l,\rho_k^{-1}q^{-\lambda_k}z_l^{-1};q\bigr)_{\tau}
\end{equation}
with $t=q^{\tau}$. In particular we have $\Delta_n(\rho q^{\lambda};t_0)=
\Delta^{(d)}(\rho q^{\lambda};t_0)$ for $\rho q^{\lambda}\in D(n)$.
The measure $d\nu_r(\rho q^{\lambda},z;t_0)$ 
is well defined on $D(r)\times {\mathfrak{C}}^{n-r}$ 
since the denominator of $\Delta_r(\rho q^{\lambda},z;t_0)$
is non zero 
by properties {\bf (ii)} and {\bf (iii)} of the $(n,t_0)$-residue pair
$(C,{\mathfrak{C}})$ (Definition \ref{Definitionpairs}). We call $d\nu_r$ 
the $r$th 
{\it measure associated with the $(n,t_0)$-residue pair} $(C,\mathfrak{C})$. 
We have the following proposition.
\begin{Prop}\label{restau=tau}
Let $(C,{\mathfrak{C}})$ be a $(n,t_0)$-residue pair
and let $\Omega=\Omega(C,{\mathfrak{C}})$ be the associated domain.
Let $d\nu_r$ be the $r$th measure associated with
$(C,\mathfrak{C})$, then 
\begin{equation}\label{eindres}
\begin{split}
\frac{1}{(2\pi i)^n}\iint\limits_{z\in C^n}f(z)&d\nu(z)=
\frac{1}{(2\pi i)^n}\iint\limits_{z\in {\mathfrak{C}}^n}f(z)d\nu(z)\\
+&\sum_{r=1}^{n}
\frac{2^r\bigl(n-r+1\bigr)_r}{(2\pi i)^{n-r}}\sum_{\omega\in D(r)}\,\,
\iint\limits_{z\in{\mathfrak{C}}^{n-r}}
f(\omega,z)d\nu_r(\omega,z)\\
\end{split}
\end{equation}
for $f\in {\mathcal{O}}_W(\Omega^n)$ where 
$\bigl(u\bigr)_r:=\prod_{i=0}^{r-1}(u+i)$ is the shifted factorial.
\end{Prop}
In the next lemma we will give the proof of Proposition \ref{restau=tau}
for $(n,t_0)$-residue pairs $(C,{\mathfrak{C}})$ such that
the interaction factor $\delta(.;t)$ is analytic on 
$\bigl(\Omega(C,{\mathfrak{C}})\bigr)^n$.
\begin{Lem}\label{restau=tau1}
Suppose that $(C,{\mathfrak{C}})$ is a $(n,t_0)$-residue pair such that
\begin{equation}\label{condcontour}
tr_{C}(x)<r_{\mathfrak{C}}(x), \qquad \forall x\in A^+(C,{\mathfrak{C}}).
\end{equation}
Then \eqref{eindres} is valid.
\end{Lem}
\begin{proof}
Fix a $(n,t_0)$-residue 
pairs $(C,{\mathfrak{C}})$ satisfying the extra condition
\eqref{condcontour}. 
We will prove by induction on $l\in \lbrace 0,\ldots,n\rbrace$ that
\begin{equation}\label{step2r}
\begin{split}
\frac{1}{(2\pi i)^n}\iint\limits_{z\in C^l\times
{\mathfrak{C}}^{n-l}}f(z)\Delta(z;t_0)\frac{dz}{z}
= &\frac{1}{(2\pi i)^n}\iint\limits_{z\in {\mathfrak{C}}^n}
f(z)\Delta(z;t_0)\frac{dz}{z}\\
+
\frac{2l}{(2\pi i)^{n-1}}
\sum_{i\in I}&\,\,\iint\limits_{z\in {\mathfrak{C}}^{n-1}}
f(t_0q^i,z)\Delta_1(t_0q^i,z;t_0)\frac{dz}{z}\\
\end{split}
\end{equation}
for $f\in {\mathcal{O}}_W(\Omega^n)$, where 
\begin{equation}\label{I0}
I=\lbrace i\in {\mathbb{N}}_0 \, | \, 
r_{C}(\alpha_0^+) > |t_0q^i| > r_{\mathfrak{C}}(\alpha_0^+)\rbrace.
\end{equation}
Then \eqref{eindres} is the special case $l=n$ in
\eqref{step2r}, since 
$D(r)=\emptyset$ for $r>1$ by \eqref{condcontour}. 

For $l=0$, \eqref{step2r} is trivial. Let $l\in \lbrace 1,\ldots,n \rbrace$.
Since $\delta(.;t)\in {\mathcal{O}}_W(\Omega^n)$ 
by Lemma \ref{deltaanalytic}, we 
can shift $C$ to ${\mathfrak{C}}$ 
for the first variable $z_1$ in \eqref{uitgangspunt}
and we obtain by \eqref{residuen}, by the $W$-invariance of $\Delta(z)$ 
and by Cauchy's Theorem for 
$f\in {\mathcal{O}}_W(\Omega^n)$,
\begin{equation}\label{step1r}
\begin{split}
\frac{1}{(2\pi i)^n}&\iint\limits_{z\in C^l\times
{\mathfrak{C}}^{n-l}}f(z)\Delta(z;t_0)\frac{dz}{z}
= \frac{1}{(2\pi i)^n}\iint\limits_{z\in C^{l-1}\times {\mathfrak{C}}^{n-l+1}}
f(z)\Delta(z;t_0)\frac{dz}{z}\\
&+
\frac{2}{(2\pi i)^{n-1}}\sum_{i\in I}\,\,\iint\limits_{z\in C^{l-1}\times 
{\mathfrak{C}}^{n-l}}
f(t_0q^i,z)\Delta_1(t_0q^i,z;t_0)\frac{dz}{z}.\\
\end{split}
\end{equation}
Here we have used that the residue at $z_1=t_0^{-1}q^{-i}$ ($i\in I$)
of $\Delta(z_1,z';t_0)/z_1$ is equal to $-\Delta_1(t_0q^i,z';t_0)$ since
\[\hbox{ res}_{x=t_0^{-1}q^{-i}}
\left(\frac{w_{c}(x;\underline{t})}{x}\right)=
-w_d(t_0q^i;t_0;t_1,t_2,t_3).
\]
The weight function $\Delta_1(t_0q^i,z;t_0)$ in \eqref{step1r}
can be rewritten as 
\begin{equation}\label{D1}
\Delta_1(t_0q^i,z;t_0)=h(t_0q^i,z;t_0)\Delta(z;tt_0q^i)
\end{equation}
with 
\begin{equation}\label{hn}
h(t_0q^i,z;t_0)=w_d(t_0q^i;t_0)
\prod_{s=1}^{n-1}\frac{\bigl(t_0^{-1}q^{-i}z_s,
t_0^{-1}q^{-i}z_s^{-1};q\bigr)_{\tau}}
{\bigl(t_0z_s,t_0z_s^{-1};q\bigr)_i} \quad (t=q^{\tau}).
\end{equation}
This follows by interchanging the factor
$\bigl(t_0z_s,t_0z_s^{-1};q\bigr)_{\infty}$ in the denominator
of $w_c(z_s;t_0)$ with the factor
$\bigl(tt_0q^iz_s,tt_0q^iz_s^{-1};q\bigr)_{\infty}$
in the denominator of $\delta_c(t_0q^i;z)$ for $s=1,\ldots,n-1$.
We have $\Delta(.;tt_0q^i)\in
{\mathcal{O}}_W(\Omega^{n-1})$ for $i\in I$ by \eqref{condcontour}
and Lemma \ref{deltaanalytic}.
We claim that $h(t_0q^i,.;t_0)\in {\mathcal{O}}_W(\Omega^{n-1})$ for $i\in I$.
Indeed, it is sufficient to check that the map
\begin{equation}\label{check1}
x\mapsto
\frac{\bigl(t_0^{-1}q^{-i}x,t_0^{-1}q^{-i}x^{-1};q\bigr)_{\infty}}
{\bigl(t_0x,t_0x^{-1};q\bigr)_i
\bigl(tt_0^{-1}q^{-i}x,tt_0^{-1}q^{-i}x^{-1};q\bigr)_{\infty}}
\end{equation}
is analytic on $\Omega$ when $i\in I$.
The zeros of the factor $\bigl(t_0x,t_0x^{-1};q\bigr)_i$
in the denominator are compensated by zeros in the numerator.
Next, we check that $\bigl( tt_0^{-1}q^{-i}x;q\bigr)_{\infty}$ is non zero
for $x\in \Omega$ and $i\in I$. Now we have that 
$\bigl( tt_0^{-1}q^{-i}x;q\bigr)_{\infty}=0$ iff $x=t^{-1}q^{i-m}t_0$
for some $m\in {\mathbb{N}}_0$. In particular, we must have
$\hbox{arg}(x)=\hbox{arg}(t_0)=2\pi \alpha_0^+$.
Since $i\in I$, we have for $m\in {\mathbb{N}}_0$,
\[
|t^{-1}q^{i-m}t_0|>t^{-1}r_{\mathfrak{C}}(\alpha_0^+)q^{-m}\geq
t^{-1}r_{\mathfrak{C}}(\alpha_0^+)>r_{C}(\alpha_0^+),
\]
where the last 
inequality is obtained from the extra condition \eqref{condcontour}.
Since $x\in \Omega$ with $\hbox{arg}(x)=2\pi \alpha_0^+$ implies that
$r_{\mathfrak{C}}(\alpha_0^+)\leq |x|\leq r_{C}(\alpha_0^+)$,
we conclude that $\bigl( tt_0^{-1}q^{-i}x;q\bigr)_{\infty}\not=0$
for $x\in \Omega$ and $i\in I$. 
Since $\Omega^{-1}=\Omega$, we then also have
$\bigl( tt_0^{-1}q^{-i}x^{-1};q\bigr)_{\infty}\not=0$
for $x\in \Omega$ and $i\in I$. Thus the map given by
\eqref{check1} is analytic on $\Omega$ if $i\in I$.
In particular, 
\[
f(t_0q^i,.)\Delta_1(t_0q^i,.;t_0)
\in {\mathcal{O}}_W(\Omega^{n-1})
\]
for $i\in I$, so we obtain by Cauchy's Theorem and \eqref{step1r},
\begin{equation}\label{step2rr}
\begin{split}
\frac{1}{(2\pi i)^n}\iint\limits_{z\in C^l\times
{\mathfrak{C}}^{n-l}}f(z)&\Delta(z;t_0)\frac{dz}{z}
= \frac{1}{(2\pi i)^n}\iint\limits_{z\in C^{l-1}\times {\mathfrak{C}}^{n-l+1}}
f(z)\Delta(z;t_0)\frac{dz}{z}\\
&+\frac{2}{(2\pi i)^{n-1}}
\sum_{i\in I}\,\,\iint\limits_{z\in {\mathfrak{C}}^{n-1}}
f(t_0q^i,z)\Delta_1(t_0q^i,z;t_0)\frac{dz}{z}\\
\end{split}
\end{equation}
for $f\in {\mathcal{O}}_W(\Omega^n)$. Then \eqref{step2r} follows by
applying the induction hypotheses on the integral over $C^{l-1}\times
{\mathfrak{C}}^{n-l+1}$ in \eqref{step2rr}.  
\end{proof}
Lemma \ref{restau=tau1} can be used to prove Proposition \ref{restau=tau}
inductively. 
The following definition will be used to formulate the induction hypotheses. 
\begin{Def}
Let $(C,{\mathfrak{C}})$ be a $(n,t_0)$-residue pair and let
$A^+(C,{\mathfrak{C}})$ \eqref{A+} be the associated
open interval.
A sequence of closed contours $(C_0,\ldots,C_s)$
is called a $(n,t_0)$-resolution for $(C,{\mathfrak{C}})$ 
if the contours $C_l$ are deformed circles
satisfying the following four conditions (we write $r_l$ for the
(radial) functions $r_{C_l}$ in the parametrization $\phi_{C_l}$ of $C_l$):\\
{\bf (i)} $C_0={\mathfrak{C}}$ and $C_s=C$;\\
{\bf (ii)} 
$r_l(x)=r_{{\mathfrak{C}}}(x)=r_C(x)$ for $x\notin A^+(C,{\mathfrak{C}})\cup
A^+({\mathfrak{C}},C)$
and $l\in\lbrace 0,\ldots,s\rbrace$;\\
{\bf (iii)}
$tr_{l+1}(x)<r_l(x)<r_{l+1}(x)$ 
for $x\in A^+(C,{\mathfrak{C}})$ and $l\in\lbrace 
0,\ldots,s-1\rbrace$;\\
{\bf (iv)} $t_0t^pq^r\notin C_l$ for $p\in\lbrace -1,\ldots,n-1\rbrace$,
$r\in {\mathbb{Z}}$ and $l\in\lbrace 1,\ldots,s-1\rbrace$.

We call $s$ the {\it length} of the resolution.
\end{Def}
Note that there 
exists a $(n,t_0)$-resolution 
for every $(n,t_0)$-residue pair $(C,{\mathfrak{C}})$.
If  $(C_0,\ldots,C_s)$ is a $(n,t_0)$-resolution  for a $(n,t_0)$-residue
pair $(C,{\mathfrak{C}})$, then $(C_l,C_{l-1})$ is a $(n,t_0)$-residue pair
satisfying the extra condition \eqref{condcontour} used to prove Lemma
\ref{restau=tau1} ($l\in \lbrace 1,\ldots,s\rbrace$).
We will prove now Proposition \ref{restau=tau} by induction on 
the length of the resolution.
\begin{proof}[Proof of Proposition \ref{restau=tau}]
Suppose that for all $n\in {\mathbb{N}}$ and all $t_0\in {\mathbb{C}}^*$ with 
$\underline{t}=(t_0,t_1,t_2,t_3)\in V$,
Proposition \ref{restau=tau} has been proved for  
$(n,t_0)$-residue pairs which have a $(n,t_0)$-resolution
of length $\leq s-1$ , where $s\geq 2$.

Fix arbitrary $n\in {\mathbb{N}}$ and
$t_0\in {\mathbb{C}}^*$ such that $\underline{t}=(t_0,t_1,t_2,t_3)
\in V$. The proposition is clear for $n=1$, so we may assume that $n>1$.
Let $(C,\mathfrak{C})$ be a 
$(n,t_0)$-residue pair 
with a $(n,t_0)$-resolution $(C_0,\ldots,C_s)$ of length $s$.
It suffices to prove \eqref{eindres} for the $(n,t_0)$-residue pair
$(C,{\mathfrak{C}})$.
We write $\Omega^{(l)}$ and $\Omega_{(l)}$ for the domains associated
with the $(n,t_0)$-residue pairs $(C_l,C_{l-1})$ and $(C_l,{\mathfrak{C}})$
respectively ($l\in\lbrace 1,\ldots,s\rbrace$). Note that
$\Omega_{(1)}\subset \Omega_{(2)}\subset\ldots\subset\Omega_{(s)}=\Omega$
where $\Omega$ is the domain associated with
the $(n,t_0)$-residue pair $(C,{\mathfrak{C}})$.
By \eqref{step2r} and \eqref{D1}, we have
\begin{equation}\label{step3n}
\begin{split}
\frac{1}{(2\pi i)^n}\iint\limits_{z\in C^n}f(z)\Delta(z;t_0)\frac{dz}{z}
= &\frac{1}{(2\pi i)^n}\iint\limits_{z\in (C_{s-1})^n}
f(z)\Delta(z;t_0)\frac{dz}{z}\\
+
\frac{2n}{(2\pi i)^{n-1}}&\sum_{i\in I_s}\,\,
\iint\limits_{z\in (C_{s-1})^{n-1}}
f_i(z)\Delta(z;tt_0q^i)\frac{dz}{z}\\
\end{split}
\end{equation}
for $f\in {\mathcal{O}}_W(\Omega^n)$, where
\begin{equation}
I_s:=\lbrace i\in {\mathbb{N}}_0 \, | \, 
r_s(\alpha_0^+)>|t_0q^i|>r_{s-1}(\alpha_0^+)
\rbrace
\end{equation}
and
$f_i(z):=f(t_0q^i,z)h(t_0q^i,z;t_0)$ with $h$ given by \eqref{hn}.
We will apply the induction hypotheses on all the terms in the right hand side
of \eqref{step3n}. For the integral over $\bigl(C_{s-1}\bigr)^n$
note that $(C_0,\ldots,C_{s-1})$ is a $(n,t_0)$-resolution of length $s-1$ for
the $(n,t_0)$-residue pair $(C_{s-1},{\mathfrak{C}})$. 
Hence, by the induction hypotheses,
\begin{equation}\label{ex1}
\begin{split}
\frac{1}{(2\pi i)^n}&\iint\limits_{z\in (C_{s-1})^n}
f(z)\Delta(z;t_0)\frac{dz}{z}=
\frac{1}{(2\pi i)^n}\iint\limits_{z\in {\mathfrak{C}}^n}
f(z)\Delta(z;t_0)\frac{dz}{z}\\
&+\sum_{r=1}^{n}\frac{2^r\bigl(n-r+1\bigr)_r}{(2\pi i)^{n-r}}\sum_{\omega\in 
D(r;C_{s-1},{\mathfrak{C}};t_0)}\,\,\,\iint\limits_{z\in{\mathfrak{C}}^{n-r}}
f(\omega,z)\Delta_r(\omega,z;t_0)\frac{dz}{z}\\
\end{split}
\end{equation}
for all $f\in {\mathcal{O}}_W(\Omega^n)$.

Now fix an $i\in I_s$. 
We have seen in the proof of Lemma \ref{restau=tau1} that
$h(t_0q^i,.;t_0)\in {\mathcal{O}}_W((\Omega^{(s)})^{n-1})$. 
In fact it follows from the proof that
$h(t_0q^i,.;t_0)\in {\mathcal{O}}_W(\Omega^{n-1})$.
In particular we have $f_i\in {\mathcal{O}}_W((\Omega_{(s-1)})^{n-1})$ 
for $f\in {\mathcal{O}}_W(\Omega^n)$. Furthermore we have that
$(tt_0t^i,t_1,t_2,t_3)\in V$ since $\hbox{arg}(tt_0t^i)=\hbox{arg}(t_0)$
and that $(C_0,\ldots,C_{s-1})$ is a $(n-1,tt_0q^i)$-resolution
of length $s-1$ 
for the $(n-1,tt_0q^i)$-residue pair $(C_{s-1},{\mathfrak{C}})$.
So we can apply the induction hypotheses for  all
the terms in the second line of \eqref{step3n}, and we obtain
\begin{equation}\label{ex2}
\begin{split}
&\frac{2n}{(2\pi i)^{n-1}}
\iint\limits_{z\in (C_{s-1})^{n-1}}f_i(z)\Delta(z;tt_0q^i)\frac{dz}{z}
=\frac{2n}{(2\pi i)^{n-1}}\iint\limits_{z\in {\mathfrak{C}}^{n-1}}f_i(z)
\Delta(z;tt_0q^i)\frac{dz}{z}\\
&+\sum_{r=2}^{n}\frac{2^{r}\bigl(n-r+1\bigr)_{r}}{(2\pi i)^{n-r}}
\sum_{\omega\in D(r-1;C_{s-1},{\mathfrak{C}};tt_0q^i)}\,\,\,\iint\limits_{z\in
{\mathfrak{C}}^{n-r}}f_i(\omega,z)\Delta_{r-1}(\omega,z;tt_0q^i)\frac{dz}{z}\\
\end{split}
\end{equation}
for $f\in {\mathcal{O}}_W(\Omega^n)$ and $i\in I_s$.

Substitution of  \eqref{ex1} and 
\eqref{ex2} in the right hand side of \eqref{step3n}
completes the proof of \eqref{eindres}, since  
\begin{equation}
\begin{split}
D(1;C,{\mathfrak{C}};t_0)&=D(1;C_{s-1},{\mathfrak{C}};t_0)\cup\lbrace
t_0q^i\rbrace_{i\in I_s},\nonumber\\
D(r;C,{\mathfrak{C}};t_0)&=
D(r;C_{s-1},{\mathfrak{C}};t_0)\cup\bigcup_{i\in I_s}
\big\lbrace (t_0q^i,\omega) \, | \, \omega\in
D(r-1;C_{s-1},{\mathfrak{C}};tt_0q^i)\big\rbrace\nonumber\\
\end{split}
\end{equation}
disjoint unions ($r\in\lbrace 2,\ldots,n\rbrace$) and
\begin{equation}\label{endform}
\begin{split}
f_i(z)\Delta(z;tt_0q^i)&=f(t_0q^i,z)\Delta_1(t_0q^i,z;t_0),\\
f_i(\omega,z)\Delta_{r-1}(\omega,z;tt_0q^i)&=
f(t_0q^i,\omega,z)\Delta_r(t_0q^i,\omega,z;t_0)\\
\end{split}
\end{equation}
for $i\in I_s$, $r\in\lbrace 2,\ldots,n\rbrace$ and
$\omega\in D(r-1;C_{s-1},{\mathfrak{C}};tt_0q^i)$
(recall the notational convention which implies that the $\Delta(z;tt_0q^i)$
respectively $\Delta_{r-1}(\omega,z;tt_0q^i)$ in the left hand
side of \eqref{endform} is with respect to $n-1$ variables $z$
respectively $(\omega,z)$, 
while the $\Delta_1(t_0q^i,z;t_0)$ respectively 
$\Delta_r(t_0q^i,\omega,z;t_0)$ in the right hand
side of \eqref{endform} is with respect to $n$ variables
$(t_0q^i,z)$ respectively $(t_0q^i,\omega,z)$).
\end{proof}   

\section{Limit transition to $q$-Racah polynomials}
In this section we will show that the orthogonality relations
and norm evaluations for the 
multivariable $q$-Racah polynomials can be obtained 
by applying the residue
calculus of the previous section to the results in Theorem \ref{concl}. 
The multivariable $q$-Racah 
polynomials are orthogonal with respect to a finite, 
discrete measure and were previously studied in \cite{vDS}.

For $\lambda\in P(r)$ we set
\begin{equation}\label{weightqRac}
\begin{split}
\Delta^{qR}\bigl(\rho q^{\lambda};&\underline{t};t\bigr):=
\prod_{i=1}^r\left(\frac{\bigl(q\rho_i^2;q\bigr)_{2\lambda_i}}
{\bigl(\rho_i^2;q\bigr)_{2\lambda_i}(q^{-1}t_0t_1t_2t_3t^{2i-2})^{\lambda_i}}
\prod_{j=0}^3\frac{\bigl(t_j\rho_i;q\bigr)_{\lambda_i}}
{\bigl(qt_j^{-1}\rho_i;q\bigr)_{\lambda_i}}\right)\\
&.\prod_{1\leq k<l\leq r}
\frac{\bigl(q\rho_k\rho_l,t\rho_k\rho_l;q\bigr)_{\lambda_k+\lambda_l}
\bigl(q\rho_k^{-1}\rho_l,t\rho_k^{-1}\rho_l;q\bigr)_{\lambda_l-\lambda_k}}
{\bigl(t^{-1}q\rho_k\rho_l,\rho_k\rho_l;q\bigr)_{\lambda_k+\lambda_l}
\bigl(t^{-1}q\rho_k^{-1}\rho_l,
\rho_k^{-1}\rho_l;q\bigr)_{\lambda_l-\lambda_k}}\\
\end{split}
\end{equation}
where $\rho_i:=t_0t^{i-1}$ and we set for $r\in {\mathbb{N}}_0$,
\begin{equation}\label{Cr}
\begin{split}
K_r(\underline{t};t):=&\prod_{i=1}^r\frac{\bigl(\rho_i^{-2};q\bigr)_{\infty}}
{\bigl(q,\rho_it_1,\rho_i^{-1}t_1,\rho_it_2,\rho_i^{-1}t_2,
\rho_it_3,\rho_i^{-1}t_3;q\bigr)_{\infty}}\\
&.\prod_{1\leq k<l\leq
r}\bigl(\rho_k^{-1}\rho_l,\rho_k^{-1}\rho_l^{-1};q\bigr)_{\tau}\\
=& \prod_{i=1}^r\frac{\bigl(\rho_i^{-2},t,t_0^{-2}t^{2-i-r};q\bigr)_{\infty}}
{\bigl(q,\rho_it_1,
\rho_i^{-1}t_1,\rho_it_2,\rho_i^{-1}t_2,\rho_it_3,\rho_i^{-1}t_3,
t^i,t_0^{-2}t^{2-2i};q\bigr)_{\infty}}
\end{split}
\end{equation}
where $t=q^{\tau}$. The discrete weights 
$\Delta^{(d)}(\rho q^{\lambda};\underline{t};t)$ \eqref{discreteweights}
can now be rewritten as follows.
\begin{Prop}\label{rewrittenqRacah} 
For $\lambda\in P(r)$ we have
\[\Delta^{(d)}(\rho q^{\lambda};\underline{t};t)=K_r(\underline{t};t)
\Delta^{qR}(\rho q^{\lambda};\underline{t};t),\]
where $\rho_i=t_0t^{i-1}$.
\end{Prop}
\begin{proof} 
We rewrite the discrete weight $\Delta^{(d)}(\rho q^{\lambda};\underline{t};t)$
using the following method. If we meet in the explicit expression for 
the discrete weight $\Delta^{(d)}(\rho q^{\lambda};\underline{t};t)$ 
a $q$-shifted factorial of the form $\bigl(aq^m;q\bigr)_b$ 
with $m\in {\mathbb{N}}_0$ depending only on $\lambda$,
we first rewrite this as a quotient of infinite products using \eqref{qshift}.
Then we replace the factors of the form $\bigl(cq^m;q\bigr)_{\infty}$ 
by $\bigl(c;q\bigr)_{\infty}\bigl(c;q\bigr)_m^{-1}$ if $m\in {\mathbb{N}}_0$,
respectively by
$\bigl(c;q\bigr)_{\infty}(-c)^{-m}q^{-\binom{1-m}{2}}
\bigl(qc^{-1};q\bigr)_{-m}$ 
if $m\in -{\mathbb{N}}$. For the case $m\in -{\mathbb{N}}$, 
we used here the formula
\begin{equation}\label{inversion}
\bigl(q^{1-l}x;q\bigr)_l=(-x)^lq^{-\binom{l}{2}}\bigl(x^{-1};q\bigr)_l,
\quad (l\in {\mathbb{N}}).
\end{equation}
Using this method we obtain for $j\in \lbrace 1,\ldots,r\rbrace$,
\begin{equation}\label{onethird}
\begin{split}
\prod_{i=1}^{j-1}&\frac{1}{\bigl(\rho_i\rho_jq^{\lambda_{i-1}+\lambda_j},
\rho_i\rho_j^{-1}q^{\lambda_{i-1}-\lambda_j};
q\bigr)_{\lambda_i-\lambda_{i-1}}}\\
&=(-1)^{\lambda_{j-1}}t^{(j-1)\lambda_j}q^{\lambda_{j-1}-
\binom{\lambda_j-\lambda_{j-1}}{2}+\binom{\lambda_j}{2}}
\frac{\bigl(q;q\bigr)_{\lambda_j-\lambda_{j-1}}
\bigl(t_0\rho_j;q\bigr)_{\lambda_j}}
{\bigl(\rho_j^2;q\bigr)_{\lambda_{j}+
\lambda_{j-1}}\bigl(qt_0^{-1}\rho_j;q\bigr)_{
\lambda_j}}\\
&.\prod_{i=1}^{j-1}\frac{\bigl(t\rho_i\rho_j;q\bigr)_{\lambda_i+\lambda_j}
\bigl(q\rho_i^{-1}\rho_j;q\bigr)_{\lambda_j-\lambda_i}}
{\bigl(\rho_i\rho_j;q\bigr)_{\lambda_i+\lambda_j}
\bigl(qt^{-1}\rho_i^{-1}\rho_j;q\bigr)_{\lambda_j-\lambda_i}}
t^{\lambda_i-\lambda_j}.\\
\end{split}
\end{equation}
Using \eqref{onethird} and applying the same method to 
the explicit expression for the weight
$w_d$ \eqref{weightfunctiondisc} gives
\begin{equation}\label{twothird}
\begin{split}
w_d(\rho_jq^{\lambda_j};\rho_jq^{\lambda_{j-1}})
&\prod_{i=1}^{j-1}\frac{1}{\bigl(\rho_{i}\rho_jq^{\lambda_{i-1}+\lambda_j},
\rho_{i}\rho_j^{-1}q^{\lambda_{i-1}-\lambda_j};
q\bigl)_{\lambda_i-\lambda_{i-1}}}
\\
=&\frac{\bigl(\rho_j^{-2};q\bigr)_{\infty}}
{\bigl(q,\rho_jt_1,\rho_j^{-1}t_1,\rho_jt_2,\rho_j^{-1}t_2,\rho_jt_3,
\rho_j^{-1}t_3;q\bigr)_{\infty}}\\
&.\frac{\bigl(q\rho_j^2;q\bigr)_{2\lambda_j}}
{\bigl(\rho_j^2;q\bigr)_{2\lambda_j}
(q^{-1}t_0t_1t_2t_3)^{\lambda_j}}\prod_{k=0}^3
\frac{\bigl(t_k\rho_j;q\bigr)_{\lambda_j}}
{\bigl(qt_k^{-1}\rho_j;q\bigr)_{\lambda_j}}\\
&.\prod_{i=1}^{j-1}
\frac{\bigl(t\rho_i\rho_j;q\bigr)_{\lambda_i+\lambda_j}
\bigl(q\rho_i^{-1}\rho_j;q\bigr)_{\lambda_j-\lambda_i}}
{\bigl(\rho_i\rho_j;q\bigr)_{\lambda_i+\lambda_j}
\bigl(qt^{-1}\rho_i^{-1}\rho_j;q\bigr)_{\lambda_j-\lambda_i}}
t^{\lambda_i-\lambda_j}\\
\end{split}
\end{equation}
for $j\in \lbrace 1,\ldots,r\rbrace$.
Now again applying the above mentioned method, one obtains
\begin{equation}\label{threethird}
\begin{split}
\prod_{i=1}^{j-1}\bigl(\rho_i^{-1}\rho_j&q^{\lambda_{j}-\lambda_i},
\rho_i^{-1}\rho_j^{-1}q^{-\lambda_{i}-\lambda_j};q\bigr)_{\tau}\\
=& \prod_{i=1}^{j-1}\bigl(\rho_i^{-1}\rho_j,
\rho_i^{-1}\rho_j^{-1};q\bigr)_{\tau}
\frac{\bigl(t\rho_i^{-1}\rho_j;q\bigr)_{\lambda_j-\lambda_i}
\bigl(q\rho_i\rho_j;q\bigr)_{\lambda_i+\lambda_j}}
{\bigl(\rho_i^{-1}\rho_j;q\bigr)_{\lambda_j-\lambda_i}
\bigl(qt^{-1}\rho_i\rho_j;q\bigr)_{\lambda_i+\lambda_j}}
t^{-\lambda_i-\lambda_j}\\
\end{split}
\end{equation}
for $j\in \lbrace 1,\ldots,r\rbrace$, where $t=q^{\tau}$.
Now the proposition follows by 
multiplying \eqref{twothird} and \eqref{threethird} and taking the product
over $j\in\lbrace 1,\ldots,r\rbrace$.
\end{proof}

In the next theorem we 
give orthogonality relations for the Askey-Wilson polynomials
when the parameters $\underline{t}$ satisfy the truncation condition
$t^{n-1}t_0t_3=q^{-N}$. 
We will formulate the theorem with the parameters considered as 
indeterminates. Set $F:={\mathbb{C}}(\underline{t},t)$, 
${\mathbb{F}}:={\mathbb{C}}(t_0,t_1,t_2,t)$ and
$\underline{t}_N:=(t_0,t_1,t_2,t^{1-n}t_0^{-1}q^{-N})$. 
Let $A_F^W$ respectively
$A_{\mathbb{F}}^W$ be the 
algebra of $W$-invariant Laurent polynomials over the field
$F$ respectively ${\mathbb{F}}$. We define the {\it Askey-Wilson polynomial} 
$P_{\lambda}(.;\underline{t};t)\in A_F^W$ over the field $F$ as
\begin{equation}
P_{\lambda}(.;\underline{t};t):=\left(\prod_{\mu<\lambda}
\frac{D_{\underline{t},t}-E_{\mu}(\underline{t},t)}
{E_{\lambda}(\underline{t},t)-E_{\mu}(\underline{t};t)}\right)m_{\lambda}.
\end{equation}
By specializing the parameters $(\underline{t},t)$ to values in $V\times
(0,1)$, we obtain the monic 
multivariable Askey-Wilson polynomial
$P_{\lambda}(.;\underline{t};t)\in A^W$ 
of degree $\lambda$ as defined in Theorem
\ref{concl} (see Definition \ref{AWdefinitie}).

Note that $P_{\lambda}(.;\underline{t}_N;t)\in A_{\mathbb{F}}^W$
is well defined since the eigenvalues
$\lbrace E(\lambda;\underline{t}_N;t)\rbrace_{\lambda\in\Lambda}$ \eqref{AWb}  
are mutually different as elements in ${\mathbb{C}}[t_0,t_1,t_2,t]$. 
\begin{Def}\label{qRdefinitie}
We call $\lbrace P_{\lambda}(.;\underline{t}_N;t)
\rbrace_{\lambda\in\Lambda_N}\subset A_{\mathbb{F}}^W$ 
with $\Lambda_N:=\lbrace \lambda\in\Lambda \, | \, \lambda_1\leq N
\rbrace$ the family of multivariable (BC type) $q$-Racah polynomials.
\end{Def}
Let $\lambda\in P(n)$, 
then the weight $\Delta^{qR}(\rho q^{\lambda};\underline{t}_N;t)
\in {\mathbb{F}}$ is well defined ($\Delta^{qR}$ given by \eqref{weightqRac})
and it is non zero if and only if $\lambda_n\leq N$ due to the factor
$\bigl(\rho_nt_3;q\bigr)_{\lambda_n}$ in the numerator of 
$\Delta^{qR}(\rho q^{\lambda};\underline{t};t)$. So 
the bilinear form 
\begin{equation}\label{bilinearformqRac}
\langle f,g\rangle_{qR,\underline{t}_N,t}:=
\sum_{\lambda\in P(n)}f(\rho q^{\lambda})g(\rho q^{\lambda})
\Delta^{qR}(\rho q^{\lambda};\underline{t}_N;t),\quad 
f,g\in A^W_{\mathbb{F}},
\end{equation} 
takes its values in the field ${\mathbb{F}}$.
Let ${\mathcal{N}}^{qR}(\lambda;\underline{t};t)$ for $\lambda\in\Lambda$
be given by
\begin{equation}\label{NqR}
{\mathcal{N}}^{qR}(\lambda;\underline{t};t):=
\frac{{\mathcal{N}}(\lambda;\underline{t};t)}{K_n(\underline{t};t)2^nn!}
\end{equation}
where ${\mathcal{N}}(\lambda)$ \eqref{normexp} is the expression for
the quadratic norms of the Askey-Wilson polynomial $P_{\lambda}$.
Substitution of the explicit expressions for ${\mathcal{N}}(\lambda)$ and
$K_n$ in \eqref{NqR} yields that 
${\mathcal{N}}^{qR}(\lambda;\underline{t}_N;t)\in {\mathbb{F}}$ and that
${\mathcal{N}}^{qR}(\lambda;\underline{t}_N;t)$ is non zero 
if and only if $\lambda\in\Lambda_N$. We have now the following theorem.
\begin{Thm}\label{qRacahorth}
Let $N\in {\mathbb{N}}$. The 
$q$-Racah polynomials 
$P_{\lambda}(.;\underline{t}_N;t)$ ($\lambda\in\Lambda_N$)
are orthogonal with respect to $\langle .,. \rangle_{qR,\underline{t}_N,t}$
and the quadratic norms are given by
\begin{equation}
\langle P_{\lambda}(.;\underline{t}_N;t),P_{\lambda}(.;\underline{t}_N;t)
\rangle_{qR,\underline{t}_N,t}=
{\mathcal{N}}^{qR}(\lambda;\underline{t}_N;t)\qquad
(\lambda\in\Lambda_N).
\end{equation}
\end{Thm}
\begin{proof}
Let $\tilde{V}\subset ({\mathbb{C}}^*)^4$ 
be the set of parameters $\underline{t}\in
({\mathbb{C}}^*)^4$ 
for which $t_0t_1t_2t_3\in {\mathbb{C}}\setminus{\mathbb{R}}$.
Note that there exists an open dense subset $I_N\subset (0,1)$ such that
$E_{\lambda}(\underline{t};t)\not=E_{\mu}(\underline{t};t)$ for all
$\underline{t}\in \tilde{V}$, $t\in I_N$ and all 
$\lambda,\mu\in \Lambda_N$ with 
$\lambda\not=\mu$.

Fix $t_0,t_1,t_2\in {\mathbb{C}}^*$ such that $\#\lbrace
\hbox{arg}(t_i),\hbox{arg}(t_i^{-1}) \, | \, i=0,1,2\rbrace=6$ 
and $t\in I_N$. Then
$\underline{t}_N\in \tilde{V}$
and there exists a sequence $\lbrace t_{3,i}\rbrace_{i\in
{\mathbb{N}}_0}\subset {\mathbb{C}}^*$ converging to $t^{1-n}t_0^{-1}q^{-N}$ 
such that $\underline{t}_i:=(t_0,t_1,t_2,t_{3,i})\in V\cap \tilde{V}$ 
for all $i$ ($V$ given in Definition \ref{WAW}). 
By considering a subsequence if necessary, we
may assume that there exist $(n,t_0)$-residue pairs $(C_i,{\mathfrak{C}})$
where $C_i$ is a $\underline{t}_i$-contour 
and where ${\mathfrak{C}}$ is a deformed
circle such that the sequences $\lbrace t_1q^j,t_2q^j,
t_{3,i}q^j\rbrace_{j\in {\mathbb{N}}_0}$ 
are in the interior of ${\mathfrak{C}}$
for all $i$ and such that 
$t^{n-1}t_0q^{N}$ is in the exterior of ${\mathfrak{C}}$. 
Then we obtain from Theorem \ref{concl},  
Proposition \ref{restau=tau} and Proposition \ref{rewrittenqRacah} that
\begin{equation}\label{bijna!}
\begin{split}
\frac{{\mathcal{N}}(\lambda;\underline{t}_i;t)}{K_n(\underline{t}_i;t)}
&\delta_{\lambda,\mu}=
\frac{1}{(2\pi i)^n}\iint\limits_{z\in {\mathfrak{C}}^n}
\bigl(P_{\lambda}P_{\mu}\bigr)(z;\underline{t}_i;t)
\frac{\Delta(z;\underline{t}_i;t)}{K_n(\underline{t}_i;t)}\frac{dz}{z}\\
+\sum_{r=1}^n&\frac{2^r\bigl(n-r+1\bigr)_r}{(2\pi i)^{n-r}}
\sum_{\omega\in D(r)}\,\,\iint\limits_{z\in {\mathfrak{C}}^{n-r}}
\bigl(P_{\lambda}P_{\mu}\bigr)(\omega,z;\underline{t}_i;t)
\frac{\Delta_r(\omega,z;\underline{t}_i;t)}{K_n(\underline{t}_i;t)}
\frac{dz}{z}\\
\end{split}
\end{equation}
where $\delta_{\lambda,\mu}$ is the Kronecker-delta and 
$D(r)=D(r;C_i,{\mathfrak{C}};t_0;t)$ \eqref{Dparameter} 
(which is independent of $i$). 
By \eqref{weightalgt} and Proposition \ref{rewrittenqRacah} we have
\begin{equation}\label{substitute}
\Delta_r(\omega,z;\underline{t}_i;t)=
K_r(\underline{t}_i;t)
\Delta^{qR}(\omega;\underline{t}_i;t)\Delta(z;\underline{t}_i;t)
\delta_c(\omega;z).
\end{equation}
After substitution of \eqref{substitute} in the right hand side of 
\eqref{bijna!} for all $r$, it follows 
from the Bounded Convergence 
Theorem that we may take the limit $i\rightarrow\infty$
within the integrals in the right hand side of \eqref{bijna!}.  
Only the completely discrete part survives the limit $i\rightarrow\infty$
in the equality \eqref{bijna!} since
\[
\lim_{i\rightarrow\infty}
\frac{K_r(\underline{t}_i;t)}{K_n(\underline{t}_i;t)}=0,
\qquad 0\leq r<n
\]
by the factor $\bigl(\rho_nt_3;q\bigr)_{\infty}$
in the denominator of $K_n(\underline{t};t)$.
The theorem follows now for the specialized parameter values $t_0,t_1,t_3,t$ 
from the fact that
\[\lbrace \rho q^{\lambda} \, | \, \lambda\in P(n), \,\, \lambda_n\leq N\rbrace
\subset D(n)\]
and the fact that $\Delta^{qR}(\rho q^{\lambda}; \underline{t}_N;t)=0$ for 
$\lambda\in P(n)$ with $\lambda_n>N$.
It is now clear that the theorem also holds over the field ${\mathbb{F}}$.
\end{proof}
The constant term identity can be simplified as follows.
\begin{Cor}
For $N\in {\mathbb{N}}$ we have the summation formula
\begin{equation}\label{norm0qR}
\langle 1,1\rangle_{qR,\underline{t}_N,t}=
\prod_{i=1}^n\frac{\bigl(qt_0^2t^{2n-i-1},qt_1^{-1}t_2^{-1}t^{i-n};q\bigr)_N}
{\bigl(qt_0t_1^{-1}t^{n-i},qt_0t_2^{-1}t^{n-i};q\bigr)_N}.
\end{equation}
\end{Cor}
\begin{proof}
First note that by 
\eqref{norm0}, \eqref{Cr} and \eqref{NqR} we have the explicit formula
\begin{equation}\label{normAWdiscrete}
{\mathcal{N}}^{qR}(0;\underline{t};t)=
\prod_{i=1}^n
\frac{\bigl(t_0t_1t_2t_3t^{2n-i-1},t_0^{-1}t_1t^{1-i}, t_0^{-1}t_2t^{1-i},
t_0^{-1}t_3t^{1-i};q\bigr)_{\infty}}
{\bigl(t_0^{-2}t^{1+i-2n},
t_1t_2t^{i-1},t_1t_3t^{i-1},t_2t_3t^{i-1};q\bigr)_{\infty}}.
\end{equation}
Then \eqref{norm0qR} follows by 
substitution of $t_3=t_0^{-1}t^{1-n}q^{-N}$ in \eqref{normAWdiscrete}
and by applying formula \eqref{inversion} repeatedly
(see also \cite[section 3]{vD2}).
\end{proof} 
The second order $q$-difference operator
$D_{\underline{t}_N,t}$ 
\eqref{secondorderqdiff} diagonalizes the $q$-Racah polynomials
$\lbrace P_{\lambda}(.;\underline{t}_N;t)\rbrace_{\lambda\in\Lambda_N}$.
By Theorem \ref{qRacahorth}
we conclude that $D_{\underline{t}_N,t}$ is symmetric
with respect to $\langle .,. \rangle_{qR,\underline{t}_N,t}$.
In \cite{vDS} the symmetry of $D_{\underline{t}_N,t}$ 
was proved by direct calculations and the orthogonality relations for 
the multivariable $q$-Racah polynomials were proved using the symmetry of 
$D_{\underline{t}_N,t}$.  
Furthermore, in \cite{vDS} the quadratic norms of the $q$-Racah polynomials 
were expressed in terms of the quadratic norm 
of the unit polynomial by studying Pieri formulas
for the $q$-Racah polynomials. The constant term identity \eqref{norm0qR} was 
recently proved by van Diejen \cite[Theorem 3]{vD2} 
by truncating a multivariable analogue of Roger's ${_6}\phi_5$-series 
\cite[Theorem 2]{vD2}, 
which in turn is closely related to an Aomoto-Ito type sum 
(cf. \cite{Ao}, \cite{I}) for the non-reduced root system $BC_n$. The proofs of
the summation formulas in \cite{vD2} are based on
a multiple ${_6}\psi_6$ summation formula of Gustafson.
 
In the one-variable case it is known that 
the Askey-Wilson integral can be rewritten as an infinite sum of residues for 
some parameter region by shifting the contour over four infinite 
sequences of poles (see \cite[Theorem 2.1]{AW}).
More generally one can ask the question whether a completely discrete orthogonality
measure for the multivariable Askey-Wilson polynomials can be obtained
by pulling the $\underline{t}$-contours over certain infinite
sequence of poles in the orthogonality relations of 
the Askey-Wilson polynomials (Theorem \ref{concl}). 

Strong indications in that direction can be found in Gustafson's paper \cite{G1}  
and the recent paper of Tarasov and Varchenko \cite{TV} where 
contours in multidimensional integrals are shifted over infinite
sequences of poles in order to arrive at (purely discrete) multidimensional
Jackson integrals. Another strong indication is the fact that
the Macdonald polynomials are orthogonal with respect to Aomoto-Ito type (cf.
\cite{Ao}, \cite{I}) weight functions (see Cherednik \cite{C2}). 
Since the $B$, $C$ and $D$ type Macdonald polynomials can be
obtained from the Askey-Wilson polynomials by suitable specialization of the
parameters we thus have orthogonality relations for these subfamilies 
of the Askey-Wilson polynomials with respect to infinite discrete measures
(and the corresponding discrete weights are directly 
related to \eqref{weightqRac}, see \cite{vD2}).

In this paper we will not consider the above mentioned questions, but instead
look at the implications of the residue calculus for certain limit cases
of the Askey-Wilson polynomials (the big and little $q$-Jacobi polynomials).
In order to study these limit cases we first need to consider
the Askey-Wilson polynomials for yet another parameter domain. 
This will be the subject of the next section.

\section{Askey-Wilson polynomials with positive orthogonality measure}
In this section we will consider the Askey-Wilson polynomials for parameters
$\underline{t}$ in the following parameter domain.
\begin{Def}\label{domain}
Let $V_{AW}$ be the set of parameters $\underline{t}=(t_0,t_1,t_2,t_3)$ 
which satisfy the following conditions:\\
{\bf (1)} The parameters $t_0,t_1,t_2,t_3$ are real, or if complex, 
then they appear in conjugate pairs.\\
{\bf (2)} $t_kt_l
\notin {\mathbb{R}}_{\geq 1}$ for all $0\leq k<l\leq 3$.
\end{Def}
Note that parameters $\underline{t}\in V_{AW}$ 
satisfy the following properties:\\
{\bf (A)} $t_i\in {\mathbb{R}}$ if $|t_i|\geq 1$;\\
{\bf (B)} There are at most two parameters with modulus $\geq 1$.
If there are two, then one is positive and the other is negative.

We will show that the multivariable
Askey-Wilson polynomials are orthogonal with respect to a
positive, partially discrete orthogonality measure 
when $t\in (0,1)$ and $\underline{t}\in V_{AW}$.
We proceed as follows. 
We apply the residue calculus of section 3 to shift the contour $C^n$ in
the integral
\begin{equation}\label{goal}
\frac{1}{(2\pi i)^n}
\iint\limits_{z\in C^n}P_{\lambda}(z)P_{\mu}(z)\Delta(z)\frac{dz}{z}
\end{equation}
to the $n$-torus $T^n$ for a specific parameter domain $V_0\subset V$
(here $C$ is a $\underline{t}$-contour (Definition \ref{contour}) and $V$
is the parameter domain given in Definition \ref{WAW}). 
We then obtain a 
partially discrete orthogonality measure which turns out to be 
well defined and positive for parameter values $\underline{t}\in V_{AW}$.
Orthogonality relations for parameter values $\underline{t}\in V_{AW}$ 
with respect to this positive, partially
discrete orthogonality measure 
can then be derived by suitable continuity arguments.

The parameter domain $V_0$ is defined as follows.
\begin{Def}\label{UAW0}
Let $V_0$ be the set of parameters $\underline{t}\in V$
for which\\
{\it (i)} at most two parameters have modulus $>1$;\\
{\it (ii)} $t_it^jq^p\notin T$ 
for $i\in\lbrace 0,\ldots,3\rbrace$, $j\in \lbrace
-1,\ldots,n-1\rbrace$ and $p\in {\mathbb{Z}}$.
\end{Def}
Fix $t\in (0,1)$, $\underline{t}\in V_0$ and 
$0\leq i\not=j\leq 3$ 
such that $|t_k|<1$ for $k\not=i,j$.
We write $\rho^{(i)}_p:=t^{p-1}t_i$ respectively $\rho^{(j)}_p:=t^{p-1}t_j$ 
for $p\in {\mathbb{Z}}$. 
Define for $r\in {\mathbb{N}}$ a finite discrete set
$D_{i}(r)=D_{i}(r;\underline{t};t)\subset {\mathbb{C}}^r$ by
\begin{equation}\label{Dir}
D_{i}(r):=\lbrace \rho^{(i)} q^{\mu} \,\, | \,\, 
\mu\in P(r),\,\, |\rho^{(i)}_rq^{\mu_r}|>1 \rbrace 
\end{equation}
and similarly for $D_{j}(r)$ (here $P(r)$ is given by \eqref{Pr}). 
We have used here the notation
$\rho^{(i)} q^{\mu}=(\rho^{(i)}_1q^{\mu_1},\ldots,\rho^{(i)}_rq^{\mu_r})$
for $\mu\in P(r)$.
Note that $D_{i}(r)=\emptyset$ if $|t_i|<1$. 
Furthermore, we write $F(r)=F(r;\underline{t};t)\subset {\mathbb{C}}^r$ 
for the disjoint union
\begin{equation}\label{Fr}
F(r):= \bigcup_{\stackrel{l+m=r}{l,m\in
{\mathbb{N}}_0}}D_{i}(l)\times D_{j}(m)\qquad (r=1,\ldots,n).
\end{equation}
(We use here the convention that $D_i(l)\times D_j(m)=\emptyset$ if
$l>0$ and $D_i(l)=\emptyset$ or if $m>0$ and $D_j(m)=\emptyset$, and that
$D_i(0)\times D_j(m)=D_j(m)$, $D_i(l)\times D_j(0)=D_i(l)$.)
Let $\omega\in F(r)$ and $z\in T^{n-r}$ and set
\begin{equation}
d\nu_r^{AW}(\omega,z;\underline{t};t):=\Delta_r^{AW}(\omega,z;\underline{t};t)
\frac{dz}{z}
\end{equation}
with weight function $\Delta_r^{AW}(\omega,z)$ for $\omega=(\omega^{(i)},
\omega^{(j)})$ with $\omega^{(i)}\in D_{i}(l)$ and
$\omega^{(j)}\in D_{j}(m)$ given by
\begin{equation}\label{DeltarAW}
\begin{split}
\Delta_r^{AW}(\omega^{(i)},\omega^{(j)},z;\underline{t};t):=
\Delta^{(d)}(\omega^{(i)};t_i)&\Delta^{(d)}(\omega^{(j)};t_j)
\Delta(z;\underline{t};t)\\
&.\delta_c(\omega^{(i)};\omega^{(j)},z)\delta_c(\omega^{(j)};z)\\
\end{split}
\end{equation}
where $\Delta^{(d)}$ is given by \eqref{discreteweights} and $\delta_c$
is given by \eqref{continuousinteraction}.
In the special case 
that $l=0$ respectively $m=0$, \eqref{DeltarAW} simplifies to
\begin{equation}\label{DeltarAWl=0}
\Delta_r^{AW}(\omega^{(j)},z;\underline{t};t)=
\Delta^{(d)}(\omega^{(j)};t_j)
\Delta(z;\underline{t};t)\delta_c(\omega^{(j)};z)=\Delta_r(\omega^{(j)},z;t_j)
\end{equation}
respectively
\begin{equation}\label{DeltarAWm=0}
\Delta_r^{AW}(\omega^{(i)},z;\underline{t};t)=
\Delta^{(d)}(\omega^{(i)};t_i)
\Delta(z;\underline{t};t)\delta_c(\omega^{(i)};z)=\Delta_r(\omega^{(i)},z;t_i).
\end{equation}
where $\Delta_r$ is given by \eqref{weightalgt}.
We obtain from the residue calculus of section 3
the following lemma.
\begin{Lem}\label{resn}
Let $t\in (0,1)$ and $\underline{t}\in V_0$.
Let $C$ be a $\underline{t}$-contour and $f\in A^W$. Then,
\begin{equation}\label{AWpos}
\begin{split}
\frac{1}{(2\pi i)^n}\iint\limits_{z\in C^n}&f(z)d\nu(z)=
\frac{1}{(2\pi i)^n}\iint\limits_{z\in T^n}f(z)d\nu(z)\\
+&\sum_{r=1}^n\frac{2^r\bigl(n-r+1\bigr)_r}{(2\pi i)^{n-r}}
\sum_{\omega\in F(r)}
\,\,\iint\limits_{z\in T^{n-r}}f(\omega,z)d\nu_r^{AW}(\omega,z).\\
\end{split}
\end{equation}
\end{Lem}
\begin{proof}
If  $|t_k|<1$ for
all $k$ then we only have
the completely continuous measure $d\nu$ on $T^n$ in the
right hand side of \eqref{AWpos} since $F(r)=\emptyset$.
Since $T$ is a $\underline{t}$-contour in this case, the lemma follows
from Lemma \ref{independent}.

Suppose that at most one parameter has modulus $> 1$. By the symmetry of
$d\nu(z;\underline{t};t)$ in the four parameters $\underline{t}$, we may assume
that $|t_0|> 1$.
By Lemma \ref{independent}, we may assume that the $\underline{t}$-contour 
$C$ satisfies the additional conditions that 
$A^+:=\lbrace x\in [0,1] \, | \, r_C(x)>1\rbrace$ is an 
open interval and that 
$\alpha_0^+\in A^+$ 
but $\alpha_i^{\pm}\notin A^+$
for $i=1,2,3$ (here $r_C$ is as in Definition \ref{contour}, 
and $\alpha_i^+$ is given by \eqref{alpha}).
Then $(C,T)$ is a $(n,t_0)$-residue pair since $\underline{t}\in V_0$
(Definition \ref{UAW0}) and $D_0(l)=D(l;C,T;t_0)$ \eqref{Dparameter}
since $t_0$ is in the interior of $C$.
The lemma is then a direct consequence
of Proposition \ref{restau=tau} and \eqref{DeltarAWm=0}.

Suppose now that two parameters have moduli $>1$. Without loss of generality,
we may assume that 
$|t_0|>1$ and $|t_1|>1$ and that the $\underline{t}$-contour 
$C$ satisfies the additional condition
that
\[\lbrace x\in [0,1] \, | \, r_C(x)>1\rbrace=
A_0^+\cup A_1^+\]
disjoint union, with $A_i^+$ open intervals such that 
$\alpha_i^+\in A_i^+$ and 
$\alpha_j^{\pm}\notin A_i^+$ for $j\not=i$ and $i=0,1$.
Let $C':=\phi_{C'}([0,1])$ be the deformed circle with parametrization
$\phi_{C'}(x)=r_{C'}(x)e^{2\pi i x}$ given by
\[ r_{C'}(x):=r_C(x)\,\,\, (x\notin A_0^+\cup A_0^-),
\quad r_{C'}(x):=1 \,\,\, (x\in A_0^+\cup A_0^-),\]
where $A_0^-:=(1-\beta, 1-\alpha)$ when $A_0^+=(\alpha,\beta)$.
Then $(C,C')$ is a $(n,t_0)$-residue pair,
$(C',T)$ is a $(n,t_1)$-residue pair and $D_0(l)=D(l;C,C';t_0)$ respectively
$D_1(m)=D(m;C',T;t_1)$ since $t_0$ and $t_1$ are in the interior of $C$.
Write $\Omega'$ for the domain associated with $(C',T)$, then
$\delta_c(\omega^{(0)};.)\in {\mathcal{O}}_W\bigl((\Omega')^{n-l}\bigr)$
for $\omega^{(0)}\in D_{0}(l)$, hence 
the lemma follows by applying Proposition \ref{restau=tau}
first to the $(n,t_0)$-residue pair $(C,C')$, 
and then to the $(n,t_1)$-residue pair $(C',T)$.
\end{proof}
For $t=q^k$ with $k\in {\mathbb{N}}$, \eqref{AWpos} can be rewritten as
\begin{equation}\label{AWposnatural}
\begin{split}
\frac{1}{(2\pi i)^n}\iint\limits_{z\in C^n}f(z)d\nu(z)=
\sum_{r=0}^n&\frac{2^r\binom{n}{r}}{(2\pi i)^{n-r}}
\sum_{e_1,\ldots,e_r}
\sum_{\stackrel{z_i\in\lbrace e_i,\ldots,e_iq^{N_{e_i}}\rbrace}{
i=1,\ldots,r}}\\
\iint\limits_{(z_{r+1},\ldots,z_n)\in T^{n-r}}&f(z)\delta(z;q^k)
\prod_{i=1}^rw_d(z_i;e_i)\prod_{j=r+1}^nw_c(z_j)\frac{dz_j}{z_j}\\
\end{split}
\end{equation}
for $f\in A^W$, where the sum is over $e_i\in\lbrace t_j \, | \,
|t_j|>1 \rbrace$ and $N_{e_i}$ is the largest positive integer
such that $|e_iq^{N_{e_i}}|>1$. 
This follows from the fact that $\delta(z;q^k)=0$
if $z_i=q^lz_j$ 
for some $i\not=j$ and some $l\in\lbrace 0,\ldots,k-1\rbrace$ and 
from the fact that
\[w_d(x;e_iq^l)=\bigl(e_ix,e_ix^{-1};q\bigr)_lw_d(x;e_i),\quad
(x=e_iq^{l+m}, m\in {\mathbb{N}}_0).\]
\begin{Rem}
Note that $\delta(.;t)\in A^W$ when $t=q^k$ with $k\in {\mathbb{N}}$.
In particular, we will only encounter residues 
of the factor $\prod_{i=1}^nw_c(z_i)$ when deforming $C^n$ for
parameter $t=q^k$ ($k\in
{\mathbb{N}}$) to the torus $T^n$
in the left hand side of \eqref{AWposnatural}.
Consequently \eqref{AWposnatural} can also be proved by induction on $n$
using the residue calculus 
for the one-variable weight functions $w_c$ and using
the $W$-invariance of the integrand in the left hand side of
\eqref{AWposnatural}. 
\end{Rem}
We define bilinear forms $\langle .,. \rangle_{r,\underline{t},t}$
on $A^W$ for $r\in\lbrace 0,\ldots,n\rbrace$,
$\underline{t}\in V_0$ and $t\in (0,1)$ by
\begin{equation}\label{symmvormAWr}
\begin{split}
\langle f,g\rangle_0&:=\iint\limits_{z\in T^n}
f(z)g(z)d\nu(z),\\
\langle f,g\rangle_{r}&:=
\sum_{\omega\in F(r)}\,\,\iint\limits_{z\in
T^{n-r}}f(\omega,z)g(\omega,z)d\nu_r^{AW}(\omega,z), \quad r\in\lbrace
1,\ldots,n\rbrace\\
\end{split}
\end{equation}
for $f,g\in A^W$ and we set 
\begin{equation}\label{symmvormAW}
\langle f,g\rangle_{\underline{t},t}:=\sum_{r=0}^n
\frac{2^r\bigl(n-r+1\bigr)_r}{(2\pi i)^{n-r}}
\langle f,g\rangle_{r,\underline{t},t},
\quad f,g\in A^W.
\end{equation}
In the following lemma we consider the symmetric bilinear form $\langle .,.
\rangle_{\underline{t},t}$ 
for parameter values $(\underline{t},t)\in V_{AW}\times
(0,1)$. 
\begin{Lem}\label{techlem}
Let $t\in (0,1)$ and $\underline{t}\in V_{AW}$.\\
{\bf (i)} The bilinear form $\langle .,. \rangle_{\underline{t},t}$
is well defined;\\
{\bf (ii)} The weight function $\Delta(z;\underline{t};t)$ respectively
$\Delta_r^{AW}(\omega,z;\underline{t};t)$ is positive for $z\in T^n$
respectively $(\omega,z)\in F(r)\times T^{n-r}$ ($r=1,\ldots,n$).
\end{Lem}
\begin{proof}
The discrete weights $w_d$ \eqref{weightfunctiondisc} appearing as factors of
the weight function $\Delta_r^{AW}(\omega,z)$ for $r>0$ are well
defined and strict positive. Indeed if $t_0t_1t_2t_3=0$,
then the factors $\left(t_iq/t_j;q\right)_kt_j^k$ in the denominator of
$w_d$ \eqref{weightfunctiondisc} should be read 
as $\prod_{l=0}^{k-1}\left(t_j-t_iq^{l+1}\right)$.
The factor $\delta(z;t)=|\delta_+(z;t)|^2$ is also well defined and positive
for $z\in T^{n-r}$.

Without loss of generality we may assume that $|t_2|,|t_3|<1$.
Fix $\omega=(\vartheta,\upsilon)\in F(r)$ with $\vartheta\in D_{0}(l)$
and $\upsilon\in D_{1}(m)$\, ($r=l+m$).
The factor $\delta_d(\vartheta)$ respectively $\delta_d(\upsilon)$
\eqref{discreteinteraction} 
appearing in the discrete weights 
$\Delta^{(d)}(\vartheta;t_0)$ respectively $\Delta^{(d)}(\upsilon;t_1)$ 
when $l>0$ respectively $m>0$
is well defined and strict positive. Indeed, if $\vartheta\in D_{0}(l)$
and $l>0$, then we have $|t_0|>1$, hence
$t_0\in {\mathbb{R}}$. Then $\delta_d(\vartheta)>0$ follows easily
from the definition of the set $D_{0}(l)$ \eqref{Dir}.

Remains to show that
$h(z):=\bigl(\prod_{l=1}^{n-r}w_c(z_l;\underline{t})\bigr)\delta_c(\vartheta;
\upsilon,z)\delta_c(\upsilon;z)$ 
is well defined and positive for $z\in T^{n-r}$. 
Let us check the case that both $t_0$ and $t_1$ have moduli $\geq 1$,
and that $t_0$ is positive real and $t_1$ negative real (see property {\bf (B)}
for parameters $\underline{t}\in V_{AW}$). The case that at most one parameter
has modulus $\geq 1$ will then also be clear.

Rewrite the 
factor $\bigl(x^2,x^{-2};q\bigr)_{\infty}$ appearing in the numerator
of $w_c(x;\underline{t})$ as
\[\bigl(x^2,x^{-2};q\bigr)_{\infty}=\bigl(x,-x,x^{-1},-x^{-1};q\bigr)_{\infty}
(qx^2,qx^{-2};q^2\bigr)_{\infty}
\]
then it is sufficient to check that the factors of the form
\begin{equation}
\begin{split}
h_0(x)&:=
\frac{\bigl(x,x^{-1};q\bigr)_{\infty}}{\bigl(t_0x,t_0x^{-1};q\bigr)_{\infty}}
\prod_{k=1}^{l}
\bigl(\vartheta_{k}x,\vartheta_{k}x^{-1},\vartheta_{k}^{-1}x,
\vartheta_{k}^{-1}x^{-1};q\bigr)_{\tau},\nonumber\\
h_1(x)&:=\frac{\bigl(-x,-x^{-1};q\bigr)_{\infty}}
{\bigl(t_1x,t_1x^{-1};q\bigr)_{\infty}}
\prod_{k=1}^{m}
\bigl(\upsilon_{k}x,\upsilon_{k}x^{-1},\upsilon_{k}^{-1}x,
\upsilon_{k}^{-1}x^{-1};q\bigr)_{\tau}\nonumber
\end{split}
\end{equation}
$(l,m\in\lbrace 0,\ldots,n-1\rbrace)$
are well defined and positive for $x\in T$ (here $t=q^{\tau}$).
Indeed, the remaining factors of $w_c(x;\underline{t})$ are easily seen
to be well defined and 
positive since $|t_2|,|t_3|<1$ and $t_2,t_3$ are both real or 
are a conjugate pair, while the remaining factors
\[
\prod_{\epsilon_i,\epsilon_j=\pm 1}\bigl(\vartheta_i^{\epsilon_i}
\upsilon_j^{\epsilon_j};q\bigr)_{\tau}\quad
(i\in \lbrace 1,\ldots, l\rbrace, \, j\in\lbrace 1,\ldots,m\rbrace)
\]
of $\delta_c(\vartheta;\upsilon,z)$ are well defined and positive
since $t_0$ is positive real and $t_1$ is negative real.
Now let $\lambda\in P(l)$ such that 
$\vartheta=\rho^{(0)} q^{\lambda}\in D_0(l)$, then
$h_0(x)=|h_0^+(x)|^2$ for $x\in T$ with $h_0^+$ given by
\begin{eqnarray}
h_0^+(x)&:=&\frac{\bigl(x;q\bigr)_{\infty}}{\bigl(t_0x;q\bigr)_{\infty}}
\prod_{k=1}^{l}\bigl(\vartheta_{k}x,\vartheta_{k}^{-1}x;q\bigr)_{\tau}
\nonumber\\
&=&\frac{\bigl(x;q\bigr)_{\infty}}{\bigl(t\vartheta_{l}x;q\bigr)_{\infty}}
\prod_{k=1}^{l}\frac{\bigl(\vartheta_{k}^{-1}x;q\bigr)_{\tau}}
{\bigl(t\vartheta_{k-1}x;q\bigr)_{\lambda_k-\lambda_{k-1}}}
\nonumber
\end{eqnarray}
where $\vartheta_0:=t^{-1}t_0$ and $\lambda_0:=0$. 
It follows that $h_0^+(x)$ is well defined for $x\in T$, since
the possible zero at $x=1$ 
of the factor $\bigl(t\vartheta_{l}x;q\bigr)_{\infty}$
in the denominator can be compensated by the zero at $x=1$ 
of the factor $\bigl(x;q\bigr)_{\infty}$.
Similarly, one deals with $h_1(x)$.
\end{proof}
Let $A_{\mathbb{R}}^W$ be the ${\mathbb{R}}$-algebra of 
$W$-invariant Laurent polynomials in the variables $z_1,\ldots,z_n$.
We obtain from Lemma \ref{techlem} the following corollary.
\begin{Cor}
Let $\underline{t}\in V_{AW}$ and $t\in (0,1)$.
Then the restriction of the bilinear form $\langle .,.
\rangle_{\underline{t},t}$ 
to $A_{\mathbb{R}}^W\times A_{\mathbb{R}}^W$ maps into ${\mathbb{R}}$ 
and is positive definite.
\end{Cor}
\begin{proof}
The monomials $m_{\lambda}$ ($\lambda\in\Lambda$)
are real valued on $F(r)\times T^{n-r}$ since $F(r)\subset {\mathbb{R}}^r$
by property {\bf (A)} for parameters in $V_{AW}$ (Definition  \ref{domain}),
so the assertion follows from Lemma \ref{techlem}{\bf (ii)}.
\end{proof}
The following theorem defines the Askey-Wilson polynomials for parameters
$\underline{t}\in V_{AW}$ and $t\in(0,1)$ as a special choice of orthogonal
basis for $A_{\mathbb{R}}^W$ with respect to the {\it positive definite}
bilinear form $\langle .,. \rangle_{\underline{t},t}: A_{\mathbb{R}}^W\times
A_{\mathbb{R}}^W\rightarrow {\mathbb{R}}$.
\begin{Thm}\label{orthopcd}
Let $t\in (0,1)$ and 
$\underline{t}\in V_{AW}$. Then there exists a unique basis
$\lbrace P_{\lambda}(.;\underline{t};t)\rbrace_{\lambda\in \Lambda}$ 
of $A_{\mathbb{R}}^W$ such that 

{\bf (i)}
$P_{\lambda}(.;\underline{t};t)=
m_{\lambda}+\sum_{\mu<\lambda}c_{\lambda,\mu}(\underline{t};t)m_{\mu}$,
some $c_{\lambda,\mu}(\underline{t};t)\in {\mathbb{R}}$;

{\bf (ii)} $\langle P_{\lambda}(.;\underline{t};t),
P_{\mu}(.;\underline{t};t)\rangle_{\underline{t},t}=0$ if $\mu\not=\lambda$.\\
Furthermore, $P_{\lambda}(.;\underline{t};t)$ is an eigenfunction of 
$D_{\underline{t},t}$ with eigenvalue $E_{\lambda}(\underline{t};t)$
and we have the explicit evaluation formula
\[
\langle P_{\lambda}(.;\underline{t};t),P_{\lambda}(.;\underline{t};t)\rangle_{
\underline{t},t}
=\mathcal{N}(\lambda;\underline{t};t),\quad \lambda\in \Lambda
\]
for the quadratic norms of the polynomials $P_{\lambda}$.
\end{Thm}
\begin{proof}
Fix $t\in (0,1)$ and $\underline{t}\in V_{AW}$.
Since $\langle .,. \rangle_{\underline{t},t}$ is positive definite on
$A_{\mathbb{R}}^W$, there exists for $\lambda\in\Lambda$ a unique $W$-invariant
Laurent polynomial $P_{\lambda}(.;\underline{t};t)\in A_{\mathbb{R}}^W$
satisfying ${\bf (i)}$ and the conditions 
$\langle P_{\lambda}(.;\underline{t};t), m_{\mu}\rangle_{\underline{t};t}=0$
for all $\mu<\lambda$. Furthermore, we have
\begin{equation}\label{GS}
P_{\lambda}(z;\underline{t};t)=m_{\lambda}(z)-\sum_{\mu<\lambda}
\frac{\langle m_{\lambda},P_{\mu}(.;\underline{t};t)\rangle_{\underline{t},t}}
{\langle P_{\mu}(.;\underline{t};t),P_{\mu}(.;\underline{t};t)
\rangle_{\underline{t},t}}
P_{\mu}(z;\underline{t};t),\quad \lambda\in\Lambda.
\end{equation}
By Lemma \ref{resn}, the polynomials 
$P_{\lambda}(z;\underline{t};t)=m_{\lambda}(z)+
\sum_{\mu<\lambda}c_{\lambda,\mu}(\underline{t};t)m_{\mu}(z)$
for  $\underline{t}\in V_0$ defined in Theorem
\ref{concl} can be given by the same formula \eqref{GS}.
Fix $\underline{t}\in V_{AW}\setminus V_{AW}^-$, where $V_{AW}^-$ is the
set of parameters $\underline{t}\in V_{AW}$ such that $t_i=\pm t^{-m}q^{-s}$
for some $i\in\lbrace 0,\ldots,3\rbrace$, $m\in\lbrace 0,\ldots,n-1\rbrace$
and $s\in {\mathbb{N}}_0$.
Let $\lbrace \underline{t}_k\rbrace_{k\in {\mathbb{N}}_0}$
be a sequence in $V_0$ converging to $\underline{t}$.
Then, by the Bounded Convergence Theorem,
\begin{equation}\label{contdependx}
\lim_{k\rightarrow\infty}\langle f,g\rangle_{\underline{t}_k,t}=
\langle f,g\rangle_{\underline{t},t},\quad \forall f,g\in A^W.
\end{equation}
Indeed, by assuming $\underline{t}\notin V_{AW}^{-}$, we have
that $F(r;\underline{\tau};t)=F(r;\underline{t};t)$ 
for $\underline{\tau}$ in an open
neighbourhood of $\underline{t}$  ($r=1,\ldots,n$) and that
no zeros in the denominator of the expression for
$\Delta_r^{AW}(\omega,.;\underline{t};t)$
($\omega\in F(r)$, $r=0,\ldots,n-1$) 
occur which need to be compensated by zeros in the numerator
(see the proof of Lemma \ref{techlem}). 
Hence the Bounded Convergence Theorem may be
applied at once.   

By induction on 
$\lambda$ we then obtain from \eqref{GS} and \eqref{contdependx} that
\begin{equation}\label{lims}
\lim_{k\rightarrow \infty}c_{\lambda,\mu}(\underline{t}_k;t)=c_{\lambda,\mu}
(\underline{t};t),\quad \mu<\lambda.
\end{equation}
By the residue calculus given in Lemma \ref{resn}, we can reformulate
Theorem \ref{concl} with respect
to the the bilinear form $\langle .,. \rangle_{\underline{\tau},t}$ for 
$\underline{\tau}\in V_0$. 
The theorem follows then for $\underline{t}\in V_{AW}\setminus V_{AW}^{-}$ 
by taking limits in the reformulated results using
Proposition \ref{triaAW} and \eqref{lims}.  

To prove the
theorem for $\underline{t}\in V_{AW}^{-}$, we use again a 
continuity argument. We treat here one typical example, the general case
is derived similarly. We assume
that $\underline{t}\in V_{AW}^{-}$ with $t_0=
t^{-m}q^{-s}$ for some $m\in\lbrace 0,\ldots,n-1\rbrace$, 
$s\in {\mathbb{N}}_0$ and that $t_i\not=t^{-l}q^{-s'}$ for all $i\in\lbrace
1,2,3\rbrace$, $l\in\lbrace 0,\ldots,n-1\rbrace$ and 
$s'\in {\mathbb{N}}_0$. Then there exists an $\epsilon>0$ 
such that $(\tau_0,t_1,t_2,t_3)\in V_{AW}\setminus V_{AW}^{-}$ 
and $F(r;\tau_0,t_1,t_2,t_3;t)=
F(r;\underline{t};t)$ for all $r\in\lbrace 1,\ldots,n
\rbrace$ and all $\tau_0\in {\mathbb{R}}_{>0}$ with $t_0-\tau_0<\epsilon$. 
We claim that
\begin{equation}\label{voorbeeld}
\lim_{\tau_0\uparrow t_0}\langle f,g\rangle_{\tau_0,t_1,t_2,t_3,t}=
\langle f,g\rangle_{\underline{t},t},\qquad \forall f,g\in A^W.
\end{equation}
We use the Bounded Convergence Theorem. 
In Lemma \ref{techlem} we have seen that zeros in the denominator of
the expression for the 
weight function $\Delta_r^{AW}(\omega,.;\underline{\tau};t)$
can occur when $\omega\in F(r;\underline{\tau};t)$ and
$\underline{\tau}\in V_{AW}^-$, and that
these zeros can be compensated by zeros in the numerator.
This amounts to replacing factors of the form
\[\frac{(1 \pm x)}{(1 \pm ux)}\]
in the weight function by $1$ when $u=1$.
For the application of the Bounded Convergence Theorem  in the limit
\eqref{voorbeeld}, it therefore suffices to note that the functions
\begin{equation}
h^{\pm}(u,x):=
\begin{cases}
\frac{(1\pm x)}{(1\pm ux)}\,\, &\hbox{ if } \,\, 
u\not= 1 \\
1 &\hbox{ if } u=1
\end{cases}
\end{equation}
are bounded on $U\times T$ where $U\subset {\mathbb{R}}_{>0}$ is some 
open set containing $1$. Now the theorem for the specific 
parameter values $\underline{t}$ follows by continuity
arguments from \eqref{voorbeeld}. 
\end{proof}
For parameters $\underline{t}\in V_{AW}$ with $|t_i|\leq 1$ the orthogonality
measure is completely continuous (i.e. we have $\langle .,.
\rangle=\langle .,. \rangle_0$) and coincides with Koornwinder's orthogonality
measure \cite{K1}. In particular the orthogonality relations  
reduce to Koornwinder's 
orthogonality relations (see \cite{K1}) and the quadratic norm
evaluations reduce to 
van Diejen's quadratic norm evaluations (see \cite{vD1}) for
parameter values $\underline{t}\in V_{AW}$ with $|t_i|\leq 1$ for all $i$.

Theorem \ref{orthopcd} implies 
that $D_{\underline{t},t}$ is symmetric with respect to 
$\langle .,. \rangle_{\underline{t},t}$.
The symmetry of $D$ and the orthogonality of the Askey-Wilson polynomials
with respect to $\langle .,. \rangle$ 
have been proved by different methods in \cite{S3}
for deformation parameter $t=q^{k}$ with $k\in {\mathbb{N}}$. 
In this case special case we can rewrite the bilinear form 
$\langle .,. \rangle_{\underline{t},q^k}$ using \eqref{AWposnatural} and we obtain
\begin{equation}\label{qSelbergAWd}
\begin{split}
\langle f,g\rangle=
\sum_{r=0}^n\frac{2^r\binom{n}{r}}{(2\pi i)^{n-r}}
\sum_{e_1,\ldots,e_r}
&\sum_{\stackrel{z_i\in\lbrace e_i,\ldots,e_iq^{N_{e_i}}\rbrace}{i=1,\ldots,r}}
\iint_{(z_{r+1},\ldots,z_n)\in T^{n-r}}\\
&f(z)g(z)\delta(z;q^k)\prod_{i=1}^rw_d(z_i;e_i)\prod_{j=r+1}^nw_c(z_j)
\frac{dz_j}{z_j}\\
\end{split}
\end{equation}
for $f,g\in A^W$. 

Theorem \ref{orthopcd} for $n=1$ reduces to the orthogonality relations
and norm evaluations stated in \cite[Theorem 2.5]{AW}.

\section{Limit transition to little $q$-Jacobi polynomials}
In this section we 
consider a limit case of the Askey-Wilson polynomials with 
positive partially discrete orthogonality measure
(Theorem \ref{orthopcd}) for which the continuous
part of the orthogonality measure disappears while the 
completely discrete part of the orthogonality
measure blows up to an infinite discrete measure.

We will obtain as limit the family of multivariable 
little $q$-Jacobi polynomials (previously introduced in \cite{S1})
which depends (besides on $q,t\in (0,1)$) on two parameters.
The parameter domain 
for the little $q$-Jacobi polynomials is defined as follows.
\begin{Def}
Let $V_L$ be the set of parameters $(a,b)$ for
which $a\in\bigl(0,1/q\bigr)$ and $b\in (-\infty,1/q\bigr)$.
\end{Def}
For functions 
$f: {\mathbb{C}}\rightarrow {\mathbb{C}}$ and $u,v\in {\mathbb{C}}$,
the Jackson $q$-integral of $f$ over $[u,v]$ is defined by
\begin{eqnarray}
\int_{u}^vf(x)d_qx&:=&\int_{0}^vf(x)d_qx-\int_0^uf(x)d_qx,\label{qJac1}\\
\int_{0}^vf(x)d_qx&:=&(1-q)\sum_{k=0}^{\infty}f(vq^k)vq^k,\label{qJac2}
\end{eqnarray}
provided that the infinite sums are absolutely convergent.
For a point $\xi\in ({\mathbb{C}}^*)^n$, we define
the Jackson integral of $f$ over the set
\begin{equation}\label{xin}
\langle \xi \rangle_n:=\lbrace \xi q^{\nu} \, | \, \nu\in P(n) \rbrace
\end{equation}
(here $\xi q^{\nu}:=(\xi_1q^{\nu_1},\ldots,\xi_nq^{\nu_n})$), by
\begin{equation}\label{multsumL}
\iint\limits_{\langle \xi\rangle_n}f(z)d_qz:=
(1-q)^n\sum_{\nu\in P(n)}f(\xi q^{\nu})\prod_{i=1}^n\xi_iq^{\nu_i}
\end{equation}
provided that the multisum is absolutely convergent.
Note that for 
special points $\xi=(\xi_1,\xi_1\gamma,\ldots,\xi_1\gamma^{n-1})\in 
\bigl({\mathbb{C}}^*\bigr)^n$, the multisum \eqref{multsumL}
can be expressed as an iterated Jackson integral by
\begin{equation}\label{iterJackint} 
\iint\limits_{\langle \xi\rangle_n}f(z)d_qz=
\int_{z_1=0}^{\xi_1}\int_{z_2=0}^{\gamma z_1}
\ldots\int_{z_n=0}^{\gamma z_{n-1}}
f(z)d_qz_n\ldots d_qz_1.
\end{equation}
Let $A_{\mathbb{R}}^S$ be the ${\mathbb{R}}$-algebra of $S$-invariant
polynomials in the variables $z_1,\ldots,z_n$.
An ${\mathbb{R}}$-basis for $A_{\mathbb{R}}^S$ is given by the set of monomials
$\lbrace \tilde{m}_{\lambda} \rbrace_{\lambda\in\Lambda}$, where
$\tilde{m}_{\lambda}(z):=\sum_{\mu\in S\lambda}z^{\mu}$.
Define a symmetric bilinear form $\langle .,. \rangle_{L,t}^{a,b}$ on
$A_{\mathbb{R}}^S$
for $t\in (0,1)$ and $(a,b)\in V_L$ by
\begin{equation}
\langle f,g\rangle_{L}:=\iint\limits_{\langle \rho_L\rangle_n}
f(z)g(z)\Delta^L(z)d_qz,\quad f,g\in A_{\mathbb{R}}^S
\end{equation}
where $\rho_{L,i}:=t^{i-1}$ and where
the weight function $\Delta^L(z)=\Delta^L(z;a,b;t)$ is given by
\begin{equation}\label{DL}
\Delta^L(z):=q^{-2\tau^2\binom{n}{3}}
t^{-(\alpha+1)\binom{n}{2}}
\left(\prod_{i=1}^nv_L(z_i)\right)\delta_{qJ}(z),\qquad
(a=q^{\alpha}, t=q^{\tau})
\end{equation}
with
\begin{equation}\label{vL}
v_L(x;a,b):=\frac{\bigl(qx;q\bigr)_{\infty}}
{\bigl( qbx;q\bigr)_{\infty}}x^{\alpha},\qquad (a=q^{\alpha})
\end{equation} 
and with interaction factor $\delta_{qJ}(z;t)$ given by
\begin{equation}\label{interactionL}
\delta_{qJ}(z;t):=
\prod_{1\leq i<j\leq n}|z_i-z_j||z_i|^{2\tau-1}\bigl(qt^{-1}z_j/z_i;q
\bigr)_{2\tau-1},\qquad (t=q^{\tau}).
\end{equation}
The function $v_L$ 
is exactly the weight function in the orthogonality measure for
the one-variable little $q$-Jacobi polynomials \cite{AA1}.
The same bilinear form $\langle .,. \rangle_{L}$ was
considered in \cite[section 5]{S1}, up to the positive constant 
$q^{-2\tau^2\binom{n}{3}}t^{-(\alpha+1)\binom{n}{2}}$. 
The weights $\Delta^L(z)$ 
in the bilinear form $\langle .,. \rangle_L$ are strict positive
for $z\in \langle \rho_L\rangle_n$ 
and $\langle f,g\rangle_L$, written out as a multidimensional infinite
sum, 
is absolutely convergent for all $f,g\in A_{\mathbb{R}}^S$ (see \cite{S1}). 

\begin{Def} Let $t\in (0,1)$ and $(a,b)\in V_L$.
The little $q$-Jacobi polynomials $\lbrace
P_{\lambda}^L(.;a,b;t)\rbrace_{\lambda\in\Lambda}$ are uniquely defined by
the two conditions\\
{\bf (a)} $P_{\lambda}^L=\tilde{m}_{\lambda}+\sum_{\mu<\lambda}
c_{\lambda,\mu}^L\tilde{m}_{\mu}$ for certain constants 
$c_{\lambda,\mu}^L=c_{\lambda,\mu}^L(a,b;t)\in
{\mathbb{R}}$;\\
{\bf (b)} 
$\langle P_{\lambda}^L,\tilde{m}_{\mu}\rangle_L=0$ for $\mu<\lambda$.\\ 
\end{Def}
The following proposition establishes the link between Askey-Wilson
polynomials with positive partially discrete orthogonality measure
and the little $q$-Jacobi polynomials. We use in the 
proposition the notation $|\lambda |:=\sum_{i=1}^n\lambda_i$
for the length of a partition $\lambda\in\Lambda$.
\begin{Prop}\label{link2}
Let $t\in (0,1)$,
$(a,b)\in V_L$ and define for $\epsilon\in {\mathbb{R}}^*$
\begin{equation}\label{tL}
\underline{t}_L(\epsilon):=\bigl(
\epsilon^{-1}q^{\frac{1}{2}}, -aq^{\frac{1}{2}}, 
\epsilon bq^{\frac{1}{2}}, -q^{\frac{1}{2}}\bigr).
\end{equation}
Then there exists 
a sequence $\lbrace \epsilon_k\rbrace_{k\in {\mathbb{N}}_0}$ in
${\mathbb{R}}_{>0}$ converging to $0$ such that
\begin{equation}\label{limmeasureL}
\begin{split}
\lim_{k\rightarrow \infty}
\left(\prod_{i=1}^n\bigl(-\epsilon_k^{-1}qt^{i-1},
-\epsilon_k^{-1}qat^{i-1};q\bigr)_{\infty}\right)&
\bigl(\epsilon_kq^{-\frac{1}{2}}\bigr)^{|\lambda |+|\mu |}
\langle m_{\lambda},m_{\mu}\rangle_{\underline{t}_L(\epsilon_k),t}\\
=2^nn!&\bigl(q;q\bigr)_{\infty}^{-2n}(1-q)^{-n}
\langle \tilde{m}_{\lambda},\tilde{m}_{\mu}\rangle_{L,t}^{a,b}\\
\end{split}
\end{equation}
for all $\lambda,\mu\in\Lambda$,
where $\langle .,. \rangle_{\underline{t},t}$ is given by \eqref{symmvormAW}.
\end{Prop}
The proof of the proposition will be given in section 8.
Note that 
$\underline{t}_L(\epsilon)\in V_{AW}$ for $\epsilon\in {\mathbb{R}}_{>0}$
sufficiently small, so $\langle .,. \rangle_{\underline{t}_L(\epsilon),t}$ 
is well defined for $\epsilon>0$ sufficiently small by Lemma \ref{techlem}.

We will use Proposition \ref{link2} to prove that the 
little $q$-Jacobi polynomials are limit cases of the Askey-Wilson polynomials
and to establish orthogonality relations and norm evaluations 
for the little $q$-Jacobi polynomials with
respect to the bilinear form $\langle .,. \rangle_L$. 

We use the following 
definition of limit transitions between $S$-invariant Laurent polynomials
(cf. \cite{SK1}). Let $f(.;u)$ ($u\in {\mathbb{R}}^*$) and $f$ be 
$S$-invariant Laurent polynomials in $n$ variables $z_1,\ldots,z_n$, 
then we write $\lim_{u\rightarrow 0}f(.;u)=f$ if 
$\lim_{u\rightarrow 0} f(z;u)=f(z)$ for all $z\in ({\mathbb{R}}^*)^n$. 
Note that the ${\mathbb{R}}$-algebra of $S$-invariant Laurent polynomials 
has as ${\mathbb{R}}$-basis
the set of monomials $\lbrace \tilde{m}_{\lambda}(z)
\rbrace_{\lambda\in\tilde{\Lambda}}$, where $\tilde{\Lambda}:=
\lbrace \lambda\in {\mathbb{Z}}^n \, | \, \lambda_1\geq\lambda_2\geq\ldots
\geq \lambda_n \rbrace$ 
and $\tilde{m}_{\lambda}(z):=\sum_{\mu\in S\lambda}z^{\mu}$.
If $f(.;u)=
\sum_{\lambda\in\tilde{\Lambda}}c_{\lambda}(u)\tilde{m}_{\lambda}$\, ($u\in
{\mathbb{R}}^*$) and 
$f=\sum_{\lambda\in\tilde{\Lambda}}c_{\lambda}\tilde{m}_{\lambda}$ 
such that 
$\lbrace \lambda\in\tilde{\Lambda} \, | \, c_{\lambda}(u)\not=0\rbrace$
is contained in some finite $u$-independent subset 
for $|u|$ sufficiently small, then we have
$\lim_{u\rightarrow 0}f(.;u)=f$
iff $\lim_{u\rightarrow 0}c_{\lambda}(u)=c_{\lambda}$ for all
$\lambda\in\tilde{\Lambda}$.   
Crucial in the limit from Askey-Wilson polynomials to
little $q$-Jacobi polynomials is a limit from rescaled monomials
$m_{\lambda}(z|u)$ to $\tilde{m}_{\lambda}(z)$, where the rescaled 
monomial $m_{\lambda}(z|u)$ for
$u\in {\mathbb{R}}^*$ is the
$S$-invariant Laurent polynomial given by
\begin{equation}\label{mlambdau}
m_{\lambda}(z|u):=u^{|\lambda |}m_{\lambda}(u^{-1}z),\quad\lambda\in\Lambda.
\end{equation}
Here we have used the notation $u^{-1}z:=(u^{-1}z_1,\ldots,u^{-1}z_n)$.
In terms of the basis $\lbrace
\tilde{m}_{\mu}\rbrace_{\mu\in\tilde{\Lambda}}$, we have
$m_{\lambda}(z|u)=\sum_{\mu\in\tilde{\Lambda}\cap
W\lambda}d_{\lambda,\mu}(u)\tilde{m}_{\mu}(z)$ for $\lambda\in\Lambda$
with $d_{\lambda,\mu}(u)$ homogeneous of degree
$|\lambda|-|\mu|$ and $d_{\lambda,\lambda}(u)\equiv 1$. 
Furthermore we have for $\lambda\in\Lambda$ that
$|\mu|\leq |\lambda|$ if $\mu\in W\lambda$ and 
$|\lambda|=|\mu|$ iff $\mu\in S\lambda$.
Hence we obtain the limit transitions
\begin{equation}\label{limmon}
\lim_{u\rightarrow 0}m_{\lambda}(z|u)=\tilde{m}_{\lambda}(z) \qquad
(\lambda\in\Lambda).
\end{equation}
We express the quadratic norm of the little $q$-Jacobi polynomials in terms of 
functions ${\mathcal{N}}^+_{qJ}(\lambda)={\mathcal{N}}^+_{qJ}(\lambda;a,b;t)$
and 
${\mathcal{N}}^-_{qJ}(\lambda)={\mathcal{N}}^-_{qJ}(\lambda;a,b;t)$ which are
defined by
\begin{equation}\label{DeltaJhat}
\begin{split}
{\mathcal{N}}^+_{qJ}(\lambda):=\prod_{i=1}^n
\frac{\Gamma_q(\lambda_i+1+(n-i)\tau+\alpha+\beta)
\Gamma_q(\lambda_i+1+(n-i)\tau+\alpha)}
{\Gamma_q(2\lambda_i+1+2(n-i)\tau+\alpha+\beta)}&\\
.\prod_{1\leq j<k\leq n}\left(
\frac{\Gamma_q(\lambda_j+\lambda_k+1+(2n-j-k+1)\tau+\alpha+\beta)}
{\Gamma_q(\lambda_j+\lambda_k+1+(2n-j-k)\tau+\alpha+\beta)}\right.&\\
\left.
.\frac{\Gamma_q(\lambda_j-\lambda_k+(k-j+1)\tau)}
{\Gamma_q(\lambda_j-\lambda_k+(k-j)\tau)}\right)&\\
\end{split}
\end{equation}
\begin{equation}\label{DeltaJtilde}
\begin{split}
{\mathcal{N}}^-_{qJ}(\lambda):=\prod_{i=1}^n
\frac{\Gamma_q(\lambda_i+1+(n-i)\tau)
\Gamma_q(\lambda_i+1+(n-i)\tau+\beta)}
{\Gamma_q(2\lambda_i+2+2(n-i)\tau+\alpha+\beta)}\qquad&\\
.\prod_{1\leq j<k\leq n}\left(
\frac{\Gamma_q(\lambda_j+\lambda_k+2+(2n-j-k-1)\tau+\alpha+\beta)}
{\Gamma_q(\lambda_j+\lambda_k+2+(2n-j-k)\tau+\alpha+\beta)}\right.&\\
\left.
.\frac{\Gamma_q(\lambda_j-\lambda_k+1+(k-j-1)\tau)}
{\Gamma_q(\lambda_j-\lambda_k+1+(k-j)\tau)}\right)&,\\
\end{split}
\end{equation}
where $a=q^{\alpha},b=q^{\beta}$, $t=q^{\tau}$ and where
\begin{equation}\label{qgamma}
\Gamma_q(u):=\frac{\bigl(q;q\bigr)_{u-1}}{(1-q)^{u-1}}\qquad (u\not\in
-{\mathbb{N}}_0)
\end{equation}
is the $q$-gamma function.
For $\lambda\in\Lambda$, $(a,b)\in V_L$ and $t\in (0,1)$ let
${\mathcal{N}}^L(\lambda)={\mathcal{N}}^L(\lambda;a,b;t)$ be given by
\begin{equation}\label{NL}
{\mathcal{N}}^L(\lambda):=
q^{\sum_{i=1}^n(\lambda_i+\alpha+2(n-i)\tau)\lambda_i}
{\mathcal{N}}^+_{qJ}(\lambda){\mathcal{N}}^-_{qJ}(\lambda).
\end{equation}
Note that ${\mathcal{N}}^L(\lambda)$ is well defined and positive.
We have now the following theorem.
\begin{Thm}\label{thmL} 
Let $(a,b)\in V_L$ and $t\in (0,1)$. There exists a sequence
$\lbrace \epsilon_k\rbrace_{k\in {\mathbb{N}}_0}$ 
in ${\mathbb{R}}_{>0}$ converging
to $0$ such that
\begin{equation}\label{limpolL}
\lim_{k\rightarrow \infty} 
\left(q^{-\frac{1}{2}}\epsilon_k\right)^{|\lambda|}
P_{\lambda}\bigl(q^{\frac{1}{2}}\epsilon_k^{-1}z;
\underline{t}_L(\epsilon_k);t\bigr)=
P_{\lambda}^L(z;a,b;t)
\end{equation}
for all $\lambda\in\Lambda$.
Furthermore, the polynomials $\lbrace P_{\lambda}^L\rbrace_{\lambda\in\Lambda}$
are orthogonal with respect to $\langle .,. \rangle_L$
and the quadratic norms of the little $q$-Jacobi polynomials are given by
\begin{equation}\label{normL}
\langle P_{\lambda}^L,P_{\lambda}^L\rangle_L={\mathcal{N}}^L(\lambda),\quad
\lambda\in\Lambda.
\end{equation}
\end{Thm}
\begin{proof} 
We write
\begin{equation}\label{ex}
\begin{split}
\bigl(q^{-\frac{1}{2}}\epsilon\bigr)^{|\lambda|}P_{\lambda}\bigl(z;
\underline{t}_L(\epsilon);t)&=\sum_{\mu\leq\lambda}c_{\lambda,\mu}(\epsilon)
\bigl(q^{-\frac{1}{2}}\epsilon\bigr)^{|\mu|}m_{\mu}(z),\\
P_{\lambda}^L(z;a,b;t)&=
\sum_{\mu\leq\lambda}c_{\lambda,\mu}^L\tilde{m}_{\mu}(z)\\
\end{split}
\end{equation}
for the expansions 
of the Askey-Wilson polynomial and the little $q$-Jacobi polynomial
in terms of monomials. In particular, we have
$c_{\lambda,\lambda}(\epsilon)=c_{\lambda,\lambda}^L=1$.
Let $\lbrace\epsilon_k\rbrace_{k\in {\mathbb{N}}_0}$ be a sequence in
${\mathbb{R}}_{>0}$ converging 
to $0$ such that \eqref{limmeasureL} is satisfied
for all $\lambda,\mu\in\Lambda$. 
We will prove that
\begin{equation}\label{limpolLx}
\lim_{k\rightarrow\infty}c_{\lambda,\mu}(\epsilon_k)=c_{\lambda,\mu}^L
\quad \forall\, \mu\leq\lambda,
\end{equation}
then the limit \eqref{limpolL} follows from \eqref{limmon}, \eqref{ex} 
and \eqref{limpolLx}. We proceed by induction on $\lambda$.
The case $\lambda=0$ is trivial. For $\lambda\not=0$, note that
$\underline{t}_L(\epsilon)\in V_{AW}$ for $\epsilon>0$ sufficiently small,  
hence by Theorem \ref{orthopcd} we can write
\begin{equation}\label{expansLL}
\bigl(q^{-\frac{1}{2}}\epsilon\bigr)^{|\lambda|}
P_{\lambda}\bigl(q^{\frac{1}{2}}z;
\underline{t}_L(\epsilon);t\bigr)=
\bigl(q^{-\frac{1}{2}}\epsilon\bigr)^{|\lambda|}m_{\lambda}(z)-
\sum_{\mu<\lambda}d_{\lambda,\mu}(\epsilon)
\bigl(q^{-\frac{1}{2}}\epsilon\bigr)^{|\mu|}P_{\mu}\bigl(z;
\underline{t}_L(\epsilon);t\bigr)
\end{equation}
with
\[
d_{\lambda,\mu}(\epsilon):=
\frac{(q^{-\frac{1}{2}}\epsilon)^{|\lambda|+|\mu|}
\langle m_{\lambda},P_{\mu}(.;\underline{t}_L(\epsilon);t)
\rangle_{\underline{t}_L(\epsilon),t}}
{(q^{-\frac{1}{2}}\epsilon)^{2|\mu|}
\langle P_{\mu}(.;\underline{t}_L(\epsilon);t),
P_{\mu}(.;\underline{t}_L(\epsilon);t)
\rangle_{\underline{t}_L(\epsilon),t}}
\]
for $\epsilon>0$ sufficiently small.
For the little $q$-Jacobi polynomials we have
\begin{equation}\label{expansLLL}
P_{\lambda}^L(z;a,b;t)=\tilde{m}_{\lambda}(z)-
\sum_{\mu<\lambda}d_{\lambda,\mu}^LP_{\mu}^L(z;a,b;t),
\end{equation}
with
\[
d_{\lambda,\mu}^L=
\frac{\langle \tilde{m}_{\lambda},P_{\mu}^L(.;a,b;t)
\rangle_{L,t}^{a,b}}
{\langle
P_{\mu}^L(.;a,b;t),P_{\mu}^L(.;a,b;t)
\rangle_{L,t}^{a,b}}.
\]
By the induction hypotheses and Proposition \ref{link2}, we have
\begin{equation}\label{partLL}
\lim_{k\rightarrow\infty}d_{\lambda,\mu}(\epsilon_k)=d_{\lambda,\mu}^L
\quad \forall \, \mu<\lambda.
\end{equation}
So \eqref{limpolLx} follows by the induction hypotheses from \eqref{ex},
\eqref{expansLL},
\eqref{expansLLL} and \eqref{partLL}.

The orthogonality relations and the norm evaluations 
for the little $q$-Jacobi polynomials
follow now from the orthogonality relations and norm evaluations for the 
Askey-Wilson polynomials with partially discrete orthogonality
measure (Theorem \ref{orthopcd}), Proposition \ref{link2}, \eqref{limpolLx}
and the observation that
\begin{equation}
\begin{split}
\lim_{\epsilon\rightarrow 0}\left(\prod_{i=1}^n
\bigl(-\epsilon^{-1}qt^{i-1},-\epsilon^{-1}qat^{i-1};q\bigr)_{\infty}\right)
&(q^{-1/2}\epsilon)^{2|\lambda|}
{\mathcal{N}}(\lambda;\underline{t}_L(\epsilon);t)\nonumber\\
=2^nn!&\bigl(q;q\bigr)_{\infty}^{-2n}(1-q)^{-n}
{\mathcal{N}}^L(\lambda;a,b;t).\nonumber\\
\end{split}
\end{equation}
\end{proof}
Full orthogonality of the little $q$-Jacobi polynomials 
was proved in \cite{S1} by means of an explicit second order $q$-difference
operator $D_L$ which diagonalizes the little $q$-Jacobi polynomials. 
The operator $D_L$ 
and the corresponding eigenvalue equations can be obtained from 
\eqref{limpolL} by taking the limit $k\rightarrow\infty$ in the
equations
\begin{equation}\label{limLformal}
\bigl(q^{-1}\epsilon_k\bigr)^{|\lambda|}
\bigl((D-E_{\lambda})P_{\lambda}\bigr)(q^{\frac{1}{2}}\epsilon_k^{-1}z;
\underline{t}_L(\epsilon_k);t)=0\quad (\lambda\in\Lambda)
\end{equation}
where $D$ is given by \eqref{secondorderqdiff} and $E_{\lambda}$
is given by \eqref{AWb} 
(see \cite[section 3]{S1} for the formal computation of the
limits of $D$ and $E_{\lambda}$).
In \cite{SK1} it was shown 
that the formal computation of the limits of $D$ and 
$E_{\lambda}$ in \cite[section 3]{S1} 
can be used to prove the limit transition 
\eqref{limpolL} for generic $t\in (0,1)$. 
See also \cite{S2} for the special case that 
$t=q^k$, $k\in {\mathbb{N}}$. 

The constant term identity for the little $q$-Jacobi polynomials
can be rewritten as follows.  
\begin{Cor}\label{gevolg}
Let $t\in (0,1)$ and $(a,b)\in V_L$. We have
\begin{equation}\label{SelbergL}
\langle 1,1\rangle_{L,t}^{a,b}=\prod_{j=1}^n\frac{
\Gamma_q(\alpha+1+(j-1)\tau)\Gamma_q(\beta+1+(j-1)\tau)\Gamma_q(j\tau)}
{\Gamma_q(\alpha+\beta+2+(n+j-2)\tau)\Gamma_q(\tau)}.
\end{equation}
\end{Cor}
The constant term identity \eqref{SelbergL} has been studied extensively in
recent years. 
It was conjectured by Askey \cite{A} for $t=q^k$, $k\in {\mathbb{N}}$
and proved in 
this case independently by Habsieger \cite{H} and Kadell \cite{K}. 
For arbitrary 
$t\in (0,1)$ the first proof appeared in Aomoto's paper \cite{Ao} 
(see also \cite{Ka} and \cite{TV} for alternative proofs). 
\section{Limit transition to big $q$-Jacobi polynomials}
In this section, we 
consider an other limit transition involving the Askey-Wilson polynomials with
positive partially discrete orthogonality measure
(Theorem \ref{orthopcd}) for which the continuous
part of the orthogonality measure disappears while the 
completely discrete part of the orthogonality
measure blows up to an infinite discrete measure.
Instead of sending 
one parameter to infinity, as we did in the previous section,
we will send now two parameters to infinity.
We will obtain the four parameter 
family of big $q$-Jacobi polynomials  
(previously introduced in \cite{S1}) as limit case of the five
parameter family of Askey-Wilson polynomials.
The parameter domain 
for the big $q$-Jacobi polynomials $V_B$ is defined as follows.
\begin{Def}
Let $V_B$ be the set of parameters $(a,b,c,d)$ for
which $c,d>0$ and $a\in\bigl(-c/dq,1/q\bigr)$, $b\in\bigl(-d/cq,1/q\bigr)$
or $a=cu$, $b=-d\overline{u}$ with $u\in {\mathbb{C}}\setminus {\mathbb{R}}$.
\end{Def}
Before we define the orthogonality measure for the big $q$-Jacobi polynomials 
we first need to introduce some more notations. We set
\[\langle \eta,\xi\rangle_n:=\bigcup_{j=0}^n\langle \xi\rangle_j\times\langle
\eta\rangle_{n-j} \subset {\mathbb{C}}^n
\]
where $\eta,\xi\in ({\mathbb{C}}^*)^n$, $\underline{c}=(c_0,\ldots,c_n)\in 
({\mathbb{C}}^*)^{n+1}$ and $\langle \xi\rangle_n$ is defined by \eqref{xin}. 
Here we use the convention that $\langle\xi\rangle_n\times
\langle \eta\rangle_0=\langle\xi\rangle_n$ and
$\langle\xi\rangle_0\times\langle\eta\rangle_n=\langle\eta\rangle_n$.
Define the $\underline{c}$-weighted Jackson integral of $f$ 
over the set $\langle \eta,\xi\rangle_n$ by
\begin{equation}\label{defJacint}
\begin{split}
\iint\limits_{\langle \eta,\xi\rangle_n}f(z)d_q^{\underline{c}}z:=
&\sum_{j=0}^n(-1)^{n-j}c_j\iint\limits_{z\in\langle \xi\rangle_j}
\,\,\,\iint\limits_{w\in\langle \eta\rangle_{n-j}}f(z,w)d_qzd_qw\\
=&(1-q)^n\sum_{j=0}^n\sum_{\stackrel{\nu\in P(j)}{\nu'\in P(n-j)}}
c_jf(\xi q^{\nu},\eta q^{\nu'})\prod_{l=1}^j\xi_lq^{\nu_l}
\prod_{m=1}^{n-j}(-\eta_mq^{\nu_m'}),\\
\end{split}
\end{equation}
where the $j=0$ respectively $j=n$ term in \eqref{defJacint} should be read as
\[(-1)^nc_0\iint\limits_{w\in \langle \eta\rangle_n}f(w)d_qw=
(1-q)^n\sum_{\nu'\in P(n)}c_0f(\eta
q^{\nu'})\prod_{m=1}^n\bigl(-\eta_mq^{\nu_m'}\bigr)\]
respectively
\[c_n\iint\limits_{z\in \langle \xi\rangle_n}f(z)d_qz=
(1-q)^n\sum_{\nu\in P(n)}c_nf(\xi q^{\nu})\prod_{l=1}^n\xi_lq^{\nu_l}.\]
If $\eta=(\eta_1,\eta_1\gamma',\ldots,\eta_1(\gamma')^{n-1})$ and
$\xi=(\xi_1,\xi_1\gamma,\ldots,\xi_1\gamma^{n-1})$, 
then the ${\underline{c}}$-weighted
Jackson integral over $\langle \eta,\xi\rangle_n$ can be rewritten as an 
iterated Jackson
integral by 
\begin{equation}
\begin{split}
\iint\limits_{\langle \eta,\xi\rangle_n}f(z)d_q^{\underline{c}}z=
\sum_{j=0}^nc_j&\int_{z_1=0}^{\xi_1}\int_{z_2=0}^{\gamma z_1}\ldots
\int_{z_j=0}^{\gamma z_{j-1}}\nonumber\\
&\int_{z_{j+1}=\eta_1}^0
\int_{z_{j+2}=
\gamma'z_{j+1}}^0\ldots\int_{z_n=\gamma'z_{n-1}}^0f(z)d_qz_n\ldots
d_qz_1.\nonumber\\
\end{split}
\end{equation}
Define a symmetric bilinear form $\langle .,. \rangle_{B,t}^{a,b,c,d}$
on $A_{\mathbb{R}}^S$ for parameters $t\in (0,1)$ and $(a,b,c,d)\in V_B$ by 
\begin{equation}
\langle f,g \rangle_B:=\iint\limits_{\langle \sigma_B,\rho_B\rangle_n}
f(z)g(z)\Delta^B(z)d_q^{{\underline{c}}_B}z,\quad
f,g\in A_{\mathbb{R}}^S,
\end{equation}
with $\rho_{B,i}:=ct^{i-1}$, $\sigma_{B,i}:=-dt^{i-1}$ and with weight function
\begin{equation}\label{DB}
\Delta^B(z):=\left(\prod_{i=1}^nv_B(z_i)\right)\delta_{qJ}(z),
\end{equation}
where $v_B$ is the weight function in the orthogonality measure for the 
one-variable big $q$-Jacobi polynomials \cite{AA2},  
\begin{equation}\label{vB}
v_B(x;a,b,c,d):=
\frac{\bigl(qx/c,-qx/d;q\bigr)_{\infty}}{\bigl(qax/c,-qbx/d;q\bigr)_{\infty}}
\end{equation}
and $\delta_{qJ}(z)=\delta_{qJ}(z;t)$ is given by \eqref{interactionL}.
The weight $\underline{c}_B$ is of the form
$c_{B,j}:=c_Bd_{B,j}$, with
\begin{equation}\label{dBj}
d_{B,j}:=\prod_{\stackrel{1\leq k<m\leq n}{k\leq j}}\Psi_t(-t^{n-m-k+1}d/c)
\end{equation}
where $\Psi_t(x)$ is defined by
\begin{equation}\label{Psit}
\Psi_t(x):=|x|^{2\tau-1}\frac{\theta(tx)}{\theta(qt^{-1}x)}
\end{equation}
with $\theta(x)$ the Jacobi theta function
\begin{equation}
\label{thetafunction}
\theta(x):=\bigl(q,x,qx^{-1};q\bigr)_{\infty},
\end{equation}
and the constant $c_B$ is defined by
\begin{equation}\label{cB}
c_B:=\frac{\bigl(q;q\bigr)_{\infty}^nq^{-2\tau^2\binom{n}{3}}
d^{-2\tau\binom{n}{2}-n}t^{-\binom{n}{2}}}
{\prod_{i=1}^n\theta(-t^{1-i}c/d)}.
\end{equation}
The positive constant $c_B$ is not essential for the definition of
$\langle .,. \rangle_B$. We have chosen to take this
constant within the definition of $\langle .,. \rangle_B$
because it will simplify formulas and notations later on. 
The Jacobi theta function satisfies the functional relation
\begin{equation}\label{theta}
\theta(q^kx)=(-x^{-1})^kq^{-\binom{k}{2}}\theta(x),\quad k\in {\mathbb{N}}_0.
\end{equation}
This implies that $\Psi_t$ is a quasi constant, 
i.e. $\Psi_t(qx)=\Psi_t(x)$. In particular, the weight 
$\underline{d}_B$ \eqref{dBj} 
is independent of $a,b$ and quasi constant in the parameters $c,d$.
We note that the bilinear form $\langle .,. \rangle_B$ 
is the same as the one defined in  
\cite{S1} up to the positive constant 
$c_B$ \eqref{cB}. This is easily verified using
the fact that $\langle .,. \rangle_B$ 
is defined as bilinear form on the space of {\it
symmetric} polynomials.

The weight $\underline{c}_B$ in the definition
of $\langle .,. \rangle_B$ is needed in order to obtain good
asymptotic behaviour. To be more precise, let
$j\in \lbrace 1,\ldots,n\rbrace$,
$\lambda\in P(j-1)$, $\mu\in P(n-j)$ and set
$\lambda^{(l)}:=(\lambda,l)\in P(j)$ for $l\geq\lambda_{j-1}$,
respectively $\mu^{(m)}:=(\mu,m)\in P(n-j+1)$ for $m\geq \mu_{n-j}$.
For $j=1$ respectively $j=n$, this should be read as $\lambda^{(l)}=l\in P(1)$,
respectively $\mu^{(m)}=m\in P(1)$.
Define
\begin{equation}
\begin{split}
z^+(l;\lambda,\mu)&:=\bigl(\rho_Bq^{\lambda^{(l)}},\sigma_Bq^{\mu}\bigr)\in
\langle\rho_B\rangle_j\times\langle\sigma_B\rangle_{n-j}\quad
(l\geq\lambda_{j-1}),\\
z^-(m;\lambda,\mu)&:=\bigl(\rho_Bq^{\lambda},\sigma_Bq^{\mu^{(m)}}\bigr)\in
\langle\rho_B\rangle_{j-1}\times\langle\sigma_B\rangle_{n-j+1}
\quad (m\geq\mu_{n-j}),\\
\end{split}
\end{equation}
then we have
\[\lim_{l\rightarrow\infty}z^+(l;\lambda,\mu)=
(\rho_Bq^{\lambda},0,\sigma_Bq^{\mu}),
\quad 
\lim_{m\rightarrow\infty}z^-(m;\lambda,\mu)=
(\rho_Bq^{\lambda},\sigma_Bq^{\mu},0)
\]
(with the obvious conventions when $j=1$ respectively $j=n$).
We have now good asymptotic behaviour in the following sense.
\begin{Lem}\label{motivation}
Let $(a,b,c,d)\in V_B$ and $t\in (0,1)$.
Then 
\begin{equation}\label{asym}
\lim_{l\rightarrow\infty}
c_{B,j}\Delta^B(z^+(l;\lambda,\mu))=
\lim_{m\rightarrow\infty}c_{B,j-1}\Delta^B(z^-(m;\lambda,\mu))
\end{equation}
for all $\lambda\in P(j-1)$, 
$\mu\in P(n-j)$ and $j\in\lbrace 1,\ldots,n\rbrace$.
The conditions \eqref{asym} for $\lambda\in P(j-1)$, $\mu\in P(n-j)$ and 
$j\in\lbrace 1,\ldots,n\rbrace$ characterize 
the weight $\underline{c}_B$ uniquely up to a constant. 
\end{Lem}
\begin{proof}
See the proof of \cite[Theorem 6.5]{S1}.
\end{proof}
In \cite{S1} it was remarked that $d_{B,j}=1$ for all $j\in\lbrace
1,\ldots,n\rbrace$ if $t=q^k$ with $k\in {\mathbb{N}}$. 
For the weight 
$\underline{c}_B$ we thus have $c_{B,j}=c_B$ for all $j$ if $t=q^k$ with
$k\in {\mathbb{N}}$, and $c_B$ can then be rewritten as
\begin{equation}\label{cjk}
c_B=\frac{q^{\binom{k}{2}\binom{n}{2}-k^2\binom{n}{3}}(c+d)^n}
{\bigl(-d/c,-c/d;q\bigr)_{\infty}^{n}
(cd)^{n+\binom{n}{2}k}}
\qquad (t=q^k,k\in {\mathbb{N}}).
\end{equation}
This follows by a straightforward calculation using the relation
$\theta(qx^{-1})=\theta(x)$, \eqref{theta} and
\begin{equation}
\sum_{i=1}^n(i-1)^2=2\binom{n}{3}+\binom{n}{2}.
\end{equation}
It was proved in \cite{S1} that the weights $\Delta^B(z)$ 
in the bilinear form $\langle .,.
\rangle_{B}$ are strict positive for $z\in \langle \sigma_B,\rho_B\rangle_n$
and that $\langle f,g\rangle_B$, written out as a multidimensional infinite
sum over $\langle \sigma_B,\rho_B\rangle_n$, 
is absolutely convergent for all $f,g\in A_{\mathbb{R}}^S$. 
\begin{Def} Let $t\in (0,1)$ and $(a,b,c,d)\in V_B$.
The big $q$-Jacobi polynomials $\lbrace
P_{\lambda}^B(.;a,b,c,d;t)\rbrace_{\lambda\in\Lambda}$ are uniquely defined by
the two conditions\\
{\bf (a)} $P_{\lambda}^B=\tilde{m}_{\lambda}+\sum_{\mu<\lambda}
c_{\lambda,\mu}^B\tilde{m}_{\mu}$ for certain constants 
$c_{\lambda,\mu}^B=c_{\lambda,\mu}^B(a,b,c,d;t)\in
{\mathbb{R}}$;\\
{\bf (b)} $\langle P_{\lambda}^B,\tilde{m}_{\mu}\rangle_B=0$ for $\mu<\lambda$.
\end{Def}
The following proposition is the analogue of Proposition \ref{link2} 
for the big $q$-Jacobi polynomials.
\begin{Prop}\label{link3}
Let $t\in (0,1)$,
$(a,b,c,d)\in V_B$ and define for $\epsilon\in {\mathbb{R}}^*$
\begin{equation}\label{tB}
\underline{t}_B(\epsilon):=\bigl(
\epsilon^{-1}(qc/d)^{\frac{1}{2}}, -\epsilon^{-1}(qd/c)^{\frac{1}{2}}, 
\epsilon a(qd/c)^{\frac{1}{2}},
-\epsilon b(qc/d)^{\frac{1}{2}} \bigr).
\end{equation}
Then there exists a sequence $\lbrace \epsilon_k\rbrace_{k\in {\mathbb{N}}_0}$
in ${\mathbb{R}}_{>0}$ converging to $0$ such that
\begin{equation}\label{limmeasureB}
\begin{split}
\lim_{k\rightarrow \infty}
\left(\prod_{i=1}^n\bigl( -\epsilon_k^{-2}qt^{i-1};q\bigr)_{\infty}\right)
&\bigl(\epsilon_k(cd/q)^{\frac{1}{2}}\bigr)^{|\lambda |+|\mu |}
\langle m_{\lambda},m_{\mu}\rangle_{\underline{t}_B(\epsilon_k),t}\\
&=2^nn!\bigl(q;q\bigr)_{\infty}^{-2n}(1-q)^{-n}
\langle \tilde{m}_{\lambda},\tilde{m}_{\mu}\rangle_{B,t}^{a,b,c,d}\\
\end{split}
\end{equation}
for all $\lambda,\mu\in\Lambda$,
where $\langle .,. \rangle_{\underline{t},t}$ is given by \eqref{symmvormAW}.
\end{Prop}
The proof will be given in section 9.
Note that $\underline{t}_B(\epsilon)\in V_{AW}$ 
for $\epsilon\in {\mathbb{R}}_{>0}$
sufficiently small, so $\langle .,. \rangle_{\underline{t}_B(\epsilon),t}$
is well defined for $\epsilon>0$ sufficiently small by Lemma \ref{techlem}.

We can repeat now the arguments of the previous section to establish
full orthogonality  of the big $q$-Jacobi polynomials 
with respect to $\langle .,. \rangle_B$
and to calculate their norms.
For $\lambda\in\Lambda$, $(a,b,c,d)\in V_B$ and $t\in (0,1)$
let ${\mathcal{N}}^B(\lambda)={\mathcal{N}}^B(\lambda;a,b,c,d;t)$ be given by
\begin{equation}
\begin{split}
{\mathcal{N}}^B(\lambda):=
(cd)^{|\lambda|}&q^{\frac{1}{2}\sum_{i=1}^n(\lambda_i-1+2(n-i)\tau)\lambda_i}
{\mathcal{N}}^+_{qJ}(\lambda;a,b;t){\mathcal{N}}^-_{qJ}(\lambda;a,b;t)\\
&. 
\prod_{i=1}^n\bigl(-q^{\lambda_i+1}bt^{n-i}c/d,-q^{\lambda_i+1}at^{n-i}d/c;q
\bigr)_{\infty}^{-1}\\
\end{split}
\end{equation}
where $t=q^{\tau}$ and ${\mathcal{N}}^+_{qJ}$ respectively 
${\mathcal{N}}^-_{qJ}$ is given by
\eqref{DeltaJhat} respectively \eqref{DeltaJtilde}.
Note that ${\mathcal{N}}^B(\lambda)$ is well defined and positive. 
We have the following theorem.
\begin{Thm}\label{thmB}
Let $t\in (0,1)$ and $(a,b,c,d)\in V_B$. There exists a sequence
$\lbrace \epsilon_k\rbrace_{k\in {\mathbb{N}}_0}$ 
in ${\mathbb{R}}_{>0}$ converging
to $0$ such that
\begin{equation}\label{limpolB}
\lim_{k\rightarrow \infty} 
\left(\epsilon_k(cd/q)^{\frac{1}{2}}\right)^{|\lambda|}
P_{\lambda}\bigl((q/cd)^{\frac{1}{2}}\epsilon_k^{-1}z;
\underline{t}_B(\epsilon_k);t
\bigr)=
P_{\lambda}^B(z;a,b,c,d;t)
\end{equation}
for all $\lambda\in\Lambda$.
Furthermore, the polynomials 
$\lbrace P_{\lambda}^B\rbrace_{\lambda\in\Lambda}$ 
are orthogonal with respect to $\langle .,. \rangle_B$ and 
the quadratic norms of the big $q$-Jacobi polynomials are given by
\begin{equation}\label{normB}
\langle P_{\lambda}^B,P_{\lambda}^B\rangle_B={\mathcal{N}}^B(\lambda),\quad
\lambda\in\Lambda.
\end{equation}
\end{Thm}
\begin{proof}
We have the limit  
\begin{equation}
\begin{split}
\lim_{\epsilon\rightarrow 0}\left(\prod_{i=1}^n
\bigl(-\epsilon^{-2}qt^{i-1};q\bigr)_{\infty}\right)&\bigl((cd/q)^{\frac{1}{2}}
\epsilon\bigr)^{2|\lambda|}
{\mathcal{N}}(\lambda;\underline{t}_B(\epsilon);t)\nonumber\\
&=2^nn!\bigl(q;q\bigr)_{\infty}^{-2n}(1-q)^{-n}
{\mathcal{N}}^B(\lambda;a,b,c,d;t).\nonumber\\
\end{split}
\end{equation}
The proof is now analogous to the proof of Theorem \ref{thmL}. 
\end{proof}
Full orthogonality of the big $q$-Jacobi polynomials
was proved in \cite{S1} by studying an explicit second order $q$-difference
operator $D_B$ which diagonalizes the big $q$-Jacobi polynomials. 
The operator $D_B$ and the corresponding eigenvalue equations can be  
obtained from \eqref{limpolB} by taking the limit $k\rightarrow \infty$
in the equations
\begin{equation}\label{limBformal}
\bigl(\epsilon_k(cd/q)^{\frac{1}{2}}\bigr)^{|\lambda|}
\bigl((D-E_{\lambda})P_{\lambda}\bigr)(\epsilon_k^{-1}(q/cd)^{\frac{1}{2}}z;
\underline{t}_B(\epsilon_k);t)=0\quad (\lambda\in\Lambda)
\end{equation}
where $D$ is given by \eqref{secondorderqdiff} and $E_{\lambda}$ is given
by \eqref{AWb} (see \cite[section 3]{S1} for the easy computation). 
In \cite{SK1} it was shown that the formal computation of the
limits of $D$ and $E_{\lambda}$ in \cite[section 3]{S1} can be used to
prove the limit transition \eqref{limpolB} for generic $t\in (0,1)$. See
also \cite{S2} for the special case that $t=q^k$, $k\in {\mathbb{N}}$. 

It follows from Theorem \ref{thmB} that $D_B$ is symmetric with respect
to $\langle .,. \rangle_B$. In \cite{S1} the symmetry of $D_B$ was established
by direct calculations in which the asymptotic behaviour of the
weight function $\Delta^B$ (see Lemma \ref{motivation}) plays a crucial role.

The quadratic norm evaluations of the big $q$-Jacobi polynomials for the
special case $a=b=0$, $c=1$ and $t=q^k$ with $k\in {\mathbb{N}}$  
are recently proved by Baker and Forrester \cite[Section 4.3]{BF} 
using Pieri formulas. In order to see that 
the quadratic norms of 
the big $q$-Jacobi polynomials in this special case are in
agreement with the quadratic norm evaluations \cite[(4.3)]{BF}, one needs
to use the evaluation 
formula for the Macdonald polynomials \cite[(6.11)]{MM} together
with \cite[Proposition 3.2]{Ka} 
and the computation preceding \cite[Remark 5.4]{S1}.

The constant term identity for the big $q$-Jacobi polynomials can be rewritten
as follows.
\begin{Cor}\label{gevolgBB}
Let $t\in (0,1)$ and $(a,b,c,d)\in V_B$. We have
\begin{equation}\label{SelbergB}
\begin{split}
\langle 1,1\rangle_{B,t}^{a,b,c,d}=&\prod_{j=1}^n\left(
\frac{
\Gamma_q(\alpha+1+(j-1)\tau)\Gamma_q(\beta+1+(j-1)\tau)\Gamma_q(j\tau)}
{\Gamma_q(\alpha+\beta+2+(n+j-2)\tau)\Gamma_q(\tau)}\right.\\
&\left.\quad\qquad\quad
.\bigl(-q^{\alpha+1+(j-1)\tau}d/c,
-q^{\beta+1+(j-1)\tau)}c/d;q\bigr)_{\infty}^{-1}
\right)\\
\end{split}
\end{equation}
where $a=q^{\alpha}, b=q^{\beta}$ and $t=q^{\tau}$.
\end{Cor}
The $q$-Selberg integral \eqref{SelbergB} for $t=q^k$, $k\in {\mathbb{N}}$
reduces to the following evaluation formula.
\begin{Cor}\label{AE}
Let $t=q^k$ with $k\in {\mathbb{N}}$ and $(a,b,c,d)\in V_B$. We have
\begin{equation}
\begin{split}
\int_{z_1=-d}^c..\int_{z_n=-d}^c
\prod_{1\leq i<j\leq n}z_i^{2k}\bigl(q^{1-k}\frac{z_j}{z_i};q\bigr)_{2k}
\prod_{i=1}^n\frac{\bigl(qz_i/c,-qz_i/d;q\bigr)_{\infty}}
{\bigl(q^{1+\alpha}z_i/c,-q^{1+\beta}z_i/d;q\bigr)_{\infty}}d_qz_i&\nonumber\\
=q^{k^2\binom{n}{3}-\binom{k}{2}\binom{n}{2}}
\prod_{i=1}^n\left(
\frac{\Gamma_q(\alpha+1+(i-1)k)\Gamma_q(\beta+1+(i-1)k)\Gamma_q(ik+1)}
{\Gamma_q(\alpha+\beta+2+(n+i-2)k)\Gamma_q(k+1)}\qquad\right.&\nonumber\\
\left..\frac{\bigl(-d/c,-c/d;q\bigr)_{\infty}(cd)^{1+(i-1)k}}
{\bigl(-q^{\alpha+1+(i-1)k}d/c,-q^{\beta+1+(i-1)k}c/d;q\bigr)_{\infty}
(c+d)}\right)&\nonumber\\
\end{split}
\end{equation}
where $t=q^{\tau}$, $a=q^{\alpha}$ and $b=q^{\beta}$.
\end{Cor}
\begin{proof}
The bilinear form $\langle .,. \rangle_B$ differs from the one considered
in \cite[section 5]{S1} by the constant $c_B$, where $c_B$ for
$t=q^k$ ($k\in {\mathbb{N}}$) is given by \eqref{cjk}. Hence the corollary 
follows from Corollary \ref{gevolgBB} and the computation preceding
\cite[Remark 5.4]{S1}.
\end{proof} 
The constant term identity for the big $q$-Jacobi polynomials 
has appeared in the literature before.
Corollary \ref{AE} was conjectured by Askey \cite{A} and proved by Evans
\cite{E}. For arbitrary 
$t\in (0,1)$ the evaluation \eqref{SelbergB} is equivalent to
Tarasov's and Varchenko's 
summation formula \cite[Theorem (E.10)]{TV} (where the 
summation over $j$ in (E.10) should be taken from $0$ to $n$). 
The proof of Tarasov and Varchenko is by computing residues for an A type
generalization of Askey-Roy's $q$-beta integral. 
The equivalence of \cite[Theorem (E.10)]{TV} with \eqref{SelbergB} 
can be seen by making the substitution of variables
$p=q$, $x=t$, $a=-d$, $b=c$, $\alpha=-qb/d$, $\beta=qa/c$ and $l=n$
in \cite[(E.10)]{TV} and by applying the formula
\[\Delta^B(\sigma_Bq^{\nu'},
\rho_Bq^{\nu})=\Delta^B(\rho_Bq^{\nu},\sigma_Bq^{\nu'})
\prod_{\stackrel{1\leq k \leq j}{1\leq m \leq n-j}}\psi_t(-t^{m-k}d/c)
\]
where $\nu\in P(j)$ 
and $\nu'\in P(n-j)$ (here $\psi_t$ is given by \eqref{Psit}).

\section{Proof of Proposition \ref{link2}.}
The proof of Proposition \ref{link2} involves tedious calculations,
for which we will give the main steps in this section. We freely
use the notations of the previous sections.
We fix in this section $(a,b)\in V_L$ with
$b\not=0$. The condition $b\not=0$ is not essential; with slight modifications,
the proof goes also through when $b=0$.

Let $\epsilon\in {\mathbb{R}}_{>0}$ and set
$\rho_{L,j}(\epsilon):=t_{L,0}(\epsilon)t^{j-1}=
\epsilon^{-1}q^{\frac{1}{2}}t^{j-1}$
for $j\in {\mathbb{Z}}$ (here $\underline{t}_L(\epsilon)$ given by \eqref{tL}).
The parameter $t_{L,1}(\epsilon)=-aq^{\frac{1}{2}}$ has modulus $> 1$
for all $\epsilon$ if $a\in \bigl(q^{-\frac{1}{2}},q^{-1}\bigr)$, 
so can also contribute to the discrete parts of the symmetric form 
$\langle .,. \rangle_{\underline{t}_L(\epsilon),t}$
\eqref{symmvormAW} in the limit \eqref{limmeasureL}.
We therefore
write $\sigma_{L,j}:=t_{L,1}t^{j-1}=-aq^{\frac{1}{2}}t^{j-1}$ 
for $j\in {\mathbb{Z}}$.
Then $F(r;\underline{t}_L(\epsilon);t)$ \eqref{Fr} 
for $\epsilon>0$ sufficiently small is given by
\[F(r;\underline{t}_L(\epsilon);t)=
\bigcup_{l+m=r} D_{0}(l;\underline{t}_L(\epsilon);t)\times
 D_{1}(m;\underline{t}_L(\epsilon);t)\subset {\mathbb{C}}^r\]
with the set 
$D_0(l;\underline{t}_L(\epsilon);t)$ \eqref{Dir} for $l>0$ given by
\begin{equation}
\begin{split}
D_{0}(l;\underline{t}_L(\epsilon);t)=&\lbrace \rho_L(\epsilon)q^{\nu}
\, | \, \nu\in P_L(l;\epsilon)\rbrace,\\
P_L(l;\epsilon):=&\lbrace\nu\in P(l) \, | \quad 
|\rho_{L,0}(\epsilon)q^{\nu_l}|>1 \rbrace\\
\end{split}
\end{equation}
and with the set $D_1(m;\underline{t}_L(\epsilon);t)$ \eqref{Dir} for $m>0$
given by
\begin{equation}
D_{1}(m;\underline{t}_L(\epsilon);t)=
\begin{cases}
&\lbrace \sigma_L^m \rbrace \,\,\hbox{ if } |\sigma_{L,m}|>1,\\
&\emptyset\qquad \hbox{ otherwise}
\end{cases}
\end{equation}
where $\sigma_L^m:=(\sigma_{L,1},\ldots,\sigma_{L,m})$.
Note in particular, that $D_{1}(m;\underline{t}_L(\epsilon);t)$
is independent of $\epsilon$.
Using the explicit definition of the symmetric form $\langle .,. \rangle$
\eqref{symmvormAW} as given in section 4, as well as the definition
for $m_{\lambda}(z|u)$ \eqref{mlambdau}, we can write
\begin{equation}\label{p2}
\begin{split}
\left(\prod_{i=1}^n\bigl(-\epsilon^{-1}qt^{i-1},
-\epsilon^{-1}qat^{i-1};q\bigr)_{\infty}
\right)&\bigl(\epsilon q^{-1/2}\bigr)^{|\lambda |+|\mu |}
\langle m_{\lambda},m_{\mu}\rangle_{\underline{t}_L(\epsilon),t}\\
=\sum_{r,l,m,\nu}\,\,\iint\limits_{x\in T^{n-r}}\bigl(m_{\lambda}m_{\mu}\bigr)
\bigl(\rho_L q^{\nu},\epsilon
q^{-\frac{1}{2}}&\sigma_L^m,\epsilon
q^{-\frac{1}{2}}x|\epsilon q^{-\frac{1}{2}}\bigr)
{\mathcal{W}}^L_{l,m;r}(\nu,x;\epsilon)\frac{dx}{x}\\
\end{split}
\end{equation}
where the sum is over four tuples $(r,l,m,\nu)$ with
$r\in\lbrace 0,\ldots,n\rbrace$, $l,m\in {\mathbb{N}}_0$ with $l+m=r$,
and $\nu\in P(l)$ (the sum over $\nu\in P(l)$ should be ignored when $l=0$). 
The renormalized weights ${\mathcal{W}}^L_{l,m;r}(\nu,x;\epsilon)$ are given by
${\mathcal{W}}^L_{0,0;0}(-;x;\epsilon):=
\Delta(x;\underline{t}_L(\epsilon);t)$ for
$r=0$, and for $r=1,\ldots,n$,
\begin{equation}
\begin{split}
{\mathcal{W}}^L_{l,m;r}(\nu,x;\epsilon)=
&\prod_{i=1}^n
\bigl(-\epsilon^{-1}qt^{i-1},-\epsilon^{-1}qat^{i-1};q\bigr)_{\infty}
\nonumber\\
&.\frac{2^r\bigl(n-r+1\bigr)_r}{(2\pi i)^{n-r}}
\Delta_r^{AW}\bigl(\rho_L(\epsilon)q^{\nu},
\sigma_L^m,x;\underline{t}_L(\epsilon);t)
\nonumber\\
\end{split}
\end{equation}
if $\bigl(\rho_L(\epsilon)q^{\nu},\sigma_L^m\bigr)\in
F(r;\underline{t}_L(\epsilon);t)$ and zero otherwise.
We split the renormalized weights in three parts,
\begin{equation}\label{renweight}
{\mathcal{W}}_{l,m;r}^L(\nu,x;\epsilon)
=\frac{2^r\bigl(n-r+1\bigr)_r}{(2\pi i)^{n-r}}
\Delta_{1,l}^{AWL}(\nu;\epsilon)\Delta_{2,l,m}^{AWL}(\nu;\epsilon)
\Delta_{3,l,m}^{AWL}(\nu,x;\epsilon)
\end{equation}
where $\Delta_{1,l}^{AWL}$, $\Delta_{2,l,m}^{AWL}$ and $\Delta_{3,l,m}^{AWL}$ 
given by 
\[
\Delta_{1,l}^{AWL}(\nu;\epsilon):=
\left(\prod_{i=1}^{l}\bigl(-\epsilon^{-1}qt^{i-1},
-\epsilon^{-1}qat^{i-1};q\bigr)_{\infty}\right)
\Delta^{(d)}\bigl(\rho_L(\epsilon)q^{\nu};t_{L,0}(\epsilon)\bigr)
\]
if $\nu\in P_L(l;\epsilon)$ and zero otherwise,
\begin{equation}
\begin{split}
\Delta_{2,l,m}^{AWL}(\nu;\epsilon):=
\prod_{i=1}^{m}\bigl(-\epsilon^{-1}qt^{l+i-1},
&-\epsilon^{-1}qat^{l+i-1};q\bigr)_{\infty}\nonumber\\
&.\Delta^{(d)}\bigl(\sigma_L^m;t_{L,1}(\epsilon)\bigr)
\delta_c\bigl(\rho_L(\epsilon)q^{\nu};\sigma_L^m\bigr)
\nonumber\\
\end{split}
\end{equation}
if $\nu\in P_L(l;\epsilon)$, 
$D_{1}(m;\underline{t}_L(\epsilon);t)\not=\emptyset$ 
and zero otherwise,
\begin{equation}
\begin{split}
\Delta_{3,l,m}^{AWL}(\nu,x;\epsilon):=
\prod_{i=1}^{n-r}\bigl(-&\epsilon^{-1}qt^{r+i-1},
-\epsilon^{-1}qat^{r+i-1};q\bigr)_{\infty}\nonumber\\
&.\Delta(x;\underline{t}_L(\epsilon);t)
\delta_c\bigl(\rho_L(\epsilon)q^{\nu};x)\delta_c(\sigma_L^m;x)\nonumber\\
\end{split}
\end{equation}
if $\nu\in P_L(l;\epsilon), D_{1}(m;\underline{t}_L(\epsilon);t)\not=\emptyset$
and zero otherwise, 
where $r=l+m$ and $\Delta$ is given by \eqref{opdeling+},
$\Delta^{(d)}$ is given by \eqref{discreteweights}
and $\delta_c$ is given by \eqref{continuousinteraction}.
The formula 
$\delta_c(z;u,v)=\delta_c(z;u)\delta_c(z;v)$ 
is used for obtaining \eqref{renweight}. 
We have used 
for the definitions of $\Delta_{1,l}^{AWL}$, $\Delta_{2,l,m}^{AWL}$ and
$\Delta_{3,l,m}^{AWL}$ 
the obvious conventions  when $l=0$ or $m=0$; for instance,
\begin{equation}
\begin{split}
\Delta_{2,0,0}^{AWL}(-;\epsilon)&=1,\quad \Delta_{2,l,0}^{AWL}(\nu;\epsilon):=1
\nonumber\\
\Delta_{2,0,m}^{AWL}(-;\epsilon)&=\prod_{i=1}^m\bigl(-\epsilon^{-1}qt^{i-1},
-\epsilon^{-1}qat^{i-1};q\bigr)_{\infty}
\Delta^{(d)}\bigl(\sigma_L^m;t_{L,1}(\epsilon);t)\nonumber\\
\end{split}
\end{equation}
for $l,m,\nu$ such that $l,m>0$, $\nu\in P(l;\epsilon)$ and 
$D(\sigma_L^m;t_{L,1}(\epsilon);t)\not=0$.

We will use Lebesgue's Dominated Convergence Theorem to
pull a limit $\epsilon_k\downarrow 0$ in the right hand side of \eqref{p2} 
through the integration over $x\in T^{n-r}$ and through 
the infinite sum over 
$\nu\in P(l)$ for some 
sequence $\lbrace \epsilon_k\rbrace_{k\in {\mathbb{N}}_0}$ 
in ${\mathbb{R}}_{>0}$ converging to
$0$. Therefore, we
need  certain estimates for the functions 
$\Delta_{1,l}^{AWL}$, $\Delta_{2,l,m}^{AWL}$ and $\Delta_{3,l,m}^{AWL}$,
which we give in the following lemma.
\begin{Lem}\label{suffcheck}
Keep the notations and conventions as above. In particular, let $l,m\in
{\mathbb{N}}_0$ with $l+m\leq n$, and write $r:=l+m$.
Then there exists a sequence
$\lbrace \epsilon_k\rbrace_{k\in {\mathbb{N}}_0}$ in ${\mathbb{R}}_{>0}$
converging to $0$ such that\\
{\bf (i)} if $l\in {\mathbb{N}}$, then for all $\nu\in P(l)$,
\[\lim_{k\rightarrow \infty}\Delta_{1,l}^{AWL}(\nu;\epsilon_k)=
\bigl(q;q\bigr)_{\infty}^{-2l}\Delta^L(\rho_Lq^{\nu};a,b;t)\prod_{i=1}^{l}
\rho_{L,i}q^{\nu_i},\]
and there exists a $K\in {\mathbb{R}}_{>0}$ independent of $\nu\in P(l)$
such that
\begin{equation}\label{partL2}
\sup_{k\in {\mathbb{N}}_0}|\Delta_{1,l}^{AWL}(\nu;\epsilon_k)|\leq
K
\Delta^L\bigl(\rho_L q^{\nu};a,b;t)\prod_{i=1}^{l}\rho_{L,i}q^{\nu_i}
\end{equation}
for all $\nu\in P(l)$.\\
{\bf (ii)} if $m\in {\mathbb{N}}$, then
$\lim_{k\rightarrow\infty}\Delta_{2,l,m}^{AWL}(\nu;\epsilon_k)=0$ 
for all $\nu\in P(l)$
and 
\[\sup_{(\nu,k)\in
P(l)\times {\mathbb{N}}_0}|\Delta_{2,l,m}^{AWL}(\nu;\epsilon_k)|<\infty.
\]
{\bf (iii)}
if $r<n$, 
then $\lim_{k\rightarrow\infty}\Delta_{3,l,m}^{AWL}(\nu,x;\epsilon_k)=0$
for all $x\in T^{n-r}$, $\nu\in P(l)$ and 
\[\sup_{(\nu,x,k)\in P(l)\times T^{n-r} \times {\mathbb{N}}_0}
|\Delta_{3,l,m}^{AWL}(\nu,x;\epsilon_k)|<\infty.
\]
\end{Lem}
Before we prove Lemma \ref{suffcheck}, we complete the proof of Proposition
\ref{link2}. 
As we have already remarked in section 5, we have that the infinite sum
\begin{equation}\label{sumexpr}
(1-q)^{-n}\langle 1,1 \rangle_{L,t}^{a,b}=\sum_{\nu\in P(n)}
\Delta^L(\rho_Lq^{\nu};a,b;t)\prod_{i=1}^n\rho_{L,i}q^{\nu_i}
\end{equation}
is absolutely convergent (cf. \cite[proof of
Proposition 6.1]{S1}) and we have
\[
\sup_{\nu,\epsilon,x}|m_{\lambda}\bigl(\rho_Lq^{\nu},\epsilon
q^{-\frac{1}{2}}\sigma_L^m,
\epsilon q^{-\frac{1}{2}}x|\epsilon q^{-\frac{1}{2}}\bigr)|<\infty\quad 
(\lambda\in\Lambda),
\]
where the supremum is taken over triples $(\nu,\epsilon,x)$ with $\nu\in
P_L(l;\epsilon)$, $\epsilon\in {\mathbb{R}}_{>0}$ and $x\in T^{n-r}$.
By Lebesgue's 
Dominated Convergence Theorem, \eqref{limmon}, \eqref{p2}, \eqref{renweight}
and Lemma \ref{suffcheck} we thus obtain
\begin{equation}
\begin{split}
\lim_{k\rightarrow \infty}
\left(\prod_{i=1}^n\bigl(-\epsilon_k^{-1}qt^{i-1},
-\epsilon_k^{-1}qat^{i-1};q\bigr)_{\infty}\right)
\bigl(\epsilon_kq^{-\frac{1}{2}}\bigr)^{|\lambda |+|\mu |}
&\langle m_{\lambda},m_{\mu}\rangle_{\underline{t}_L(\epsilon_k),t}\nonumber\\
=\sum_{r,l,m,\nu}\,\,
\iint\limits_{x\in T^{n-r}}
\lim_{k\rightarrow \infty}\bigl(m_{\lambda}m_{\mu}\bigr)
\bigl(\rho_Lq^{\nu},\epsilon_k
q^{-\frac{1}{2}}\sigma_L^m,\epsilon_k
q^{-\frac{1}{2}}x&|\epsilon_k q^{-\frac{1}{2}}\bigr)
{\mathcal{W}}^L_{l,m;r}(\nu,x;\epsilon_k)\frac{dx}{x}\nonumber\\
=2^nn!\bigl(q;q\bigr)_{\infty}^{-2n}\sum_{\nu\in
P(n)}\bigl(\tilde{m}_{\lambda}\tilde{m}_{\mu}\bigr)
(\rho_Lq^{\nu})\Delta^L(\rho_Lq^{\nu}&;a,b;t)\prod_{i=1}^n\rho_{L,i}q^{\nu_i}
\nonumber\\
=2^nn!(1-q)^{-n}\bigl(q;q\bigr)_{\infty}^{-2n}\langle \tilde{m}_{\lambda},
\tilde{m}_{\mu} 
\rangle_{L,t}^{a,b}&\nonumber\\
\end{split}
\end{equation}
for some sequence $\lbrace \epsilon_k\rbrace_{k\in {\mathbb{N}}_0}$ in
${\mathbb{R}}_{>0}$ converging to $0$, 
where the 
sum in the second line is over four tuples $(r,l,m,\nu)$ with $r\in\lbrace
0,\ldots,n\rbrace$, 
$l,m\in {\mathbb{N}}_0$ with $l+m=r$, and $\nu\in P(l)$.
So for the proof of Proposition \ref{link2}, it suffices to prove
Lemma \ref{suffcheck}.
We use the following elementary lemma.
\begin{Lem}\label{standlim}\cite[Lemma 3.1]{SK2}
For given 
$\epsilon_0\in {\mathbb{R}}_{>0}$, we set $\epsilon_k:=\epsilon_0q^k$.\\
{\bf (a)} Let $c\in {\mathbb{C}}$. For $\epsilon_0\in {\mathbb{R}}_{>0}$
with $|c|\epsilon_0\not\in\lbrace q^{-l}\rbrace_{l\in {\mathbb{N}}_0}$
there exist positive constants 
$K^{\pm}>0$ which only depend on $\epsilon_0$ and
$|c|$, such that
$K^-\leq |\bigl(c\epsilon_k ;q\bigr)_{\infty}|\leq K^+$ for all
$k\in {\mathbb{N}}_0$. Furthermore, we have $\lim_{k\rightarrow
\infty}\bigl(c\epsilon_k;q\bigr)_{\infty}=1$. \\
{\bf (b)}
Let $a,b \in {\mathbb{C}}^*$, and set
\begin{equation}\label{fm}
f_{\lbrace l,m\rbrace}(\epsilon;a,b):=
\frac{\bigl(\epsilon^{-1}aq^{1-m};q\bigr)_m}
{\bigl(\epsilon^{-1}bq^{1-m-l};q\bigr)_m},\qquad (l,m\in {\mathbb{N}}_0).
\end{equation}
Let $\epsilon_0\in {\mathbb{R}}_{>0}$ such that 
$\epsilon_0^{-1}|b|\not\in\lbrace q^k\rbrace_{k\in
{\mathbb{N}}_0}$. 
Then there exists a positive constant $K>0$ which depends only on 
$\epsilon_0$, $|a|$ and $|b|$,
such that
$|f_{\lbrace l,m\rbrace}(\epsilon_k;a,b)|\leq K|q^la/b|^m$ 
for all $k,l,m\in {\mathbb{N}}_0$. Furthermore, we have 
$\lim_{k\rightarrow\infty}f_{\lbrace l,m\rbrace}(\epsilon_k;a,b)=(q^la/b)^m$.\\
{\bf (c)}
Let $u_i,v_j\in {\mathbb{C}}^*$ 
for $i\in\lbrace 1,\ldots,r\rbrace$,
$j\in\lbrace 1,\ldots,s\rbrace$ and assume that $r<s$, or that $r=s$ and
$|u_1\ldots u_r|<|v_1\ldots v_r|$.
Set
\begin{equation}\label{g}
g(\epsilon):=\frac{\bigl(\epsilon^{-1}u_1,\ldots,
\epsilon^{-1}u_r;q\bigr)_{\infty}}
{\bigl(\epsilon^{-1}v_1,\ldots,\epsilon^{-1}v_s;q\bigr)_{\infty}}.
\end{equation}
Let $\epsilon_0\in {\mathbb{R}}_{>0}$ 
such that $\epsilon_0^{-1}|v_j|\notin\lbrace q^l\rbrace_{l\in
{\mathbb{Z}}}$ for $j\in\lbrace 1,\ldots,s\rbrace$.
Then there exist a positive constant $K>0$ which depends only on $\epsilon_0$, 
$|u_i|$ and $|v_j|$, 
such that $\sup_{k\in {\mathbb{N}}_0}|g(\epsilon_k)|\leq K$.
Furthermore, we have $\lim_{k\rightarrow \infty}g(\epsilon_k)=0$.
\end{Lem}
The proofs of {\bf (b)} 
and {\bf (c)} are based on the formula \eqref{inversion}
for $q$-shifted factorials. See \cite[Lemma 3.1]{SK2} for details.

We proceed with the proof of Lemma \ref{suffcheck}. We use the notation
$\epsilon_k:=\epsilon_0q^k$ for given $\epsilon_0\in {\mathbb{R}}_{>0}$.\\
{\bf (i)}
By \eqref{discreteweights} one has 
\begin{equation}\label{fL}
\begin{split}
\Delta_{1,l}^{AWL}(\nu;\epsilon)&=\delta_d(\rho_L(\epsilon)q^{\nu})\\
&.\prod_{i=1}^l\Bigl(
\bigl(-\epsilon^{-1}qt^{i-1},
-\epsilon^{-1}qat^{i-1};q\bigr)_{\infty}w_d(\rho_{L,i}(\epsilon)q^{\nu_i};
\rho_{L,i}(\epsilon)q^{\nu_{i-1}})\Bigr)\\
\end{split}
\end{equation}
with $\delta_d$ given by \eqref{discreteinteraction} and $w_d$ given by
\eqref{weightfunctiondisc}.
By \eqref{discreteinteraction} and \eqref{tL}, we have
\begin{equation}\label{dd}
\delta_d(\rho_L(\epsilon)q^{\nu})=
F_1(\nu)G_1(\nu;\epsilon)
\end{equation}
with
\begin{equation}\label{FF1}
\begin{split}
F_1(\nu):=&\prod_{1\leq i<j\leq l}
\frac{\bigl(t^{j-i}q^{\nu_j-\nu_i};q\bigr)_{\tau}}
{\bigl(t^{i-j}q^{\nu_{i-1}-\nu_j};q\bigr)_{\nu_i-\nu_{i-1}}},\\
G_1(\nu;\epsilon):=&\prod_{1\leq i<j\leq l}
\frac{\bigl(\epsilon^2t^{2-i-j}q^{-\nu_i-\nu_j-1};q\bigr)_{\tau}}
{\bigl(\epsilon^{-2}t^{i+j-2}q^{\nu_{i-1}+\nu_j+1};q\bigr)_{\nu_i-
\nu_{i-1}}},\\
\end{split}
\end{equation}
for $\nu\in P(l)$, where $\nu_0=0$ and $t=q^{\tau}$.
By applying \eqref{inversion}
to the $q$-shifted factorials in the denominator of $F_1$
and using the formula
\begin{equation}\label{lover3}
\sum_{i=1}^l(i-1)(l-i)=\binom{l}{3},
\end{equation}
we obtain
\begin{equation}\label{expresL1}
\begin{split}
F_1(\nu)=&\delta_{qJ}(\rho_Lq^{\nu})q^{-2\tau^2\binom{l}{3}}\prod_{j=1}^l
\frac{\bigl(t^{j-1}q^{\nu_j+1};q\bigr)_{\infty}}
{\bigl(q^{\nu_j-\nu_{j-1}+1};q\bigr)_{\infty}}t^{-2(l-j)\nu_j}\\
&.\prod_{1\leq i<j\leq l}
\bigl(-q^{\nu_j-\nu_i+1}t^{j-i}\bigr)^{\nu_i-\nu_{i-1}}
q^{\binom{\nu_i-\nu_{i-1}}{2}},\\
\end{split}
\end{equation}
where $\delta_{qJ}$ \eqref{interactionL} 
is the interaction factor for the weight function
of the little $q$-Jacobi polynomials.
On the other hand, we have for $i=1,\ldots,l$ by \eqref{weightfunctiondisc},
\begin{equation}\label{dd1}
\bigl(-\epsilon^{-1}qt^{i-1},
-\epsilon^{-1}qat^{i-1};q\bigr)_{\infty}w_d(\rho_{L,i}(\epsilon)q^{\nu_i};
\rho_{L,i}(\epsilon)q^{\nu_{i-1}})
=I_{1,i}(\nu)J_{1,i}(\nu;\epsilon)
\end{equation}
with
\begin{eqnarray}
I_{1,i}(\nu)&:=&\frac{
\bigl(t^{i-1}q^{\nu_{i-1}+1}ab\bigr)^{\nu_{i-1}-\nu_i}}
{\bigl(q,qbt^{i-1}q^{\nu_i};q\bigr)_{\infty}
\bigl(q;q\bigr)_{\nu_i-\nu_{i-1}}}\nonumber\\
&=& v_L(\rho_{L,i}q^{\nu_i})\frac{a^{(1-i)\tau-\nu_i}
\bigl(t^{i-1}q^{\nu_{i-1}+1}ab\bigr)^{\nu_{i-1}-\nu_i}}
{\bigl(q,t^{i-1}q^{\nu_i+1};q\bigr)_{\infty}
\bigl(q;q\bigr)_{\nu_i-\nu_{i-1}}}\nonumber
\end{eqnarray}
(here $v_L$ \eqref{vL} 
is the one-variable weight function of the little $q$-Jacobi
polynomials)
and with $J_{1,i}(\nu;\epsilon)$ given by
\begin{equation}
\begin{split}
J_{1,i}(\nu;\epsilon):=&
\frac{\bigl(\epsilon^2t^{2-2i}q^{-2\nu_{i-1}-1};q\bigr)_{\infty}}
{\bigl(\epsilon^2bt^{1-i}q^{-\nu_{i-1}},-\epsilon t^{1-i}q^{-\nu_{i-1}},
-\epsilon at^{1-i}q^{-\nu_{i-1}};q\bigr)_{\infty}}\nonumber\\
&.\frac{\bigl(\epsilon^{-2}q^{1+2\nu_{i-1}}t^{2i-2};q\bigr)_{\nu_i-
\nu_{i-1}}\bigl(-\epsilon^{-1}t^{i-1}q;q\bigr)_{\nu_{i-1}}
\bigl(-\epsilon^{-1}at^{i-1}q;q\bigr)_{\nu_i}}
{\bigl(\epsilon^{-2}b^{-1}t^{i-1}q^{\nu_{i-1}+1},
-\epsilon^{-1}a^{-1}t^{i-1}q^{\nu_{i-1}+1};q\bigr)_{\nu_i-\nu_{i-1}}}
\nonumber\\
&.
\frac{\bigl(1-\epsilon^{-2}q^{1+2\nu_i}t^{2i-2}\bigr)}
{\bigl(1-\epsilon^{-2}q^{1+2\nu_{i-1}}t^{2i-2}\bigr)}.\nonumber\\
\end{split}
\end{equation}
So we have by \eqref{fL}, \eqref{dd}, \eqref{expresL1} and \eqref{dd1},
$\Delta_{1,l}^{AWL}(\nu;\epsilon)=M_1(\nu)N_1(\nu;\epsilon)$ 
for $\epsilon\in {\mathbb{R}}_{>0}$ and $\nu\in P_L(l;\epsilon)$
with
\begin{equation}\label{M1L}
\begin{split}
M_1(\nu):=&F_1(\nu)\prod_{i=1}^lI_{1,i}(\nu)\\
=&\Delta^L(\rho_Lq^{\nu})\bigl(q;q\bigr)_{\infty}^{-2l}t^{\binom{l}{2}}
\prod_{i=1}^l\bigl(t^{-2(l-i)}a^{-1}\bigr)^{\nu_i}
\bigl(t^{i-1}q^{\nu_{i-1}+1}ab
\bigr)^{\nu_{i-1}-\nu_i}\\
&.\prod_{1\leq i<j\leq l}
\bigl(-q^{\nu_j-\nu_i+1}t^{j-i}\bigr)^{\nu_i-\nu_{i-1}}
q^{\binom{\nu_i-\nu_{i-1}}{2}},\\
\end{split}
\end{equation}
($\Delta^L$ given by \eqref{DL})
and with
\[
N_1(\nu;\epsilon):=
G_1(\nu;\epsilon)\prod_{i=1}^lJ_{1,i}(\nu;\epsilon).\]
Now replace the factor 
$\bigl(-\epsilon^{-1}at^{i-1}q;q\bigr)_{\nu_i}$ in $J_{1,i}(\nu;\epsilon)$
by
\[\bigl(-\epsilon^{-1}at^{i-1}q^{\nu_{i-1}+1};q\bigr)_{\nu_i-\nu_{i-1}}
\bigl(-\epsilon^{-1}at^{i-1}q;q\bigr)_{\nu_{i-1}}\]
for $i\in\lbrace 1,\ldots,l\rbrace$, then $N_1(\nu;\epsilon)$ can explicitly
be given by
\begin{equation}\label{N1}
N_1(\nu;\epsilon)=N_1^1(\nu;\epsilon)N_1^2(\nu;\epsilon)N_1^3(\nu;\epsilon)
\end{equation}
with   
\begin{equation}\label{N11}
\begin{split}
N_1^1(\nu;\epsilon):=&\prod_{i=1}^l\frac{
\bigl(\epsilon^2t^{2-2i}q^{-2\nu_{i-1}-1};q\bigr)_{\infty}}
{\bigl(\epsilon^2bt^{1-i}q^{-\nu_{i-1}},-\epsilon t^{1-i}q^{-\nu_{i-1}},
-\epsilon at^{1-i}q^{-\nu_{i-1}};q\bigr)_{\infty}}\\
&.\prod_{1\leq i<j\leq l}
\bigl(\epsilon^2t^{2-i-j}q^{-\nu_i-\nu_j-1};q\bigr)_{\tau}\\
\end{split}
\end{equation}
if $\nu\in P_L(l;\epsilon)$ and zero otherwise,
\begin{equation}\label{N12}
\begin{split}
N_1^2(\nu;\epsilon):=
\prod_{i=1}^l\left(\frac{\bigl(\epsilon^{-2}q^{1+2\nu_{i-1}}t^{2i-2},
-\epsilon^{-1}at^{i-1}q^{1+\nu_{i-1}};q\bigr)_{\nu_i-\nu_{i-1}}}
{\bigl(\epsilon^{-2}b^{-1}t^{i-1}q^{1+\nu_{i-1}},
-\epsilon^{-1}a^{-1}t^{i-1}q^{1+\nu_{i-1}};q\bigr)_{\nu_i-\nu_{i-1}}}\right.&\\
\left..\frac{(\epsilon^{-2}q^{1+2\nu_i}t^{2i-2};q\bigr)_1}
{(\epsilon^{-2}q^{1+2\nu_{i-1}}t^{2i-2};q\bigr)_1}\right)&\\
\end{split}
\end{equation}
if $\nu\in P_L(l;\epsilon)$ and zero otherwise, and
\begin{equation}\label{N13}
N_1^3(\nu;\epsilon)=
\frac{\prod_{i=1}^l\bigl(-\epsilon^{-1}t^{i-1}q,-\epsilon^{-1}at^{i-1}q;
q\bigr)_{\nu_{i-1}}}
{\prod_{1\leq i<j\leq l}\bigl(\epsilon^{-2}t^{i+j-2}q^{\nu_{i-1}+\nu_j+1}
;q\bigr)_{\nu_i-\nu_{i-1}}}
\end{equation}
if $\nu\in P_L(l;\epsilon)$ and zero otherwise.
For the factor $N_1^1$ we have for generic $\epsilon_0>0$,
\begin{equation}\label{limN11}
\lim_{k\rightarrow\infty}N_1^1(\nu;\epsilon_k)=1 \qquad (\nu\in P(l))
\end{equation}
by Lemma \ref{standlim}{\bf (a)}.
For the factor $N_1^2$ we can use Lemma \ref{standlim}{\bf (b)} to
calculate the limit. We obtain for generic $\epsilon_0>0$, 
\begin{equation}\label{limN12}
\lim_{k\rightarrow\infty}N_1^2(\nu;\epsilon_k)=\prod_{i=1}^l
\bigl(q^{\nu_{i-1}+2}t^{i-1}a^2b\bigr)^{\nu_i-\nu_{i-1}}
\end{equation}
for all $\nu\in P(l)$.
As an example, let us calculate explicitly the limit of a factor of $N_1^2$,
using Lemma \ref{standlim}{\bf (b)}. Consider
the factor
\begin{equation}\label{exampleX}
N^{i,2}_1(\nu;\epsilon):=
\frac{\bigl(\epsilon^{-2}q^{1+2\nu_{i-1}}t^{2i-2};q\bigr)_{
\nu_i-\nu_{i-1}}}{\bigl(\epsilon^{-2}q^{1+\nu_{i-1}}b^{-1}t^{i-1};q
\bigr)_{\nu_i-\nu_{i-1}}}
\end{equation}
of $N_1^2(\nu;\epsilon)$ for some
$i\in\lbrace 1,\ldots,l\rbrace$. 
Then for generic $\epsilon_0>0$, 
we obtain by Lemma \ref{standlim}{\bf (b)},
\begin{equation}\label{example1}
\begin{split}
\lim_{k\rightarrow \infty}N^{i,2}_1(\nu;\epsilon_k)&=\lim_{k\rightarrow\infty}
N^{i,2}_1(\nu;\epsilon_{k+\nu_l})\\
&=\lim_{k\rightarrow\infty}N^{i,2}_1(\nu;q^{\nu_l}\epsilon_k)\\
&=\lim_{k\rightarrow\infty}f_{\lbrace \nu_{i-1},\nu_i-\nu_{i-1}\rbrace}\bigl(
\epsilon_{2\nu_l-\nu_i-\nu_{i-1}+2k};\epsilon_0^{-1}t^{2i-2},
\epsilon_0^{-1}b^{-1}t^{i-1}\bigr)\\
&=\bigl(q^{\nu_{i-1}}t^{i-1}b\bigr)^{\nu_i-\nu_{i-1}}.\\
\end{split}
\end{equation}
Similarly, one deals with the other factors of $N_1^2$ and one obtains
\eqref{limN12}.
Finally, we have for generic $\epsilon_0>0$, 
\begin{equation}\label{limN13}
\begin{split}
\lim_{k\rightarrow \infty}
N_1^3(\nu;\epsilon_k)=
&\prod_{i=1}^l\bigl(at^{2(i-1)}q^{\nu_{i-1}+1}\bigr)^{\nu_{i-1}}\\
&.\prod_{1\leq i<j\leq l}
\bigl(-t^{i+j-2}q^{\nu_{i-1}+\nu_j+1}\bigr)^{\nu_{i-1}-\nu_i}
q^{-\binom{\nu_i-\nu_{i-1}}{2}}\\
\end{split}
\end{equation}
since 
\begin{equation}\label{identBB}
\sum_{i=1}^l\nu_{i-1}=\sum_{1\leq i<j\leq
l}(\nu_i-\nu_{i-1})\qquad (\nu\in P(l)).
\end{equation}
We thus obtain for generic $\epsilon_0>0$ 
by \eqref{N1}, \eqref{limN11}, \eqref{limN12}
and \eqref{limN13},  
\begin{equation}\label{limN1}
\begin{split}
\lim_{k\rightarrow \infty}N_1(\nu;\epsilon_k)=
&\prod_{i=1}^lt^{(i-1)(\nu_i+\nu_{i-1})}q^{\nu_{i-1}\nu_i-\nu_{i-1}
+2\nu_i}a^{2\nu_i-\nu_{i-1}}b^{\nu_i-\nu_{i-1}}\\
&.\prod_{1\leq i<j\leq l}
\bigl(-t^{i+j-2}q^{\nu_{i-1}+\nu_j+1}\bigr)^{\nu_{i-1}-\nu_i}
q^{-\binom{\nu_i-\nu_{i-1}}{2}}\\
\end{split}
\end{equation}
for all $\nu\in P(l)$.
By \eqref{M1L} and \eqref{limN1} we obtain for generic $\epsilon_0>0$,
\begin{equation}
\begin{split}
\lim_{k\rightarrow\infty}\Delta_{1,l}^{AWL}(\nu;\epsilon_k)&=M_1(\nu)
\lim_{k\rightarrow\infty}N_1(\nu;\epsilon_k)\nonumber\\
&=\bigl(q;q\bigr)_{\infty}^{-2l}\Delta^L(\rho_Lq^{\nu})t^{\binom{l}{2}}
q^{|\nu|}\nonumber\\
&=\bigl(q;q\bigr)_{\infty}^{-2l}\Delta^L(\rho_Lq^{\nu})
\prod_{i=1}^l\rho_{L,i}q^{\nu_i}\nonumber\\
\end{split}
\end{equation}
for all $\nu\in P(l)$, where $\Delta^L$ is given by \eqref{DL} and
$|\nu|:=\nu_1+\cdots +\nu_l$ for $\nu\in P(l)$.

To prove the estimate \eqref{partL2}, we use the estimates of
Lemma \ref{standlim} {\bf (a)} and {\bf (b)} for (factors of) $N_1$.
For $N_1^1$, we use Lemma \ref{standlim}{\bf (a)} and the condition that
$N_1^1(\nu;\epsilon)=0$ if $\nu\not\in P_L(l;\epsilon)$ to prove that
$\sup_{\nu,k}|N_1^1(\nu;\epsilon_k)|<
\infty$ for generic $\epsilon_0>q^{\frac{1}{2}}$.
Indeed, since 
$\nu\in P_L(l;\epsilon)$ implies $\epsilon<q^{\frac{1}{2}+\nu_l}$,
we have for $\epsilon_0>q^{\frac{1}{2}}$,
\begin{equation}\label{estimateN11}
\sup_{(\nu,k)\in P(l)\times {\mathbb{N}}_0}|N_1^1(\nu;\epsilon_k)|=
\sup_{(\nu,k)\in P(l)\times 
{\mathbb{N}}_0}|N_1^1(\nu;\epsilon_kq^{\nu_l})|<\infty
\end{equation}
by \eqref{N11} and Lemma \ref{standlim}{\bf (a)}. 

For $N_1^2(\nu;\epsilon)$ 
we want to establish the estimate
\begin{equation}\label{estimateN12}
\sup_{k\in {\mathbb{N}}_0}|N_1^2(\nu;\epsilon_k)|\leq
K_1^2\prod_{i=1}^l\bigl(q^{\nu_{i-1}+2}t^{i-1}a^2|b|\bigr)^{\nu_i-\nu_{i-1}}
\end{equation}
for generic $\epsilon_0>q^{\frac{1}{2}}$ 
with $K_1^2>0$ independent of $\nu\in P(l)$,
in view of the limit \eqref{limN12}.
This can be done with the help of Lemma \ref{standlim}{\bf (b)}.
As an example, we consider the
factor $N_1^{i,2}(\nu;\epsilon)$ \eqref{exampleX}.
In view of the limit \eqref{example1}, we want to prove the estimate
\[\sup_{k\in {\mathbb{N}}_0}|N_1^{i,2}(\nu;\epsilon_k)|\leq K_1^{i,2}
\bigl(q^{\nu_{i-1}}t^{i-1}|b|\bigr)^{\nu_i-\nu_{i-1}}
\]
for generic $\epsilon_0>q^{\frac{1}{2}}$ with $K_1^{i,2}>0$ independent
of $\nu\in P(l)$.
This follows for generic $\epsilon_0>q^{\frac{1}{2}}$, using the fact that
$N_1^{i,2}(\nu;\epsilon)=0$ if $\nu\notin P_L(l;\epsilon)$, by the estimates
\begin{equation}
\begin{split}
\sup_{k\in {\mathbb{N}}_0}|N_{1}^{i,2}(\nu;\epsilon_k)|&=
\sup_{k\in {\mathbb{N}}_0}
|N_1^{i,2}(\nu;q^{\nu_l}\epsilon_{k})|\nonumber\\
&=\sup_{k\in {\mathbb{N}}_0}|f_{\lbrace \nu_{i-1},\nu_i-\nu_{i-1}
\rbrace}\bigl(
\epsilon_{2\nu_l-\nu_i-\nu_{i-1}+2k};\epsilon_0^{-1}t^{2i-2},
\epsilon_0^{-1}b^{-1}t^{i-1}\bigr)|\nonumber\\
&\leq\sup_{k\in {\mathbb{N}}_0}|f_{\lbrace \nu_{i-1},
\nu_i-\nu_{i-1}\rbrace}\bigl(
\epsilon_{k};\epsilon_0^{-1}t^{2i-2},
\epsilon_0^{-1}b^{-1}t^{i-1}\bigr)|\nonumber\\
&\leq
K_1^{i,2}\bigl(q^{\nu_{i-1}}t^{i-1}|b|\bigr)^{\nu_i-\nu_{i-1}}
\nonumber\\
\end{split}
\end{equation}
with $K_1^{i,2}$ independent of $\nu$ by Lemma \ref{standlim}{\bf (b)}. 
Similarly, one deals with the
other factors. 

For $N_1^3$, we want to prove that
\begin{equation}\label{estimateN13}
\begin{split}
\sup_{k\in {\mathbb{N}}_0}|N_1^3(\nu;\epsilon_k)|
\leq
&K_1^3\prod_{i=1}^l\bigl(at^{2(i-1)}q^{\nu_{i-1}+1}\bigr)^{\nu_{i-1}}\\
&.\prod_{1\leq i<j\leq l}
\bigl(t^{i+j-2}q^{\nu_{i-1}+\nu_j+1}\bigr)^{\nu_{i-1}-\nu_i}
q^{-\binom{\nu_i-\nu_{i-1}}{2}}\\
\end{split}
\end{equation}
for generic $\epsilon_0>q^{\frac{1}{2}}$ with $K_1^3>0$ independent of 
$\nu\in P(l)$, in view of the limit \eqref{limN13}. 
This follows by straightforward estimates, using  \eqref{inversion} and the
fact that $N_1^3(\nu;\epsilon)=0$ if $\nu\not\in P_L(l;\epsilon)$. 
Hence by \eqref{N1}, \eqref{estimateN11}, \eqref{estimateN12} and
\eqref{estimateN13} we have the estimate
\begin{equation}
\begin{split}
\sup_{k\in {\mathbb{N}}_0}|N_1(\nu;\epsilon_k)|\leq
&K_1\prod_{i=1}^lt^{(i-1)(\nu_i+\nu_{i-1})}q^{\nu_{i-1}\nu_i-\nu_{i-1}
+2\nu_i}a^{2\nu_i-\nu_{i-1}}|b|^{\nu_i-\nu_{i-1}}\\
&.\prod_{1\leq i<j\leq l}
\bigl(t^{i+j-2}q^{\nu_{i-1}+\nu_j+1}\bigr)^{\nu_{i-1}-\nu_i}
q^{-\binom{\nu_i-\nu_{i-1}}{2}}\\
\end{split}
\end{equation}
for generic $\epsilon_0>q^{\frac{1}{2}}$, 
with $K_1>0$ independent of $\nu\in P(l)$, so in particular,
\[
\sup_{k\in {\mathbb{N}}_0}|\Delta_{1,l}^{AWL}(\nu;\epsilon_k)|=
|M_1(\nu)|\sup_{k\in {\mathbb{N}}_0}|N_1(\nu;\epsilon_k)|
\leq K\Delta^L(\rho_Lq^{\nu})\prod_{i=1}^l\rho_{L,i}q^{\nu_i}
\]
with $K>0$ independent of $\nu\in P(l)$.
This completes the proof of  {\bf (i)}.\\
For {\bf (ii)}, note that 
$\delta_c(\rho_L(\epsilon);\sigma_L^m)$ ($\delta_c$ given
by \eqref{continuousinteraction})
can be rewritten as
\begin{equation}
\begin{split}
\delta_c(\rho_L(\epsilon)q^{\nu};\sigma_L^m)&=
\prod_{\stackrel{1\leq i\leq l}{1\leq j\leq m}}
\bigl(-\epsilon aq^{-\nu_i}t^{j-i},-\epsilon a^{-1}q^{-1-\nu_i}t^{2-i-j};
q\bigr)_{\tau}\\
&.\prod_{j=1}^m\frac{\bigl(-\epsilon^{-1}aqt^{j-1},
-\epsilon^{-1}a^{-1}t^{1-j};q
\bigr)_{\infty}}
{\bigl(-\epsilon^{-1}aqt^{l+j-1},-\epsilon^{-1}a^{-1}t^{l-j+1};q
\bigr)_{\infty}}\\
&.\prod_{\stackrel{0\leq i\leq l-1}{1\leq j\leq m}}
\frac{\bigl(-\epsilon^{-1}aq^{1+\nu_i}t^{l+j-1},
-\epsilon^{-1}a^{-1}q^{\nu_i}t^{l-j+1};q\bigr)_{\nu_{i+1}-\nu_i}}
{\bigl(-\epsilon^{-1}aq^{1+\nu_i}t^{i+j-1},
-\epsilon^{-1}a^{-1}q^{\nu_i}t^{i-j+1};q\bigr)_{\nu_{i+1}-\nu_{i}}}\\
\end{split}
\end{equation}
where $\nu_0=0$ and $t=q^{\tau}$.
By the explicit expression for $\Delta^{(d)}$
\eqref{discreteweights} and for the parameters $\underline{t}_L(\epsilon)$
\eqref{tL} we then obtain
\begin{equation}
\Delta_{2,l,m}^{AWL}(\nu;\epsilon)=
N_2^1N_2^2(\nu;\epsilon)N_2^3(\epsilon)N_2^4(\nu;\epsilon)
\end{equation} with
\[N_2^1:=\prod_{i=1}^m\frac{\bigl(a^{-2}q^{-1}t^{2(1-i)};q\bigr)_{\infty}}
{\bigl(q,aqt^{i-1},a^{-1}t^{1-i};q\bigr)_{\infty}}
\prod_{1\leq i<j\leq m}\bigl(t^{j-i},a^{-2}t^{2-i-j}q^{-1};q\bigr)_{\tau},\]
\[N_2^2(\nu;\epsilon):=
\prod_{j=1}^m\frac{\prod_{i=1}^l\bigl(-\epsilon aq^{-\nu_i}t^{j-i},
-\epsilon a^{-1}q^{-1-\nu_i}t^{2-i-j};q\bigr)_{\tau}}
{\bigl(-\epsilon a^{-1}bt^{1-j},-\epsilon abqt^{j-1};q\bigr)_{\infty}}\]
if $\nu\in P_L(l;\epsilon)$ and zero otherwise,
\[N_2^3(\epsilon):=\prod_{j=1}^m
\frac{\bigl(-\epsilon^{-1}t^{l+j-1}q;q\bigr)_{\infty}}
{\bigl(-\epsilon^{-1}t^{l-j+1}a^{-1};q\bigr)_{\infty}},
\]
\[N_2^4(\nu;\epsilon):=
\prod_{\stackrel{0\leq i\leq l-1}{1\leq j\leq m}}
\frac{\bigl(-\epsilon^{-1}aq^{1+\nu_i}t^{l+j-1},
-\epsilon^{-1}a^{-1}q^{\nu_i}t^{l-j+1};q\bigr)_{\nu_{i+1}-\nu_i}}
{\bigl(-\epsilon^{-1}aq^{1+\nu_i}t^{i+j-1},-\epsilon^{-1}
a^{-1}q^{\nu_i}t^{i-j+1};q
\bigr)_{\nu_{i+1}-\nu_i}}
\]
if $\nu\in P_L(l;\epsilon)$ and zero otherwise.
We have for generic $\epsilon_0>q^{\frac{1}{2}}$,
\[\lim_{k\rightarrow\infty}N_2^2(\nu;\epsilon_k)=1,
\quad \sup_{(\nu,k)\in P(l)\times {\mathbb{N}}_0}|N_2^2(\nu;\epsilon_k)|<
\infty\]
by Lemma \ref{standlim}{\bf (a)} and by the fact that $N_2^2(\nu;\epsilon)=0$
if $\nu\not\in P_L(l;\epsilon)$,
\[\lim_{k\rightarrow\infty}N_2^3(\epsilon_k)=0,\quad
\sup_{k\in {\mathbb{N}}_0}|N_2^3(\epsilon_k)|<\infty\]
by Lemma \ref{standlim}{\bf (c)} since $0<a<1/q$, and
\[\lim_{k\rightarrow\infty}N_2^4(\nu;\epsilon_k)=t^{2m|\nu|},\quad
\sup_{k\in {\mathbb{N}}_0}|N_2^4(\nu;\epsilon_k)|\leq 
K_2^4t^{2m|\nu|}\leq K_2^4\]
with $K_2^4>0$ independent of $\nu\in P(l)$ by Lemma \ref{standlim}{\bf (b)}
and by the fact that $N_2^4(\nu;\epsilon)=0$
if $\nu\not\in P_L(l;\epsilon)$.
This completes the proof of {\bf (ii)}.
Remains to prove {\bf (iii)}.
We can use the explicit
formulas for the weight function $\Delta$ \eqref{opdeling+},
\eqref{wc}, \eqref{deltatau} 
and for $\delta_c$ \eqref{continuousinteraction} as well
as the definition of $\underline{t}_L(\epsilon)$ \eqref{tL} to give an explicit
expression for $\Delta_{3,l,m}^{AWL}(\nu,x;\epsilon)$. 
This explicit expression can be written as
\[\Delta_{3,l,m}^{AWL}(\nu,x;\epsilon)=
N_3^1(x)N_3^2(\nu,x;\epsilon)N_3^3(x;\epsilon)
N_3^4(\nu,x;\epsilon)\]
with
\begin{equation}\label{F}
N_3^1(x):=\delta(x;t)\delta_c(\sigma_L^m;x)
\prod_{i=1}^{n-r}\frac{\bigl(x_i^2,x_i^{-2};q\bigr)_{\infty}}
{\bigl(-q^{\frac{1}{2}}x_i,-q^{\frac{1}{2}}x_i^{-1},-q^{\frac{1}{2}}ax_i,
-q^{\frac{1}{2}}ax_i^{-1};q\bigr)_{\infty}},
\end{equation}
if $D_1(m;\underline{t}_L(\epsilon);t)\not=0$ or $m=0$ and zero otherwise,
where $\delta$ is given by \eqref{deltatau},
\begin{equation}\label{G}
N_3^2(\nu,x;\epsilon):=
\prod_{j=1}^{n-r}\frac{\prod_{i=1}^l
\bigl(\epsilon q^{-\frac{1}{2}-\nu_i}t^{1-i}x_j,\epsilon q^{-\frac{1}{2}-\nu_i}
t^{1-i}x_j^{-1};q\bigr)_{\tau}}
{\bigl(\epsilon bq^{\frac{1}{2}}x_j,\epsilon bq^{\frac{1}{2}}x_j^{-1};
q\bigr)_{\infty}}\quad (t=q^{\tau})
\end{equation}
if $\nu\in P_L(l;\epsilon)$ and zero otherwise,
\begin{equation}\label{H}
N_3^3(x;\epsilon):=
\prod_{i=1}^{n-r}
\frac{\bigl(-\epsilon^{-1}qt^{r+i-1},-\epsilon^{-1}qat^{r+i-1},
\epsilon^{-1}q^{\frac{1}{2}}x_i,
\epsilon^{-1}q^{\frac{1}{2}}x_i^{-1};q\bigr)_{\infty}}
{\bigl(\epsilon^{-1}q^{\frac{1}{2}}x_i,\epsilon^{-1}q^{\frac{1}{2}}x_i^{-1},
\epsilon^{-1}q^{\frac{1}{2}}t^lx_i,\epsilon^{-1}q^{\frac{1}{2}}t^lx_i^{-1};
q\bigr)_{\infty}}
\end{equation}
and 
\begin{equation}\label{I}
N_3^4(\nu,x;\epsilon):=\prod_{\stackrel{0\leq i\leq l-1}{1\leq j\leq n-r}}
\frac{\bigl(\epsilon^{-1}q^{\frac{1}{2}+\nu_i}t^{l}x_j,
\epsilon^{-1}q^{\frac{1}{2}+\nu_i}t^{l}x_j^{-1};q\bigr)_{\nu_{i+1}-\nu_i}}
{\bigl(\epsilon^{-1}q^{\frac{1}{2}+\nu_i}t^{i}x_j,
\epsilon^{-1}q^{\frac{1}{2}+\nu_i}t^{i}x_j^{-1};q\bigr)_{\nu_{i+1}-\nu_{i}}}
\end{equation}
if $\nu\in P_L(l;\epsilon)$ and zero otherwise, where $\nu_0=0$.
We can proceed now as in the proof of {\bf (ii)}. 
Note that $N_3^1$ is bounded on $T^{n-r}$.
For $N_3^2$ it follows from Lemma \ref{standlim}
{\bf (a)} that $\lim_{k\rightarrow\infty}N_3^2(\nu,x;\epsilon_k)=1$ for
all $\nu\in P(l)$, $x\in T^{n-r}$ and that
\[
\sup_{(\nu,x,k)\in P(l)\times T^{n-r}\times
{\mathbb{N}}_0}|N_3^2(\nu,x;\epsilon_k)|
<\infty
\]
for  generic $\epsilon_0>q^{\frac{1}{2}}$.
By Lemma \ref{standlim}{\bf (c)} and the fact that $0<a<1/q$, we have
for generic 
$\epsilon_0>0$ that $\lim_{k\rightarrow\infty}N_3^3(x;\epsilon_k)=0$
for all $x\in T^{n-r}$ and 
\[
\sup_{(x,k)\in T^{n-r}\times
{\mathbb{N}}_0}|N_3^3(x;\epsilon_k)|<\infty.
\]
Finally, we can use Lemma \ref{standlim}{\bf (b)} to
prove that $\lim_{k\rightarrow\infty}N_3^4(\nu,x;\epsilon_k)=t^{2(n-r)|\nu|}$
for all $\nu\in P(l)$, $x\in T^{n-r}$ and that
\[
\sup_{(\nu,x,k)\in P(l)\times T^{n-r}\times
{\mathbb{N}}_0}|N_3^4(\nu,x;\epsilon_k)|
<\infty
\]
for generic $\epsilon_0>q^{\frac{1}{2}}$. 
This completes the proof of {\bf (iii)}.

\section{Proof of Proposition \ref{link3}.}
We first give a new expression for the weight $\underline{c}_B$ 
appearing in the definition of $\langle .,. \rangle_B$.
\begin{Lem}\label{newexpcBj}
The weight ${\underline{c}}_B\in \bigl({\mathbb{C}}^*\bigr)^{n+1}$ can be
rewritten as
\begin{equation}\label{cjtau}
\begin{split}
c_{B,j}=&\bigl(q;q\bigr)_{\infty}^n\prod_{i=1}^j\frac{\theta(-t^{i+j-n}c/d)}
{\theta(-t^{1-i}d/c)\theta(-t^ic/d)}
\prod_{i=1}^{n-j}\frac{1}{\theta(-t^{1-i}c/d)}\\
&.
q^{-2\tau^2\bigl((n-j)\binom{j}{2}+\binom{j}{3}+\binom{n-j}{3}\bigr)}
t^{-\binom{j}{2}-\binom{n-j}{2}}\\
&.c^{-2\tau\bigl(j(n-j)+\binom{j}{2}\bigr)-j}d^{-2\tau\binom{n-j}{2}+j-n}\\
\end{split}
\end{equation}
for $j=0,1,\ldots,n$.
\end{Lem}
\begin{proof}
For $j=0$, \eqref{cjtau} 
follows from \eqref{cB} since $c_{B,j}=c_Bd_{B,j}$ with
$d_{B,0}=1$. For $j\in\lbrace 1,\ldots, n\rbrace$, write 
$\tilde{c}_{B,j}$ for the right hand side of \eqref{cjtau}, then by
the explicit expression \eqref{dBj}
for $d_{B,j}$ it remains to prove that 
\[\frac{\tilde{c}_{B,j}}{\tilde{c}_{B,j-1}}=
\prod_{m=j+1}^n\Psi_t(-t^{n-m-j+1}d/c)\]
for $j\in\lbrace 1,\ldots,n\rbrace$, with $\Psi_t$ given by \eqref{Psit}.
This follows by a direct calculation.
\end{proof}

The proof of Proposition \ref{link3} is similar to the proof 
of Proposition \ref{link2}. 
Again, we indicate the main steps in the computation.
We freely use the notations of previous sections.
We fix in this section $(a,b,c,d)\in V_B$ with $a,b\not=0$.
The conditions $a,b\not=0$ are not essential; with slight modifications,
the proof goes also through when $a=0$ or $b=0$.

For $\epsilon\in {\mathbb{R}}_{>0}$ we set
$\rho_{B,j}(\epsilon):=t_{B,0}(\epsilon)t^{j-1}=(qc/d)^{\frac{1}{2}}
\epsilon^{-1}t^{j-1}$ 
and $\sigma_{B,j}(\epsilon):=t_{B,1}(\epsilon)t^{j-1}=-(qd/c)^{\frac{1}{2}}
\epsilon^{-1}t^{j-1}$ for
$j\in {\mathbb{Z}}$, where $\underline{t}_B(\epsilon)$ 
is given by \eqref{tB}.
Then we have for $\epsilon>0$ sufficiently small,
\[F(r;\underline{t}_B(\epsilon);t)=
\bigcup_{\stackrel{l+m=r}{l,m\in {\mathbb{N}}_0}} 
D_{0}(l;\underline{t}_B(\epsilon);t)\times
 D_{1}(m;\underline{t}_B(\epsilon);t)\subset {\mathbb{C}}^r\]
where $F(r)$ is given by \eqref{Fr} and
\begin{equation}
\begin{split}
D_{0}(l;\underline{t}_B(\epsilon);t)&=\lbrace \rho_B(\epsilon)q^{\nu}
\, | \, \nu\in P_B^{(0)}(l;\epsilon)\rbrace,\\
P_B^{(0)}(l;\epsilon)&:=\lbrace\nu\in P(l) \, | \, 
|\rho_{B,l}(\epsilon)q^{\nu_l}|>1 \rbrace,\\
\end{split}
\end{equation}
if $l>0$, respectively 
\begin{equation}
\begin{split}
D_{1}(m;\underline{t}_B(\epsilon);t)&=\lbrace \sigma_B(\epsilon)q^{\nu}
\, | \, \nu\in P_B^{(1)}(m;\epsilon)\rbrace,\\ 
P_B^{(1)}(m;\epsilon)&:=\lbrace\nu\in P(m) \, | \, 
|\sigma_{B,m}(\epsilon)q^{\nu_m}|>1 \rbrace\\
\end{split}
\end{equation}
if $m>0$. We write
\begin{equation}\label{p3}
\begin{split}
&\left(\prod_{i=1}^n\bigl(-\epsilon^{-2}qt^{i-1};q\bigr)_{\infty}
\right)\bigl((cd/q)^{\frac{1}{2}}\epsilon\bigr)^{|\lambda |+|\mu |}
\langle m_{\lambda},m_{\mu}\rangle_{\underline{t}_B(\epsilon),t}\\
&=\sum_{r,l,m,\nu,\nu'}\,\,\iint\limits_{x\in T^{n-r}}
\bigl(m_{\lambda}m_{\mu}\bigr)
\bigl(\rho_Bq^{\nu},\sigma_Bq^{\nu'}, (cd/q)^{\frac{1}{2}}\epsilon x
|(cd/q)^{\frac{1}{2}}\epsilon\bigr)
{\mathcal{W}}^B_{l,m;r}(\nu,\nu',x;\epsilon)\frac{dx}{x}\\
\end{split}
\end{equation}
where $\rho_{B,i}:=ct^{i-1}$, $\sigma_{B,i}:=-dt^{i-1}$ 
and $m_{\lambda}(z|u)$ is given by \eqref{mlambdau},
and the sum
is over five tuples $(r,l,m,\nu,\nu')$ with $r\in\lbrace 0,\ldots,n\rbrace$,
$l,m\in {\mathbb{N}}_0$ with $l+m=r$, $\nu\in P(l)$, $\nu'\in P(m)$, and with
the renormalized weight ${\mathcal{W}}^B_{l,m;r}(\nu,\nu',x;\epsilon)$ given by
\begin{equation}
{\mathcal{W}}^B_{l,m;r}(\nu,\nu',x;\epsilon):=
\frac{2^r\bigl(n-r+1\bigr)_r}{(2\pi
i)^{n-r}}\Delta_{1,l,m}^{AWB}(\nu,\nu';\epsilon)
\Delta_{2,l,m}^{AWB}(\nu,\nu',x;\epsilon)
\end{equation}
when $r=l+m$, with
\begin{equation}
\begin{split}
\Delta_{1,l,m}^{AWB}(\nu,\nu';\epsilon):=
&\left(\prod_{i=1}^{r}\bigl(-\epsilon^{-2}qt^{i-1};q\bigr)_{\infty}\right)
\delta_c\bigl(\rho_B(\epsilon)q^{\nu};\sigma_B(\epsilon)q^{\nu'}\bigr)\\
&.\Delta^{(d)}\bigl(\rho_B(\epsilon)q^{\nu};t_{B,0}(\epsilon)\bigr)
\Delta^{(d)}\bigl(\sigma_B(\epsilon)q^{\nu'};t_{B,1}(\epsilon)\bigr)\\
\end{split}
\end{equation}
if $\nu\in P_B^{(0)}(l;\epsilon)$, 
$\nu'\in P_B^{(1)}(m;\epsilon)$ and zero otherwise,
\begin{equation}
\begin{split}
\Delta_{2,l,m}^{AWB}(\nu,\nu',x;\epsilon):=&\prod_{i=1}^{n-r}
\bigl(-\epsilon^{-2}qt^{r+i-1};q\bigr)_{\infty}\\
&.\Delta(x;\underline{t}_B(\epsilon);t)
\delta_c\bigl(\rho_B(\epsilon)q^{\nu};x)
\delta_c(\sigma_B\bigl(\epsilon)q^{\nu'};
x\bigr)\\
\end{split}
\end{equation}
if $\nu\in P_B^{(0)}(l;\epsilon)$, $\nu'\in P_B^{(1)}(m;\epsilon)$ 
and zero otherwise, with $\delta_c$ given by \eqref{continuousinteraction}.
We use the obvious conventions when $l=0$, $m=0$ or $r=n$
(compare with the little $q$-Jacobi case in section 8).
In particular, we have $\Delta_{2,l,n-l}^{AWB}(\nu,\nu';\epsilon)=1$ for 
$\nu\in P_B^{(0)}(l;\epsilon)$, $\nu'\in P_B^{(1)}(n-l;\epsilon)$
and $l\in\lbrace 0,\ldots,n\rbrace$.

We will use Lebesgue's Dominated Convergence Theorem to
pull a limit $\epsilon_k\downarrow 0$ in the right hand side of \eqref{p3} 
through the integration over $x\in T^{n-r}$ and through 
the infinite sums over 
$\nu\in P(l)$ and $\nu'\in P(m)$
for some sequence $\lbrace \epsilon_k\rbrace_{k\in {\mathbb{N}}_0}$ 
in ${\mathbb{R}}_{>0}$ converging to
$0$.
We use the following
lemma.
\begin{Lem}\label{suffcheckB}
Keep the notations and conventions as above. Let $l,m\in {\mathbb{N}}_0$
with $l+m\in\lbrace 0,\ldots,n\rbrace$ and write $r:=l+m$. 
Then there exists a sequence
$\lbrace \epsilon_k\rbrace_{k\in {\mathbb{N}}_0}$ in ${\mathbb{R}}_{>0}$
converging to $0$ such that\\
{\bf (i)} for all $\nu\in P(l)$, $\nu'\in P(m)$ we have
\[
\lim_{k\rightarrow \infty}\Delta_{1,l,m}^{AWB}(\nu,\nu';\epsilon_k)
=\bigl(q;q\bigr)_{\infty}^{-2r}c_{B,l}\Delta^B(\rho_Bq^{\nu}, 
\sigma_Bq^{\nu'})\prod_{i=1}^{l}
\rho_{B,i}q^{\nu_i}\prod_{j=1}^m|\sigma_{B,j}|q^{\nu_j'},\]
and there exists a $K\in {\mathbb{R}}_{>0}$ independent of $\nu\in P(l)$
and $\nu'\in P(m)$ such that
\[\sup_{k\in {\mathbb{N}}_0}|\Delta_{1,l,m}^{AWB}(\nu,\nu';\epsilon_k)|
\leq Kc_{B,l}\Delta^B(\rho_Bq^{\nu}, \sigma_Bq^{\nu'})
\prod_{i=1}^{l}\rho_{B,i}q^{\nu_i}\prod_{j=1}^m|\sigma_{B,j}|q^{\nu_j'}.
\]
for all $\nu\in P(l)$ and all $\nu'\in P(m)$, where
$\Delta^B(z)=\Delta^B(z;a,b,c,d;t)$ is given by \eqref{DB}.\\
{\bf (ii)} if $r<n$, then
$\lim_{k\rightarrow\infty}\Delta_{2,l,m}^{AWB}(\nu,\nu',x;\epsilon_k)=0$ 
for all $\nu\in P(l)$, $\nu'\in P(m)$, $x\in T^{n-r}$
and
\[
\sup_{(k,\nu,\nu',x)\in {\mathbb{N}}_0\times
P(l)\times P(m)\times
T^{n-r}}|\Delta_{2,l,m}^{AWB}(\nu,\nu',x;\epsilon_k)|<\infty.\]
\end{Lem}
Before we prove Lemma \ref{suffcheckB}, we complete the proof of Proposition
\ref{link3}. As we have already remarked before, we have that the infinite sum
\begin{equation}\label{sumexprB}
(1-q)^{-n}\langle 1,1 \rangle_B=\sum_{(\nu,\nu',l)}
c_{B,l}\Delta^B(\rho_Bq^{\nu},\sigma_Bq^{\nu'})
\prod_{i=1}^l\rho_{B,i}q^{\nu_i}\prod_{j=1}^{n-l}
|\sigma_{B,j}|q^{\nu_j'}
\end{equation}
is absolutely convergent, where
the sum is taken over the three tuples $(\nu,\nu',l)$ with
$\nu\in P(l)$ and $\nu'\in P(n-l)$ and $l\in\lbrace
0,\ldots,n\rbrace$. 
Since
\[\sup_{\nu,\nu',x,\epsilon}
|m_{\lambda}\bigl(\rho_Bq^{\nu},
\sigma_Bq^{\nu'},(cd/q)^{\frac{1}{2}}\epsilon x|
(cd/q)^{\frac{1}{2}}\epsilon\bigr)|<\infty,
\]
where the supremum is taken over the four tuples $(\nu,\nu',x,\epsilon)$
with $\nu\in P_B^{(0)}(l;\epsilon)$, 
$\nu'\in P_B^{(1)}(m;\epsilon)$, $x\in T^{n-r}$ 
($r=l+m$) and $\epsilon>0$, we obtain by Lebesgue's 
Dominated Convergence Theorem, \eqref{limmon}, \eqref{p3}, \eqref{sumexprB}
and Lemma \ref{suffcheckB} 
\begin{equation}
\begin{split}
&\lim_{k\rightarrow\infty}
\left(\prod_{i=1}^n\bigl(-\epsilon_k^{-2}qt^{i-1};q\bigr)_{\infty}
\right)\bigl((cd/q)^{\frac{1}{2}}\epsilon_k\bigr)^{|\lambda |+|\mu |}
\langle m_{\lambda},m_{\mu}\rangle_{\underline{t}_B(\epsilon_k),t}\nonumber\\
&=\frac{2^nn!}{\bigl(q;q\bigr)_{\infty}^{2n}}\sum_{l=0}^n\sum_{
\stackrel{\nu\in P(l)}{\nu'\in P(n-l)}}
c_{B,l}\bigl(\tilde{m}_{\lambda}\tilde{m}_{\mu}\Delta^B\bigr)
(\rho_Bq^{\nu}, \sigma_Bq^{\nu'})\prod_{i=1}^l\rho_{B,i}q^{\nu_i}
\prod_{j=1}^{n-l}|\sigma_{B,j}|q^{\nu_j'}\nonumber\\
&=2^nn!(1-q)^{-n}\bigl(q;q\bigr)_{\infty}^{-2n}\langle \tilde{m}_{\lambda},
\tilde{m}_{\mu} \rangle_{B,t}^{a,b,c,d},\nonumber\\
\end{split}
\end{equation}
for some sequence $\lbrace \epsilon_k\rbrace_{k\in {\mathbb{N}}_0}$ in
${\mathbb{R}}_{>0}$ converging to $0$.
So for the proof of Proposition \ref{link3}, it suffices to prove
Lemma \ref{suffcheckB}.\\
{\it Proof of Lemma } \ref{suffcheckB}.\\
Using the explicit expressions for
$\Delta^{(d)}$ 
\eqref{discreteweights}, $\delta_c$ \eqref{continuousinteraction}
and $\underline{t}_B(\epsilon)$ \eqref{tB}, we can write
\begin{equation}\label{opsplitsingB}
\Delta_{1,l,m}^{AWB}(\nu,\nu';\epsilon)
:=U_0(\nu,\nu';l,m)U_+(\epsilon;\nu,\nu';l,m)U_-(\epsilon;\nu,\nu';l,m)
\end{equation}
with $U_0,U_+$ respectively 
$U_-$ the factor of $\Delta_{1,l,m}^{AW}$ consisting
of products of $q$-shifted factorials of the form $(e;q\bigr)_{s}$,
$\bigl(\epsilon^2e;q\bigr)_{s}$ respectively $\bigl(\epsilon^{-2}e;q\bigr)_s$
($s\in {\mathbb{N}}_0\cup\lbrace\infty\rbrace$).
By a straightforward computation, the factors $U_0,U_+$ and $U_-$ can
be explicitly given by
\begin{equation}\label{U0}
\begin{split}
U_0(\nu,\nu';l,m):=&\Psi_0(\nu,l;a,b,c,d)\Psi_0(\nu',m;b,a,d,c)\\
&.\prod_{\stackrel{1\leq i\leq l}{1\leq j\leq m}}
\bigl(-t^{i-j}q^{\nu_i-\nu_j'}c/d,-t^{j-i}q^{\nu_j'-\nu_i}d/c;q\bigr)_{\tau}\\
\end{split}
\end{equation}
if $(\nu,\nu')\in P_0(l;\epsilon)\times P_1(m;\epsilon)$, and zero otherwise, 
with $t=q^{\tau}$ and with
\begin{equation}\label{Psi0}
\begin{split}
\Psi_0(\nu,l;a,&b,c,d):=F_1(\nu)\\
&.\prod_{i=1}^l\frac{1}
{\bigl(q,-t^{1-i}q^{-\nu_{i-1}}d/c,at^{i-1}q^{\nu_{i-1}+1},
-bt^{i-1}q^{\nu_{i-1}+1}c/d;q\bigr)_{\infty}}\\
&.\prod_{i=1}^l\frac{\bigl(at^{i-1}q^{\nu_{i-1}+1},
-bt^{i-1}q^{\nu_{i-1}+1}c/d;q
\bigr)_{\nu_i-\nu_{i-1}}}
{\bigl(q,-t^{i-1}q^{\nu_{i-1}+1}c/d;q\bigr)_{\nu_i-\nu_{i-1}}
\bigl(abt^{i-1}q^{\nu_{i-1}+1}\bigr)^{\nu_i-\nu_{i-1}}}
\end{split}
\end{equation}
where $\nu_0=0$ and $F_1(\nu)$ is given by \eqref{FF1}, and 
\begin{equation}\label{U+}
\begin{split}
U_+(\epsilon;\nu,\nu';l,m):=&
\Psi_+(\epsilon;\nu,l;a,b,c,d)\Psi_+(\epsilon;\nu',m;b,a,d,c)\\
&.\prod_{\stackrel{1\leq i\leq l}{1\leq j\leq m}}
\bigl(-\epsilon^2t^{2-i-j}q^{-1-\nu_i-\nu_j'};q\bigr)_{\tau}\\
\end{split}
\end{equation}
if $(\nu,\nu')\in P_B^{(0)}(l;\epsilon)\times P_B^{(1)}(m;\epsilon)$ 
and zero otherwise, 
with
\begin{equation}\label{Psi+}
\begin{split}
\Psi_+(\epsilon;\nu,l;a,b,c,d)&:=
\prod_{i=1}^l
\frac{\bigl(\epsilon^2t^{2(1-i)}q^{-2\nu_{i-1}-1}d/c;q\bigr)_{\infty}}
{\bigl(\epsilon^2at^{1-i}q^{-\nu_{i-1}}d/c,
-\epsilon^2bt^{1-i}q^{-\nu_{i-1}};q\bigr)_{\infty}}\\
&.\prod_{1\leq i<j\leq l}\bigl(\epsilon^2t^{2-i-j}q^{-\nu_i-\nu_j-1}d/c;q
\bigr)_{\tau}\\
\end{split}
\end{equation}
and 
\begin{equation}\label{U-}
\begin{split}
U_-(\epsilon;\nu,&\nu';l,m):=\Psi_-(\epsilon;\nu,l;a,b,c,d)
\Psi_-(\epsilon;\nu',m;b,a,d,c)\\
&.\left(\prod_{j=1}^m\frac{\bigl(-\epsilon^{-2}qt^{l+j-1};q
\bigr)_{\infty}}
{\bigl(-\epsilon^{-2}qt^{j-1};q\bigr)_{\infty}}\right)
\prod_{\stackrel{1\leq i\leq l}{1\leq j\leq m}}
\bigl(-\epsilon^{-2}t^{i+j-2}q^{\nu_i+\nu_j'+1};q\bigr)_{\tau}\\
\end{split}
\end{equation}
if $(\nu,\nu')\in 
P_B^{(0)}(l;\epsilon)\times P_B^{(1)}(m;\epsilon)$ and zero otherwise, 
with
\begin{equation}\label{Psi-}
\begin{split}
\Psi_-(\epsilon;\nu,l;a,&b,c,d):=\prod_{i=1}^l\frac{\bigl(
\epsilon^{-2}t^{2(i-1)}q^{2\nu_i+1}c/d;q\bigr)_1}
{\bigl( \epsilon^{-2}t^{2(i-1)}q^{2\nu_{i-1}+1}c/d;q\bigr)_1}\\
&.\prod_{i=1}^l\frac{\bigl(\epsilon^{-2}t^{2(i-1)}q^{2\nu_{i-1}+1}c/d,
-\epsilon^{-2}t^{i-1}q^{\nu_{i-1}+1};q\bigr)_{\nu_i-\nu_{i-1}}}
{\bigl(\epsilon^{-2}a^{-1}t^{i-1}q^{\nu_{i-1}+1}c/d,
-\epsilon^{-2}b^{-1}t^{i-1}q^{\nu_{i-1}+1};q\bigr)_{\nu_i-\nu_{i-1}}}\\
&.\frac{\prod_{i=1}^l\bigl(-\epsilon^{-2}qt^{i-1};q\bigr)_{\nu_{i-1}}}
{\prod_{1\leq i<j\leq l}\bigl(\epsilon^{-2}t^{i+j-2}
q^{\nu_{i-1}+\nu_j+1}c/d;q\bigr)_{\nu_i-\nu_{i-1}}}.\\
\end{split}
\end{equation}
For given $\epsilon_0\in {\mathbb{R}}^*$, we write $\epsilon_k:=\epsilon_0q^k$.
Then for generic $\epsilon_0>0$ we have
\begin{equation}\label{limU+}
\lim_{k\rightarrow\infty}
U_+(\epsilon_k;\nu,\nu';l,m)=1
\end{equation}
for all $(\nu,\nu')\in P(l)\times P(m)$ by Lemma \ref{standlim}{\bf (a)}.
By \eqref{identBB}, we have for generic $\epsilon_0>0$
\begin{equation}\label{factor3line}
\begin{split}
&\lim_{k\rightarrow\infty}
\frac{\prod_{i=1}^l\bigl(-\epsilon_k^{-2}qt^{i-1};q\bigr)_{\nu_{i-1}}}
{\prod_{1\leq i<j\leq l}\bigl(\epsilon_k^{-2}t^{i+j-2}
q^{\nu_{i-1}+\nu_j+1}c/d;q\bigr)_{\nu_i-\nu_{i-1}}}\\
&=\prod_{i=1}^lq^{\binom{\nu_{i-1}+1}{2}}t^{(i-1)\nu_{i-1}}
\prod_{1\leq i<j\leq l}
\bigl(-t^{i+j-2}q^{\nu_{i-1}+\nu_j+1}c/d\bigr)^{\nu_{i-1}-\nu_{i}}
q^{-\binom{\nu_i-\nu_{i-1}}{2}}\\
\end{split}
\end{equation}
for the factor of $\Psi_-$ in the third line of \eqref{Psi-}.
The factor of $U_-$ in the second line of \eqref{U-} can be rewritten as
\begin{equation}\label{factor2line}
\begin{split}
\left(\prod_{j=1}^m\frac{\bigl(-\epsilon^{-2}qt^{l+j-1};q\bigr)_{\infty}}
{\bigl(-\epsilon^{-2}qt^{j-1};q\bigr)_{\infty}}\right)
\prod_{\stackrel{1\leq i\leq l}{1\leq j\leq m}}
&\bigl(-\epsilon^{-2}t^{i+j-2}q^{\nu_i+\nu_j'+1};q\bigr)_{\tau}\\
&=\prod_{\stackrel{1\leq i\leq l}{1\leq j\leq m}}
\frac{\bigl(-\epsilon^{-2}qt^{i+j-1};q\bigl)_{\nu_i+\nu_j'}}
{\bigl(-\epsilon^{-2}qt^{i+j-2};q\bigr)_{\nu_i+\nu_j'}}.\\
\end{split}
\end{equation}
It follows then from \eqref{U-}, 
\eqref{Psi-}, \eqref{factor3line}, \eqref{factor2line}
and Lemma \ref{standlim}{\bf (b)} that for generic $\epsilon_0>0$,
\begin{equation}\label{limU-}
\lim_{k\rightarrow\infty} U_-(\epsilon_k;\nu,\nu';l,m)=t^{m|\nu|+l|\nu'|}
\Psi_-^{\infty}(\nu,l;a,b,c,d)\Psi_-^{\infty}(\nu',m;b,a,d,c)
\end{equation}
with
\begin{equation}\label{Psi-infty}
\begin{split}
\Psi_-^{\infty}(\nu,l;a,b,c,d):=&\prod_{i=1}^l\bigl(t^{i-1}q^{\nu_{i-1}+2}ab
\bigr)^{\nu_i-\nu_{i-1}}q^{\binom{\nu_{i-1}+1}{2}}t^{(i-1)\nu_{i-1}}\\
&.\prod_{1\leq i<j\leq l}\bigl(-t^{i+j-2}q^{\nu_{i-1}+\nu_j+1}c/d
\bigr)^{\nu_{i-1}-\nu_i}q^{-\binom{\nu_i-\nu_{i-1}}{2}}.\\
\end{split}
\end{equation}
We will rewrite now $U_0$ in the form
\begin{equation}\label{U0rewritten}
\begin{split}
U_0(\nu,\nu';&l,m)=\bigl(q;q\bigr)_{\infty}^{-r}
\prod_{i=1}^l\frac{\theta(-t^{i-m}c/d)}
{\theta(-t^{1-i}d/c)\theta(-t^ic/d)}
\prod_{j=1}^m\frac{1}{\theta(-t^{1-j}c/d)}\\
&.C_0(\nu,\nu';l,m)\Delta^B(\rho_Bq^{\nu},\sigma_Bq^{\nu'};a,b,c,d;t)\\
\end{split}
\end{equation}
and we determine the factor $C_0(\nu,\nu';l,m)$ explicitly.
Using \eqref{theta} and the formula $\theta(x)=\theta(qx^{-1})$
for the Jacobi theta function $\theta(x)$ 
\eqref{thetafunction}, we can rewrite the factor 
$\Psi_0(\nu,l;a,b,c,d)$ \eqref{Psi0} as
\begin{equation}\label{Psi0rewritten}
\begin{split}
\Psi_0(\nu,&l;a,b,c,d)=F_1(\nu)\\
&.\prod_{i=1}^l
\frac{v_B(\rho_{B,i}q^{\nu_i};a,b,c,d)
(t^{i-1}c/d)^{\nu_{i-1}}q^{\binom{\nu_{i-1}+1}{2}}}
{\theta(-t^{1-i}d/c)
\bigl(q^{\nu_i+1}t^{i-1};q\bigr)_{\infty}\bigl(q;q\bigr)_{\nu_i-\nu_{i-1}}
(abt^{i-1}q^{\nu_{i-1}+1})^{\nu_i-\nu_{i-1}}}\\
\end{split}
\end{equation}
where $v_B$ \eqref{vB} is the one-variable 
weight function for the big $q$-Jacobi
polynomials.
Since $v_B(-dx;a,b,c,d)=v_B(dx;b,a,d,c)$, we obtain from \eqref{Psi0rewritten}
\begin{equation}\label{Psi0rewritten2}
\begin{split}
\Psi_0&(\nu',m;b,a,d,c)=F_1(\nu')\\
&.\prod_{j=1}^m
\frac{v_B(\sigma_{B,j}q^{\nu_j'};a,b,c,d)
(t^{j-1}d/c)^{\nu_{j-1}'}q^{\binom{\nu_{j-1}'+1}{2}}}
{\theta(-t^{1-j}c/d)
\bigl(q^{\nu_j'+1}t^{j-1};q\bigr)_{\infty}\bigl(q;q\bigr)_{\nu_j'-\nu_{j-1}'}
(abt^{j-1}q^{\nu_{j-1}'+1})^{\nu_j'-\nu_{j-1}'}}.\\
\end{split}
\end{equation}
The factor $F_1(\nu)$ \eqref{FF1} 
of $\Psi_0(\nu,l;a,b,c,d)$ 
can be rewritten as
\begin{equation}\label{deelPsi0rewritten}
\begin{split}
F_1(\nu)
=&\delta_{qJ}(\rho_Bq^{\nu})q^{-2\tau^2\binom{l}{3}}c^{-l(l-1)\tau}
\prod_{i=1}^l\frac{\bigl(t^{i-1}q^{1+\nu_i};
q\bigr)_{\infty}}
{\bigl( q^{1+\nu_i-\nu_{i-1}};q\bigr)_{\infty}}t^{-2(l-i)\nu_i}\\
&.\prod_{1\leq i<j\leq l}
\bigl(-q^{\nu_j-\nu_i+1}t^{j-i}\bigr)^{\nu_i-\nu_{i-1}}
q^{\binom{\nu_i-\nu_{i-1}}{2}}\\
\end{split}
\end{equation}
for $\nu\in P(l)$ where 
$\delta_{qJ}$ is given by \eqref{interactionL} and $t=q^{\tau}$.
This follows from \eqref{expresL1} since 
$\delta_{qJ}(\rho_Bq^{\nu})=
c^{l(l-1)\tau}\delta_{qJ}(\rho_Lq^{\nu})$ for $\nu\in P(l)$
(here $\rho_{L,i}=t^{i-1}$).
Similarly, we have for the factor $F_1(\nu')$ of
$\Psi_0(\nu',m;b,a,d,c)$, 
\begin{equation}\label{deelPsi0rewritten2}
\begin{split}
F_1(\nu')=&\delta_{qJ}(\sigma_Bq^{\nu'})q^{-2\tau^2\binom{m}{3}}d^{-m(m-1)\tau}
\prod_{j=1}^m\frac{\bigl(t^{j-1}q^{1+\nu_j'};
q\bigr)_{\infty}}
{\bigl( q^{1+\nu_j'-\nu_{j-1}'};q\bigr)_{\infty}}t^{-2(m-j)\nu_j'}\\
&.\prod_{1\leq i<j\leq m}
\bigl(-q^{\nu_j'-\nu_i'+1}t^{j-i}\bigr)^{\nu_i'-\nu_{i-1}'}
q^{\binom{\nu_i'-\nu_{i-1}'}{2}},\\
\end{split}
\end{equation}
since $\delta_{qJ}(\sigma_Bq^{\nu'})=
d^{m(m-1)\tau}\delta_{qJ}(\rho_Lq^{\nu'})$ for $\nu'\in P(m)$.

Finally, we set for $z=(z_1,\ldots,z_r)$ with $r:=l+m$,
\[\tilde{\delta}_{qJ}^l(z):=\prod_{\stackrel{1\leq i\leq l}{l+1\leq j\leq r}}
|z_i-z_j||z_i|^{2\tau-1}\bigl(qt^{-1}z_j/z_i;q\bigr)_{2\tau-1},
\]
then the factor of $U_0$ in the second line of \eqref{U0} can be rewritten as
\begin{equation}\label{deelU0rewritten}
\begin{split}
\prod_{\stackrel{1\leq i\leq l}{1\leq j\leq m}}
\bigl(-t^{i-j}&q^{\nu_i-\nu_j'}c/d,-t^{j-i}q^{\nu_j'-\nu_i}d/c;q\bigr)_{\tau}\\
&=\tilde{\delta}_{qJ}^l(\rho_Bq^{\nu},\sigma_Bq^{\nu'})t^{m|\nu|-l|\nu'|}
\prod_{i=1}^l
\frac{\theta(-t^{i-m}c/d)}{\theta(-t^ic/d)
\bigl(ct^{i-1}q^{\nu_i}\bigr)^{2m\tau}}\\
\end{split}
\end{equation}
for $\nu\in P(l)$ and $\nu'\in P(m)$, since we have for $i\in\lbrace
1,\ldots,l\rbrace$, $j\in\lbrace 1,\ldots,m\rbrace$ that
\begin{equation}\label{hulpje}
\begin{split}
\bigl(-t^{i-j}&q^{\nu_i-\nu_j'}c/d,-t^{j-i}q^{\nu_j'-\nu_i}d/c;q\bigr)_{\tau}\\
&=\frac{\theta(-t^{i-j}c/d)
\bigl(-t^{j-i-1}q^{1+\nu_j'-\nu_i}d/c;q\bigr)_{\infty}}
{\theta(-t^{i-j+1}c/d)\bigl(-t^{j-i+1}q^{\nu_j'-\nu_i}d/c;q\bigr)_{\infty}}
\bigl(1+t^{j-i}q^{\nu_j'-\nu_i}d/c\bigr)t^{\nu_i-\nu_j'}\\
\end{split}
\end{equation}
(formula 
\eqref{hulpje} follows from a straightforward computation using
\eqref{inversion}, \eqref{theta} and $\theta(x)=\theta(qx^{-1})$).
Since we have
\[
\delta_{qJ}(\rho_Bq^{\nu},\sigma_Bq^{\nu'})=\delta_{qJ}(\rho_Bq^{\nu})
\delta_{qJ}(\sigma_Bq^{\nu'})
\tilde{\delta}_{qJ}^l(\rho_Bq^{\nu},\sigma_Bq^{\nu'})
\]
for $\nu\in P(l)$ and $\nu'\in P(m)$, we obtain from
\eqref{DB}, \eqref{U0},
\eqref{Psi0rewritten}, \eqref{Psi0rewritten2}, \eqref{deelPsi0rewritten}, 
\eqref{deelPsi0rewritten2}
and \eqref{deelU0rewritten} that \eqref{U0rewritten} holds with
\begin{equation}\label{C0}
\begin{split}
C_0(\nu,\nu';l,m)=&t^{-m|\nu|-l|\nu'|}c^{-2lm\tau}q^{-2m\binom{l}{2}\tau^2}\\
&.{\hat{C}}_0(\nu,l;a,b,c,d){\hat{C}}_0(\nu',m;b,a,d,c)\\
\end{split}
\end{equation}
where
\begin{equation}\label{hatC0}
\begin{split}
{\hat{C}}_0(\nu,l;a,b,c,d)&:=q^{-2\tau^2\binom{l}{3}}c^{-2\binom{l}{2}\tau}\\
&.\prod_{i=1}^l\bigl(t^{i-1}c/d\bigr)^{\nu_{i-1}}q^{\binom{\nu_{i-1}+1}{2}}
t^{-2(l-i)\nu_i}\bigl(abt^{i-1}q^{\nu_{i-1}+1}\bigr)^{\nu_{i-1}-\nu_i}\\
&.\prod_{1\leq i<j\leq l}
\bigl(-q^{\nu_j-\nu_i+1}t^{j-i}\bigr)^{\nu_i-\nu_{i-1}}
q^{\binom{\nu_i-\nu_{i-1}}{2}}.\\
\end{split}
\end{equation}
We have by Lemma \ref{newexpcBj}
(with $n$ in the right hand side 
of \eqref{cjtau} equal to $r=l+m$ in this situation)
by \eqref{opsplitsingB}, 
\eqref{limU+}, \eqref{limU-}, \eqref{U0rewritten} and \eqref{C0} 
that for generic $\epsilon_0>0$,
\begin{equation}\label{almostthere}
\begin{split}
\lim_{k\rightarrow \infty}\Delta_{1,l,m}^{AWB}(\nu,\nu';&\epsilon_k)=
U_0(\nu,\nu';l,m)\lim_{k\rightarrow\infty}U_-(\epsilon_k;\nu,\nu';l,m)\\
=&\bigl(q;q\bigr)_{\infty}^{-2r}c_{B,l}\Delta^B(\rho_Bq^{\nu}, \sigma_B
q^{\nu'};a,b,c,d;t)\\
&.E(\nu,l;a,b,c,d)E(\nu',m;b,a,d,c)\\
\end{split}
\end{equation}
for $\nu\in P(l)$ and $\nu'\in P(m)$, with
\begin{equation}
E(\nu,l;a,b,c,d):=q^{2\tau^2\binom{l}{3}}t^{\binom{l}{2}}
c^{2\binom{l}{2}\tau+l}
{\hat{C}}_0(\nu,l;a,b,c,d)\Psi_-^{\infty}(\nu,l;a,b,c,d).
\end{equation}
By \eqref{Psi-infty} and \eqref{hatC0}, we obtain
\[E(\nu,l;a,b,c,d)=
c^lt^{\binom{l}{2}}q^{|\nu|}=\prod_{i=1}^l\rho_{B,i}q^{\nu_i},
\quad \nu\in P(l).\]
In particular we have  
$E(\nu',m;b,a,d,c)=\prod_{j=1}^m|\sigma_{B,j}|q^{\nu_j'}$ for $\nu'\in P(m)$, 
hence by \eqref{almostthere}
\[
\lim_{k\rightarrow \infty}\Delta_{1,l,m}^{AWB}(\nu,\nu';\epsilon_k)
=\bigl(q;q\bigr)_{\infty}^{-2r}c_{B,l}\Delta^B(\rho_Bq^{\nu}, 
\sigma_Bq^{\nu'})\prod_{i=1}^{l}
\rho_{B,i}q^{\nu_i}\prod_{j=1}^m|\sigma_{B,j}|q^{\nu_j'}\]
for all $\nu\in P(l)$, $\nu'\in P(m)$. To complete the proof of Lemma
\ref{suffcheckB}{\bf (i)}, it suffices to prove that
for generic $\epsilon_0>\hbox{max}\bigl((qc/d)^{\frac{1}{2}},
(qd/c)^{\frac{1}{2}}\bigr)$,
\[
\sup_{k\in {\mathbb{N}}_0}|\Delta_{1,l,m}^{AWB}(\nu,\nu';\epsilon_k)|
\leq Kc_{B,l}\Delta^B(\rho_Bq^{\nu}, \sigma_Bq^{\nu'};a,b,c,d;t)\prod_{i=1}^{l}
\rho_{B,i}q^{\nu_i}\prod_{i=1}^m|\sigma_{B,j}|q^{\nu_j'}
\]
for all $\nu\in P(l)$ and all 
$\nu'\in P(m)$, with $K>0$ independent of $\nu$ and
$\nu'$. This can be proved by similar arguments as in the little
$q$-Jacobi case (see proof of Lemma \ref{suffcheck}{\bf (i)}).
In particular, the estimates for almost all factors of $\Delta_{1,l,m}^{AWB}$
can be obtained from one of the three estimates of Lemma \ref{standlim}.
Only for the factor in the third line of the expression of $\Psi_-$ 
\eqref{Psi-} one needs a seperate argument to establish the desired estimate.
We may assume that this factor is zero unless
$\nu\in P_B^{(0)}(l;\epsilon)$ and $\nu'\in P_B^{(1)}(m;\epsilon)$.
Since we have the limit \eqref{factor3line}, we would like to establish
the estimate
\begin{equation}\label{factor3lineestimate}
\begin{split}
&\sup_{\lbrace k\in {\mathbb{N}}_0| \nu\in P_B^{(0)}(l;\epsilon_k) \hbox{ and }
\nu'\in P_B^{(1)}(m;\epsilon_k)\rbrace}\left|
\frac{\prod_{i=1}^l\bigl(-\epsilon_k^{-2}qt^{i-1};q\bigr)_{\nu_{i-1}}}
{\prod_{1\leq i<j\leq l}\bigl(\epsilon_k^{-2}t^{i+j-2}
q^{\nu_j+\nu_{i-1}+1}c/d;q\bigr)_{\nu_i-\nu_{i-1}}}\right|\nonumber\\
&\leq K'\prod_{i=1}^lq^{\binom{\nu_{i-1}+1}{2}}t^{(i-1)\nu_{i-1}}
\prod_{1\leq i<j\leq l}
\bigl(t^{i+j-2}q^{\nu_{i-1}+\nu_j+1}c/d\bigr)^{\nu_{i-1}-\nu_{i}}
q^{-\binom{\nu_i-\nu_{i-1}}{2}}
\nonumber\\
\end{split}
\end{equation}
for some $K'>0$ independent of $\nu\in P(l)$ and $\nu'\in P(m)$,
for generic $\epsilon_0>\hbox{max}\bigl((qc/d)^{\frac{1}{2}},
(qd/c)^{\frac{1}{2}}\bigr)$.
This can be done similarly as we have done for the factor $N_1^3$ \eqref{N13}
in the little $q$-Jacobi case.\\
The proof of Lemma \ref{suffcheckB}{\bf (ii)} is similar to the proof of 
Lemma \ref{suffcheck}{\bf (iii)}.\\

{\it Acknowledgements.} It is a pleasure to thank Prof. Masatoshi Noumi and 
Dr. Mathijs S. Dijkhuizen for their hospitality during my stay at Kobe
university. The author thanks Prof. Tom H. Koornwinder 
for many valuable comments.

\bibliographystyle{amsplain}
 
\end{document}